\documentclass{article}

\usepackage[a4paper]{geometry}
\usepackage{xcolor}

\usepackage{fullpage} 
\usepackage{authblk}
\usepackage{amsmath}
\usepackage{amsthm}
\usepackage{amssymb}
\usepackage{mathrsfs}
\usepackage{listings}
\usepackage{mathtools}
\usepackage{caption}
\usepackage{subfigure}
\usepackage{graphicx}
\usepackage{hyperref}
\usepackage{tikz-cd}
\usepackage{standalone}
\usepackage{float}
\usetikzlibrary{patterns}

\usepackage[capitalise]{cleveref}

\hypersetup{linktocpage}
\usepackage [english]{babel}
\usetikzlibrary{babel}
\usepackage [autostyle, english = american]{csquotes}
\MakeOuterQuote{"}

\newtheorem{theorem}{Theorem}
\newtheorem{proposition}{Proposition}
\newtheorem{definition}{Definition}
\newtheorem{remark}{Remark}
\newtheorem{lemma}{Lemma}
\newtheorem{corollary}{Corollary}

\newcommand{\R}{\mathbb{R}}
\newcommand{\N}{\mathbb{N}}
\newcommand{\Z}{\mathbb{Z}}

\newcommand{\Id}{\mathrm{I}}

\newcommand{\curl}{\mathrm{curl}\,}
\renewcommand{\div}{\mathrm{div}\,}
\renewcommand{\d}{\mathrm{d}}
\newcommand{\BS}{\mathrm{BS}\,}
\renewcommand{\H}{\mathrm{H}\,}

\usepackage{tocbasic} 
\usepackage[multiuser,draft,inline,marginclue]{fixme}
\fxusetheme{colorsig}
\FXRegisterAuthor{rou}{arou}{Roussel}
\FXRegisterAuthor{rr}{arr}{Remi}

\usepackage[
    backend=biber,
    maxbibnames=4,
    maxcitenames=4,
    style=alphabetic,
    sorting=nyt,
    isbn=false,
    url=false,
    firstinits=true,
    date=year]{biblatex}
\makeatletter
\def\blx@bblfile@biber{%
  \blx@secinit
  \begingroup
  \blx@bblstart
\begingroup
\makeatletter
\@ifundefined{ver@biblatex.sty}
  {\@latex@error
     {Missing 'biblatex' package}
     {The bibliography requires the 'biblatex' package.}
      \aftergroup }
  {}
\endgroup

\refsection{0}
  \datalist[entry]{nyt/global//global/global}
    \entry{alnaesUnifiedFormLanguage2014}{article}{}
      \name{author}{5}{}{%
        {{hash=904acb8bf2a8298efc38ab0f9740bc36}{%
           family={Alnæs},
           familyi={A\bibinitperiod},
           given={Martin\bibnamedelima S.},
           giveni={M\bibinitperiod\bibinitdelim S\bibinitperiod}}}%
        {{hash=bae2791f700c5805db2f0fda7977a949}{%
           family={Logg},
           familyi={L\bibinitperiod},
           given={Anders},
           giveni={A\bibinitperiod}}}%
        {{hash=4b17ef0cc79d17415951d391575aa713}{%
           family={Ølgaard},
           familyi={Ø\bibinitperiod},
           given={Kristian\bibnamedelima B.},
           giveni={K\bibinitperiod\bibinitdelim B\bibinitperiod}}}%
        {{hash=fd9b078969e262e2b3456fa16701b5ce}{%
           family={Rognes},
           familyi={R\bibinitperiod},
           given={Marie\bibnamedelima E.},
           giveni={M\bibinitperiod\bibinitdelim E\bibinitperiod}}}%
        {{hash=997effee24e95b11152025bbad8e00ad}{%
           family={Wells},
           familyi={W\bibinitperiod},
           given={Garth\bibnamedelima N.},
           giveni={G\bibinitperiod\bibinitdelim N\bibinitperiod}}}%
      }
      \strng{namehash}{af7397fd05974c8336c566aff2e513f7}
      \strng{fullhash}{0e81b43a4ae9667e216853ac4eda3a05}
      \strng{bibnamehash}{af7397fd05974c8336c566aff2e513f7}
      \strng{authorbibnamehash}{af7397fd05974c8336c566aff2e513f7}
      \strng{authornamehash}{af7397fd05974c8336c566aff2e513f7}
      \strng{authorfullhash}{0e81b43a4ae9667e216853ac4eda3a05}
      \field{labelalpha}{Aln+14}
      \field{sortinit}{A}
      \field{sortinithash}{a3dcedd53b04d1adfd5ac303ecd5e6fa}
      \field{labelnamesource}{author}
      \field{labeltitlesource}{shorttitle}
      \field{abstract}{We present the Unified Form Language (UFL), which is a domain-specific language for representing weak formulations of partial differential equations with a view to numerical approximation. Features of UFL include support for variational forms and functionals, automatic differentiation of forms and expressions, arbitrary function space hierarchies for multifield problems, general differential operators and flexible tensor algebra. With these features, UFL has been used to effortlessly express finite element methods for complex systems of partial differential equations in near-mathematical notation, resulting in compact, intuitive and readable programs. We present in this work the language and its construction. An implementation of UFL is freely available as an open-source software library. The library generates abstract syntax tree representations of variational problems, which are used by other software libraries to generate concrete low-level implementations. Some application examples are presented and libraries that support UFL are highlighted.}
      \field{day}{5}
      \field{issn}{0098-3500}
      \field{journaltitle}{ACM Transactions on Mathematical Software}
      \field{month}{3}
      \field{number}{2}
      \field{shortjournal}{ACM Trans. Math. Softw.}
      \field{shorttitle}{Unified form language}
      \field{title}{Unified form language: {{A}} domain-specific language for weak formulations of partial differential equations}
      \field{urlday}{24}
      \field{urlmonth}{10}
      \field{urlyear}{2023}
      \field{volume}{40}
      \field{year}{2014}
      \field{dateera}{ce}
      \field{urldateera}{ce}
      \field{pages}{9:1\bibrangedash 9:37}
      \range{pages}{-1}
      \verb{doi}
      \verb 10.1145/2566630
      \endverb
      \verb{file}
      \verb /home/rrobin/Zotero/storage/7UWP4W8Q/Alnæs et al. - 2014 - Unified form language A domain-specific language .pdf
      \endverb
      \verb{urlraw}
      \verb https://dl.acm.org/doi/10.1145/2566630
      \endverb
      \verb{url}
      \verb https://dl.acm.org/doi/10.1145/2566630
      \endverb
    \endentry
    \entry{alonso-rodriguezFiniteElementApproximation2018}{article}{}
      \name{author}{5}{}{%
        {{hash=86abaf634bb2d42d7a1bda9b96993027}{%
           family={Alonso-Rodríguez},
           familyi={A\bibinithyphendelim R\bibinitperiod},
           given={A.},
           giveni={A\bibinitperiod}}}%
        {{hash=546191dcb9a64962f0992e004d132f2a}{%
           family={Camaño},
           familyi={C\bibinitperiod},
           given={J.},
           giveni={J\bibinitperiod}}}%
        {{hash=3820c6bbfd84d613887b653aaba31086}{%
           family={Rodríguez},
           familyi={R\bibinitperiod},
           given={R.},
           giveni={R\bibinitperiod}}}%
        {{hash=31171371a77aefab1503a9ea99f99271}{%
           family={Valli},
           familyi={V\bibinitperiod},
           given={A.},
           giveni={A\bibinitperiod}}}%
        {{hash=cbeaa96391df8432c623ceb3d50885f2}{%
           family={Venegas},
           familyi={V\bibinitperiod},
           given={P.},
           giveni={P\bibinitperiod}}}%
      }
      \strng{namehash}{5ad08f30dd3d50a1c2bb2246e0332152}
      \strng{fullhash}{aa0b4b1521c78031d0a3d3a5171f9224}
      \strng{bibnamehash}{5ad08f30dd3d50a1c2bb2246e0332152}
      \strng{authorbibnamehash}{5ad08f30dd3d50a1c2bb2246e0332152}
      \strng{authornamehash}{5ad08f30dd3d50a1c2bb2246e0332152}
      \strng{authorfullhash}{aa0b4b1521c78031d0a3d3a5171f9224}
      \field{labelalpha}{Alo+18}
      \field{sortinit}{A}
      \field{sortinithash}{a3dcedd53b04d1adfd5ac303ecd5e6fa}
      \field{labelnamesource}{author}
      \field{labeltitlesource}{title}
      \field{abstract}{In this paper we are concerned with two topics: the formulation and analysis of the eigenvalue problem for the \$\$\textbackslash mathop \{\textbackslash mathbf \{curl\}\}\textbackslash nolimits \$\$operator in a multiply connected domain and its numerical approximation by means of finite elements. We prove that the \$\$\textbackslash mathop \{\textbackslash mathbf \{curl\}\}\textbackslash nolimits \$\$operator is self-adjoint on suitable Hilbert spaces, all of them being contained in the space for which \$\$\textbackslash mathop \{\textbackslash mathbf \{curl\}\}\textbackslash nolimits \textbackslash varvec\{v\}\textbackslash cdot \textbackslash varvec\{n\}=0\$\$on the boundary. Additional constraints must be imposed when the physical domain is not topologically trivial: we show that a viable choice is the vanishing of the line integrals of \$\$\textbackslash varvec\{v\}\$\$on suitable homological cycles lying on the boundary. A saddle-point variational formulation is devised and analyzed, and a finite element numerical scheme is proposed. It is proved that eigenvalues and eigenfunctions are efficiently approximated and some numerical results are presented in order to assess the performance of the method.}
      \field{day}{1}
      \field{issn}{1615-3383}
      \field{journaltitle}{Foundations of Computational Mathematics}
      \field{langid}{english}
      \field{month}{12}
      \field{number}{6}
      \field{shortjournal}{Found Comput Math}
      \field{title}{Finite {{Element Approximation}} of the {{Spectrum}} of the {{Curl Operator}} in a {{Multiply Connected Domain}}}
      \field{urlday}{1}
      \field{urlmonth}{3}
      \field{urlyear}{2023}
      \field{volume}{18}
      \field{year}{2018}
      \field{dateera}{ce}
      \field{urldateera}{ce}
      \field{pages}{1493\bibrangedash 1533}
      \range{pages}{41}
      \verb{doi}
      \verb 10.1007/s10208-018-9373-4
      \endverb
      \verb{file}
      \verb /home/rrobin/Zotero/storage/YBCVW3VH/Alonso-Rodríguez et al. - 2018 - Finite Element Approximation of the Spectrum of th.pdf
      \endverb
      \verb{urlraw}
      \verb https://doi.org/10.1007/s10208-018-9373-4
      \endverb
      \verb{url}
      \verb https://doi.org/10.1007/s10208-018-9373-4
      \endverb
      \keyw{65N15,65N25,65N30,76M10,78M10,Beltrami fields,Finite element approximation,Force-free fields,Multiply connected domain,Spectrum of curl curl operator}
    \endentry
    \entry{amestoyFullyAsynchronousMultifrontal2001}{article}{}
      \name{author}{4}{}{%
        {{hash=067d758421e08941aee47418ed7998ec}{%
           family={Amestoy},
           familyi={A\bibinitperiod},
           given={P.\bibnamedelimi R.},
           giveni={P\bibinitperiod\bibinitdelim R\bibinitperiod}}}%
        {{hash=b3f9ecf41ae97b32ccc5361e3f213817}{%
           family={Duff},
           familyi={D\bibinitperiod},
           given={I.\bibnamedelimi S.},
           giveni={I\bibinitperiod\bibinitdelim S\bibinitperiod}}}%
        {{hash=d702bfb17935a9626bc0c37201bb4ead}{%
           family={Koster},
           familyi={K\bibinitperiod},
           given={J.},
           giveni={J\bibinitperiod}}}%
        {{hash=a0205248e5b9727198823971807661e0}{%
           family={L'Excellent},
           familyi={L\bibinitperiod},
           given={J.-Y.},
           giveni={J\bibinithyphendelim Y\bibinitperiod}}}%
      }
      \strng{namehash}{5273fd116fa18da4962850bf09c148dc}
      \strng{fullhash}{5273fd116fa18da4962850bf09c148dc}
      \strng{bibnamehash}{5273fd116fa18da4962850bf09c148dc}
      \strng{authorbibnamehash}{5273fd116fa18da4962850bf09c148dc}
      \strng{authornamehash}{5273fd116fa18da4962850bf09c148dc}
      \strng{authorfullhash}{5273fd116fa18da4962850bf09c148dc}
      \field{labelalpha}{Ame+01}
      \field{sortinit}{A}
      \field{sortinithash}{a3dcedd53b04d1adfd5ac303ecd5e6fa}
      \field{labelnamesource}{author}
      \field{labeltitlesource}{title}
      \field{journaltitle}{SIAM Journal on Matrix Analysis and Applications}
      \field{number}{1}
      \field{title}{A fully asynchronous multifrontal solver using distributed dynamic scheduling}
      \field{volume}{23}
      \field{year}{2001}
      \field{dateera}{ce}
      \field{pages}{15\bibrangedash 41}
      \range{pages}{27}
    \endentry
    \entry{amestoyHybridSchedulingParallel2006}{article}{}
      \name{author}{4}{}{%
        {{hash=067d758421e08941aee47418ed7998ec}{%
           family={Amestoy},
           familyi={A\bibinitperiod},
           given={P.\bibnamedelimi R.},
           giveni={P\bibinitperiod\bibinitdelim R\bibinitperiod}}}%
        {{hash=fd7092fd6435d6aaa80bdee340f48c04}{%
           family={Guermouche},
           familyi={G\bibinitperiod},
           given={A.},
           giveni={A\bibinitperiod}}}%
        {{hash=a0205248e5b9727198823971807661e0}{%
           family={L'Excellent},
           familyi={L\bibinitperiod},
           given={J.-Y.},
           giveni={J\bibinithyphendelim Y\bibinitperiod}}}%
        {{hash=94f1d77ffc0f677a14f1510af83ce945}{%
           family={Pralet},
           familyi={P\bibinitperiod},
           given={S.},
           giveni={S\bibinitperiod}}}%
      }
      \strng{namehash}{4d8936a83c29237a4ac9c4027c7c94b9}
      \strng{fullhash}{4d8936a83c29237a4ac9c4027c7c94b9}
      \strng{bibnamehash}{4d8936a83c29237a4ac9c4027c7c94b9}
      \strng{authorbibnamehash}{4d8936a83c29237a4ac9c4027c7c94b9}
      \strng{authornamehash}{4d8936a83c29237a4ac9c4027c7c94b9}
      \strng{authorfullhash}{4d8936a83c29237a4ac9c4027c7c94b9}
      \field{labelalpha}{Ame+06}
      \field{sortinit}{A}
      \field{sortinithash}{a3dcedd53b04d1adfd5ac303ecd5e6fa}
      \field{labelnamesource}{author}
      \field{labeltitlesource}{title}
      \field{journaltitle}{Parallel Computing}
      \field{number}{2}
      \field{title}{Hybrid scheduling for the parallel solution of linear systems}
      \field{volume}{32}
      \field{year}{2006}
      \field{dateera}{ce}
      \field{pages}{136\bibrangedash 156}
      \range{pages}{21}
    \endentry
    \entry{amroucheVectorPotentialThreedimensional1998}{article}{}
      \name{author}{4}{}{%
        {{hash=7ac0edcd37bfe5dd60c707e390df11b9}{%
           family={Amrouche},
           familyi={A\bibinitperiod},
           given={Cherif},
           giveni={C\bibinitperiod}}}%
        {{hash=48ff98f52ce53f6a6911d314d7052ac0}{%
           family={Bernardi},
           familyi={B\bibinitperiod},
           given={Christine},
           giveni={C\bibinitperiod}}}%
        {{hash=2bc7d57aa6c07734f2ea532bebd442da}{%
           family={Dauge},
           familyi={D\bibinitperiod},
           given={Monique},
           giveni={M\bibinitperiod}}}%
        {{hash=52484497cc24f688e1ca110fb5cf0675}{%
           family={Girault},
           familyi={G\bibinitperiod},
           given={Vivette},
           giveni={V\bibinitperiod}}}%
      }
      \strng{namehash}{1afd2e0bd471e7870b83424c5cd625e3}
      \strng{fullhash}{1afd2e0bd471e7870b83424c5cd625e3}
      \strng{bibnamehash}{1afd2e0bd471e7870b83424c5cd625e3}
      \strng{authorbibnamehash}{1afd2e0bd471e7870b83424c5cd625e3}
      \strng{authornamehash}{1afd2e0bd471e7870b83424c5cd625e3}
      \strng{authorfullhash}{1afd2e0bd471e7870b83424c5cd625e3}
      \field{labelalpha}{Amr+98}
      \field{sortinit}{A}
      \field{sortinithash}{a3dcedd53b04d1adfd5ac303ecd5e6fa}
      \field{labelnamesource}{author}
      \field{labeltitlesource}{title}
      \field{abstract}{This paper presents several results concerning the vector potential which can be associated with a divergence-free function in a bounded three-dimensional domain. Different types of boundary conditions are given, for which the existence, uniqueness and regularity of the potential are studied. This is applied firstly to the finite element discretization of these potentials and secondly to a new formulation of incompressible viscous flow problems.On présente dans cet article un certain nombre de résultats concernant le potentiel vecteur associé à une fonction à divergence nulle dans un ouvert borné de dimension trois. En particulier, plusieurs types de conditions aux limites sont proposés, pour lesquels on discute l'existence, l'unicité et la régularité du potentiel vecteur. On analyse la convergence d'une discrétisation par éléments finis de ces potentiels et on indique une application concernant l'approximation de fluides visqueux incompressibles.© 1998 B. G. Teubner Stuttgart—John Wiley \& Sons, Ltd.}
      \field{day}{1}
      \field{journaltitle}{Mathematical Methods in The Applied Sciences - MATH METH APPL SCI}
      \field{month}{6}
      \field{shortjournal}{Mathematical Methods in The Applied Sciences - MATH METH APPL SCI}
      \field{title}{Vector potential in three-dimensional nonsmooth domains}
      \field{volume}{21}
      \field{year}{1998}
      \field{dateera}{ce}
      \field{pages}{823\bibrangedash 864}
      \range{pages}{42}
      \verb{doi}
      \verb 10.1002/(SICI)1099-1476(199806)21:93.0.CO;2-B
      \endverb
      \verb{file}
      \verb /home/rrobin/Zotero/storage/QEIRMGVG/Amrouche et al. - 1998 - Vector potential in three-dimensional nonsmooth do.pdf
      \endverb
    \endentry
    \entry{arnoldFiniteElementExterior2010}{article}{}
      \name{author}{3}{}{%
        {{hash=f95d6102f6f20db57c2e5e3505a0c374}{%
           family={Arnold},
           familyi={A\bibinitperiod},
           given={Douglas},
           giveni={D\bibinitperiod}}}%
        {{hash=e149ed0cc036a811c5d99777be1830fc}{%
           family={Falk},
           familyi={F\bibinitperiod},
           given={Richard},
           giveni={R\bibinitperiod}}}%
        {{hash=94a4fd4da337a1d51e2e55b3ddb54002}{%
           family={Winther},
           familyi={W\bibinitperiod},
           given={Ragnar},
           giveni={R\bibinitperiod}}}%
      }
      \strng{namehash}{65cadc6ff66c31c150bdeba754fd6fa6}
      \strng{fullhash}{65cadc6ff66c31c150bdeba754fd6fa6}
      \strng{bibnamehash}{65cadc6ff66c31c150bdeba754fd6fa6}
      \strng{authorbibnamehash}{65cadc6ff66c31c150bdeba754fd6fa6}
      \strng{authornamehash}{65cadc6ff66c31c150bdeba754fd6fa6}
      \strng{authorfullhash}{65cadc6ff66c31c150bdeba754fd6fa6}
      \field{extraname}{1}
      \field{labelalpha}{AFW10}
      \field{sortinit}{A}
      \field{sortinithash}{a3dcedd53b04d1adfd5ac303ecd5e6fa}
      \field{labelnamesource}{author}
      \field{labeltitlesource}{shorttitle}
      \field{abstract}{Advancing research. Creating connections.}
      \field{issn}{0273-0979, 1088-9485}
      \field{journaltitle}{Bulletin of the American Mathematical Society}
      \field{month}{4}
      \field{number}{2}
      \field{shorttitle}{Finite element exterior calculus}
      \field{title}{Finite element exterior calculus: from {Hodge} theory to numerical stability}
      \field{urlday}{25}
      \field{urlmonth}{1}
      \field{urlyear}{2024}
      \field{volume}{47}
      \field{year}{2010}
      \field{urldateera}{ce}
      \field{pages}{281\bibrangedash 354}
      \range{pages}{74}
      \verb{doi}
      \verb 10.1090/S0273-0979-10-01278-4
      \endverb
      \verb{file}
      \verb Full Text PDF:C\:\\Users\\rorousse\\Zotero\\storage\\DPE2UYQH\\Arnold et al. - 2010 - Finite element exterior calculus from Hodge theor.pdf:application/pdf
      \endverb
      \verb{urlraw}
      \verb https://www.ams.org/bull/2010-47-02/S0273-0979-10-01278-4/
      \endverb
      \verb{url}
      \verb https://www.ams.org/bull/2010-47-02/S0273-0979-10-01278-4/
      \endverb
      \keyw{de Rham cohomology,exterior calculus,Finite element exterior calculus,Hodge Laplacian,Hodge theory,mixed finite elements}
    \endentry
    \entry{arnoldFiniteElementExterior2006}{article}{}
      \name{author}{3}{}{%
        {{hash=2c0d5891c854db0dc6743a2f93208c49}{%
           family={Arnold},
           familyi={A\bibinitperiod},
           given={Douglas\bibnamedelima N.},
           giveni={D\bibinitperiod\bibinitdelim N\bibinitperiod}}}%
        {{hash=aea4e093a0914a227674c914b4692a6d}{%
           family={Falk},
           familyi={F\bibinitperiod},
           given={Richard\bibnamedelima S.},
           giveni={R\bibinitperiod\bibinitdelim S\bibinitperiod}}}%
        {{hash=94a4fd4da337a1d51e2e55b3ddb54002}{%
           family={Winther},
           familyi={W\bibinitperiod},
           given={Ragnar},
           giveni={R\bibinitperiod}}}%
      }
      \strng{namehash}{f7166307ec336227a17b50cf80131d50}
      \strng{fullhash}{f7166307ec336227a17b50cf80131d50}
      \strng{bibnamehash}{f7166307ec336227a17b50cf80131d50}
      \strng{authorbibnamehash}{f7166307ec336227a17b50cf80131d50}
      \strng{authornamehash}{f7166307ec336227a17b50cf80131d50}
      \strng{authorfullhash}{f7166307ec336227a17b50cf80131d50}
      \field{extraname}{2}
      \field{labelalpha}{AFW06}
      \field{sortinit}{A}
      \field{sortinithash}{a3dcedd53b04d1adfd5ac303ecd5e6fa}
      \field{labelnamesource}{author}
      \field{labeltitlesource}{title}
      \field{abstract}{Finite element exterior calculus is an approach to the design and understanding of finite element discretizations for a wide variety of systems of partial differential equations. This approach brings to bear tools from differential geometry, algebraic topology, and homological algebra to develop discretizations which are compatible with the geometric, topological, and algebraic structures which underlie well-posedness of the PDE problem being solved. In the finite element exterior calculus, many finite element spaces are revealed as spaces of piecewise polynomial differential forms. These connect to each other in discrete subcomplexes of elliptic differential complexes, and are also related to the continuous elliptic complex through projections which commute with the complex differential. Applications are made to the finite element discretization of a variety of problems, including the Hodge Laplacian, Maxwell’s equations, the equations of elasticity, and elliptic eigenvalue problems, and also to preconditioners.}
      \field{issn}{0962-4929, 1474-0508}
      \field{journaltitle}{Acta Numerica}
      \field{langid}{english}
      \field{month}{5}
      \field{shortjournal}{Acta Numerica}
      \field{title}{Finite element exterior calculus, homological techniques, and applications}
      \field{urlday}{2}
      \field{urlmonth}{5}
      \field{urlyear}{2023}
      \field{volume}{15}
      \field{year}{2006}
      \field{dateera}{ce}
      \field{urldateera}{ce}
      \field{pages}{1\bibrangedash 155}
      \range{pages}{155}
      \verb{doi}
      \verb 10.1017/S0962492906210018
      \endverb
      \verb{file}
      \verb /home/rrobin/Zotero/storage/RU3HPI5C/Arnold et al. - 2006 - Finite element exterior calculus, homological tech.pdf
      \endverb
      \verb{urlraw}
      \verb https://www.cambridge.org/core/product/identifier/S0962492906210018/type/journal_article
      \endverb
      \verb{url}
      \verb https://www.cambridge.org/core/product/identifier/S0962492906210018/type/journal_article
      \endverb
    \endentry
    \entry{arnoldTopologyThreedimensionalSteady1966}{article}{}
      \name{author}{1}{}{%
        {{hash=ff4159d339a8205fcaf83c1c04fc4275}{%
           family={Arnold},
           familyi={A\bibinitperiod},
           given={V.\bibnamedelimi I.},
           giveni={V\bibinitperiod\bibinitdelim I\bibinitperiod}}}%
      }
      \strng{namehash}{ff4159d339a8205fcaf83c1c04fc4275}
      \strng{fullhash}{ff4159d339a8205fcaf83c1c04fc4275}
      \strng{bibnamehash}{ff4159d339a8205fcaf83c1c04fc4275}
      \strng{authorbibnamehash}{ff4159d339a8205fcaf83c1c04fc4275}
      \strng{authornamehash}{ff4159d339a8205fcaf83c1c04fc4275}
      \strng{authorfullhash}{ff4159d339a8205fcaf83c1c04fc4275}
      \field{extraname}{1}
      \field{labelalpha}{Arn66}
      \field{sortinit}{A}
      \field{sortinithash}{a3dcedd53b04d1adfd5ac303ecd5e6fa}
      \field{labelnamesource}{author}
      \field{labeltitlesource}{title}
      \field{day}{1}
      \field{issn}{0021-8928}
      \field{journaltitle}{Journal of Applied Mathematics and Mechanics}
      \field{month}{1}
      \field{number}{1}
      \field{shortjournal}{Journal of Applied Mathematics and Mechanics}
      \field{title}{On the topology of three-dimensional steady flows of an ideal fluid}
      \field{urlday}{21}
      \field{urlmonth}{11}
      \field{urlyear}{2023}
      \field{volume}{30}
      \field{year}{1966}
      \field{dateera}{ce}
      \field{urldateera}{ce}
      \field{pages}{223\bibrangedash 226}
      \range{pages}{4}
      \verb{doi}
      \verb 10.1016/0021-8928(66)90070-0
      \endverb
      \verb{file}
      \verb /home/rrobin/Zotero/storage/WBEPJYDK/Arnol'd - 1966 - On the topology of three-dimensional steady flows .pdf;/home/rrobin/Zotero/storage/XLSBGUVW/0021892866900700.html
      \endverb
      \verb{urlraw}
      \verb https://www.sciencedirect.com/science/article/pii/0021892866900700
      \endverb
      \verb{url}
      \verb https://www.sciencedirect.com/science/article/pii/0021892866900700
      \endverb
    \endentry
    \entry{Arnold2014}{incollection}{}
      \name{author}{1}{}{%
        {{hash=92d7efabd161dd52a00501a33ccb997d}{%
           family={Arnold},
           familyi={A\bibinitperiod},
           given={Vladimir\bibnamedelima I.},
           giveni={V\bibinitperiod\bibinitdelim I\bibinitperiod}}}%
      }
      \name{editor}{5}{}{%
        {{hash=761613632efc6050cb4d3d161ccdb749}{%
           family={Givental},
           familyi={G\bibinitperiod},
           given={Alexander\bibnamedelima B.},
           giveni={A\bibinitperiod\bibinitdelim B\bibinitperiod}}}%
        {{hash=8dd0b5fcc69b4af017b6f79d00f0534e}{%
           family={Khesin},
           familyi={K\bibinitperiod},
           given={Boris\bibnamedelima A.},
           giveni={B\bibinitperiod\bibinitdelim A\bibinitperiod}}}%
        {{hash=93d2e0f87d0f9726a79e089336ec16ae}{%
           family={Varchenko},
           familyi={V\bibinitperiod},
           given={Alexander\bibnamedelima N.},
           giveni={A\bibinitperiod\bibinitdelim N\bibinitperiod}}}%
        {{hash=360eee41486b7fd090506a244d6bfacd}{%
           family={Vassiliev},
           familyi={V\bibinitperiod},
           given={Victor\bibnamedelima A.},
           giveni={V\bibinitperiod\bibinitdelim A\bibinitperiod}}}%
        {{hash=d9e219194ba6049646322e4a4277176b}{%
           family={Viro},
           familyi={V\bibinitperiod},
           given={Oleg\bibnamedelima Ya.},
           giveni={O\bibinitperiod\bibinitdelim Y\bibinitperiod}}}%
      }
      \list{location}{1}{%
        {Berlin, Heidelberg}%
      }
      \list{publisher}{1}{%
        {Springer Berlin Heidelberg}%
      }
      \strng{namehash}{92d7efabd161dd52a00501a33ccb997d}
      \strng{fullhash}{92d7efabd161dd52a00501a33ccb997d}
      \strng{bibnamehash}{92d7efabd161dd52a00501a33ccb997d}
      \strng{authorbibnamehash}{92d7efabd161dd52a00501a33ccb997d}
      \strng{authornamehash}{92d7efabd161dd52a00501a33ccb997d}
      \strng{authorfullhash}{92d7efabd161dd52a00501a33ccb997d}
      \strng{editorbibnamehash}{e9138ba87cdd9f83217923e65088dff9}
      \strng{editornamehash}{e9138ba87cdd9f83217923e65088dff9}
      \strng{editorfullhash}{e193e9ae6d058d6a092a5bab45edbef0}
      \field{extraname}{2}
      \field{labelalpha}{Arn14}
      \field{sortinit}{A}
      \field{sortinithash}{a3dcedd53b04d1adfd5ac303ecd5e6fa}
      \field{labelnamesource}{author}
      \field{labeltitlesource}{title}
      \field{abstract}{The classical Hopf invariant distinguishes among the homotopy classes of continuous mappings from the three-sphere to the two-sphere and is equal to the linking number of the two curves that are the preimages of any two regular points of the two-sphere.}
      \field{booktitle}{Vladimir {{I}}. {{Arnold}} - collected works: {{Hydrodynamics}}, bifurcation theory, and algebraic geometry 1965-1972}
      \field{isbn}{978-3-642-31031-7}
      \field{title}{The asymptotic {{Hopf}} invariant and its applications}
      \field{year}{2014}
      \field{dateera}{ce}
      \field{pages}{357\bibrangedash 375}
      \range{pages}{19}
      \verb{doi}
      \verb 10.1007/978-3-642-31031-7_32
      \endverb
      \verb{urlraw}
      \verb https://doi.org/10.1007/978-3-642-31031-7_32
      \endverb
      \verb{url}
      \verb https://doi.org/10.1007/978-3-642-31031-7_32
      \endverb
    \endentry
    \entry{arnoldTopologicalMethodsHydrodynamics2021}{book}{}
      \name{author}{2}{}{%
        {{hash=92d7efabd161dd52a00501a33ccb997d}{%
           family={Arnold},
           familyi={A\bibinitperiod},
           given={Vladimir\bibnamedelima I.},
           giveni={V\bibinitperiod\bibinitdelim I\bibinitperiod}}}%
        {{hash=8dd0b5fcc69b4af017b6f79d00f0534e}{%
           family={Khesin},
           familyi={K\bibinitperiod},
           given={Boris\bibnamedelima A.},
           giveni={B\bibinitperiod\bibinitdelim A\bibinitperiod}}}%
      }
      \list{location}{1}{%
        {Cham}%
      }
      \list{publisher}{1}{%
        {Springer International Publishing}%
      }
      \strng{namehash}{b56c9851fea889f81e228860dc34db81}
      \strng{fullhash}{b56c9851fea889f81e228860dc34db81}
      \strng{bibnamehash}{b56c9851fea889f81e228860dc34db81}
      \strng{authorbibnamehash}{b56c9851fea889f81e228860dc34db81}
      \strng{authornamehash}{b56c9851fea889f81e228860dc34db81}
      \strng{authorfullhash}{b56c9851fea889f81e228860dc34db81}
      \field{labelalpha}{AK21}
      \field{sortinit}{A}
      \field{sortinithash}{a3dcedd53b04d1adfd5ac303ecd5e6fa}
      \field{labelnamesource}{author}
      \field{labeltitlesource}{title}
      \field{isbn}{978-3-030-74277-5 978-3-030-74278-2}
      \field{langid}{english}
      \field{series}{Applied {{Mathematical Sciences}}}
      \field{title}{Topological {{Methods}} in {{Hydrodynamics}}}
      \field{urlday}{21}
      \field{urlmonth}{11}
      \field{urlyear}{2023}
      \field{volume}{125}
      \field{year}{2021}
      \field{dateera}{ce}
      \field{urldateera}{ce}
      \verb{doi}
      \verb 10.1007/978-3-030-74278-2
      \endverb
      \verb{file}
      \verb /home/rrobin/Zotero/storage/8LUMZHC2/Arnold et Khesin - 2021 - Topological Methods in Hydrodynamics.pdf
      \endverb
      \verb{urlraw}
      \verb https://link.springer.com/10.1007/978-3-030-74278-2
      \endverb
      \verb{url}
      \verb https://link.springer.com/10.1007/978-3-030-74278-2
      \endverb
    \endentry
    \entry{balay1998petsc}{article}{}
      \name{author}{4}{}{%
        {{hash=0661bc042688479e939e81e7f9185f66}{%
           family={Balay},
           familyi={B\bibinitperiod},
           given={Satish},
           giveni={S\bibinitperiod}}}%
        {{hash=232270861ce704a354417bfda0acdb1b}{%
           family={Gropp},
           familyi={G\bibinitperiod},
           given={William},
           giveni={W\bibinitperiod}}}%
        {{hash=65f9dd1eb923bd39cbd7ae2d5bc637cd}{%
           family={McInnes},
           familyi={M\bibinitperiod},
           given={Lois\bibnamedelima Curfman},
           giveni={L\bibinitperiod\bibinitdelim C\bibinitperiod}}}%
        {{hash=db20e2f433b4eee5280b269cfbfcc84f}{%
           family={Smith},
           familyi={S\bibinitperiod},
           given={Barry\bibnamedelima F},
           giveni={B\bibinitperiod\bibinitdelim F\bibinitperiod}}}%
      }
      \strng{namehash}{87719e982f89ddec1928a8e37cd2ebcf}
      \strng{fullhash}{87719e982f89ddec1928a8e37cd2ebcf}
      \strng{bibnamehash}{87719e982f89ddec1928a8e37cd2ebcf}
      \strng{authorbibnamehash}{87719e982f89ddec1928a8e37cd2ebcf}
      \strng{authornamehash}{87719e982f89ddec1928a8e37cd2ebcf}
      \strng{authorfullhash}{87719e982f89ddec1928a8e37cd2ebcf}
      \field{labelalpha}{Bal+98}
      \field{sortinit}{B}
      \field{sortinithash}{8de16967003c7207dae369d874f1456e}
      \field{labelnamesource}{author}
      \field{labeltitlesource}{title}
      \field{journaltitle}{Argonne National Laboratory}
      \field{number}{17}
      \field{title}{{{PETSc}}, the portable, extensible toolkit for scientific computation}
      \field{volume}{2}
      \field{year}{1998}
      \field{dateera}{ce}
    \endentry
    \entry{bevirRelaxationFluxConsumption1980}{inproceedings}{}
      \name{author}{2}{}{%
        {{hash=323239847509f2eb0c54fada2f162ccf}{%
           family={Bevir},
           familyi={B\bibinitperiod},
           given={M.\bibnamedelimi K.},
           giveni={M\bibinitperiod\bibinitdelim K\bibinitperiod}}}%
        {{hash=d8ac9302e7fe577b1efcd6546c17c9ce}{%
           family={Gray},
           familyi={G\bibinitperiod},
           given={J.\bibnamedelimi W.},
           giveni={J\bibinitperiod\bibinitdelim W\bibinitperiod}}}%
      }
      \list{location}{1}{%
        {Los Alamos}%
      }
      \strng{namehash}{96c73bae844de8365f9379db17b319c7}
      \strng{fullhash}{96c73bae844de8365f9379db17b319c7}
      \strng{bibnamehash}{96c73bae844de8365f9379db17b319c7}
      \strng{authorbibnamehash}{96c73bae844de8365f9379db17b319c7}
      \strng{authornamehash}{96c73bae844de8365f9379db17b319c7}
      \strng{authorfullhash}{96c73bae844de8365f9379db17b319c7}
      \field{labelalpha}{BG80}
      \field{sortinit}{B}
      \field{sortinithash}{8de16967003c7207dae369d874f1456e}
      \field{labelnamesource}{author}
      \field{labeltitlesource}{title}
      \field{booktitle}{Proceedings of the {{Reversed-Field}}}
      \field{eventtitle}{Pinch {{Theory Workshop}}}
      \field{title}{Relaxation, flux consumption and quasi steady state pinches}
      \field{year}{1980}
      \field{dateera}{ce}
    \endentry
    \entry{buffaTracesCurlLipschitz2002}{article}{}
      \name{author}{3}{}{%
        {{hash=81aee8796fd7c7c099cdc144d9485292}{%
           family={Buffa},
           familyi={B\bibinitperiod},
           given={A.},
           giveni={A\bibinitperiod}}}%
        {{hash=ac38abf304a3808a0422ee52c5d0168e}{%
           family={Costabel},
           familyi={C\bibinitperiod},
           given={M.},
           giveni={M\bibinitperiod}}}%
        {{hash=5bf1e2fe336311a81f58c1322e12cf7e}{%
           family={Sheen},
           familyi={S\bibinitperiod},
           given={D.},
           giveni={D\bibinitperiod}}}%
      }
      \strng{namehash}{3f5a36a3876183b8f194b3eb19af6589}
      \strng{fullhash}{3f5a36a3876183b8f194b3eb19af6589}
      \strng{bibnamehash}{3f5a36a3876183b8f194b3eb19af6589}
      \strng{authorbibnamehash}{3f5a36a3876183b8f194b3eb19af6589}
      \strng{authornamehash}{3f5a36a3876183b8f194b3eb19af6589}
      \strng{authorfullhash}{3f5a36a3876183b8f194b3eb19af6589}
      \field{labelalpha}{BCS02}
      \field{sortinit}{B}
      \field{sortinithash}{8de16967003c7207dae369d874f1456e}
      \field{labelnamesource}{author}
      \field{labeltitlesource}{title}
      \field{abstract}{We study tangential vector fields on the boundary of a bounded Lipschitz domain Ω in R3. Our attention is focused on the definition of suitable Hilbert spaces corresponding to fractional Sobolev regularities and also on the construction of tangential differential operators on the non-smooth manifold. The theory is applied to the characterization of tangential traces for the space {$H(\curl, \Omega)$}. Hodge decompositions are provided for the corresponding trace spaces, and an integration by parts formula is proved.}
      \field{issn}{0022247X}
      \field{journaltitle}{Journal of Mathematical Analysis and Applications}
      \field{langid}{english}
      \field{month}{12}
      \field{number}{2}
      \field{shortjournal}{Journal of Mathematical Analysis and Applications}
      \field{title}{On traces for {$H(\curl, \Omega)$} in {{Lipschitz}} domains}
      \field{urlday}{5}
      \field{urlmonth}{6}
      \field{urlyear}{2023}
      \field{volume}{276}
      \field{year}{2002}
      \field{dateera}{ce}
      \field{urldateera}{ce}
      \field{pages}{845\bibrangedash 867}
      \range{pages}{23}
      \verb{doi}
      \verb 10.1016/S0022-247X(02)00455-9
      \endverb
      \verb{file}
      \verb /home/rrobin/Zotero/storage/HNK6A4QU/Buffa et al. - 2002 - On traces for {$H(\curl, \Omega)$} in Lipschitz domains.pdf
      \endverb
      \verb{urlraw}
      \verb https://linkinghub.elsevier.com/retrieve/pii/S0022247X02004559
      \endverb
      \verb{url}
      \verb https://linkinghub.elsevier.com/retrieve/pii/S0022247X02004559
      \endverb
    \endentry
    \entry{cantarellaBiotSavartOperator2001}{article}{}
      \name{author}{3}{}{%
        {{hash=f33510f68f461a8bfc76644a4b7dc6a4}{%
           family={Cantarella},
           familyi={C\bibinitperiod},
           given={Jason},
           giveni={J\bibinitperiod}}}%
        {{hash=95d8d67e73e48bba52a5faca2cf33505}{%
           family={DeTurck},
           familyi={D\bibinitperiod},
           given={Dennis},
           giveni={D\bibinitperiod}}}%
        {{hash=734ee77d7ef9701f6e25a4db78add6b2}{%
           family={Gluck},
           familyi={G\bibinitperiod},
           given={Herman},
           giveni={H\bibinitperiod}}}%
      }
      \strng{namehash}{0bde656bc3ccf5afb860533bc13e829a}
      \strng{fullhash}{0bde656bc3ccf5afb860533bc13e829a}
      \strng{bibnamehash}{0bde656bc3ccf5afb860533bc13e829a}
      \strng{authorbibnamehash}{0bde656bc3ccf5afb860533bc13e829a}
      \strng{authornamehash}{0bde656bc3ccf5afb860533bc13e829a}
      \strng{authorfullhash}{0bde656bc3ccf5afb860533bc13e829a}
      \field{labelalpha}{CDG01}
      \field{sortinit}{C}
      \field{sortinithash}{4c244ceae61406cdc0cc2ce1cb1ff703}
      \field{labelnamesource}{author}
      \field{labeltitlesource}{title}
      \field{issn}{00222488}
      \field{journaltitle}{Journal of Mathematical Physics}
      \field{langid}{english}
      \field{number}{2}
      \field{shortjournal}{J. Math. Phys.}
      \field{title}{The {{Biot}}–{{Savart}} operator for application to knot theory, fluid dynamics, and plasma physics}
      \field{urlday}{12}
      \field{urlmonth}{5}
      \field{urlyear}{2020}
      \field{volume}{42}
      \field{year}{2001}
      \field{dateera}{ce}
      \field{urldateera}{ce}
      \verb{doi}
      \verb 10.1063/1.1329659
      \endverb
      \verb{file}
      \verb /home/rrobin/Zotero/storage/9FEC9366/Cantarella et al. - 2001 - The Biot–Savart operator for application to knot t.pdf;/home/rrobin/Zotero/storage/AGMPD2QT/Cantarella et al. - 2001 - The Biot–Savart operator for application to knot t.pdf;/home/rrobin/Zotero/storage/BAHUHPUJ/Cantarella et al. - 2001 - The Biot–Savart operator for application to knot t.pdf;/home/rrobin/Zotero/storage/NAT39WZK/Cantarella et al. - 2001 - The Biot–Savart operator for application to knot t.pdf;/home/rrobin/Zotero/storage/WZ6IBVTE/Cantarella et al. - 2001 - The Biot–Savart operator for application to knot t.pdf;/home/rrobin/Zotero/storage/XSXGN53F/Cantarella et al. - 2001 - The Biot–Savart operator for application to knot t.pdf
      \endverb
      \verb{urlraw}
      \verb http://scitation.aip.org/content/aip/journal/jmp/42/2/10.1063/1.1329659
      \endverb
      \verb{url}
      \verb http://scitation.aip.org/content/aip/journal/jmp/42/2/10.1063/1.1329659
      \endverb
    \endentry
    \entry{cantarellaIsoperimetricProblemsHelicity2000}{article}{}
      \name{author}{4}{}{%
        {{hash=f33510f68f461a8bfc76644a4b7dc6a4}{%
           family={Cantarella},
           familyi={C\bibinitperiod},
           given={Jason},
           giveni={J\bibinitperiod}}}%
        {{hash=95d8d67e73e48bba52a5faca2cf33505}{%
           family={DeTurck},
           familyi={D\bibinitperiod},
           given={Dennis},
           giveni={D\bibinitperiod}}}%
        {{hash=734ee77d7ef9701f6e25a4db78add6b2}{%
           family={Gluck},
           familyi={G\bibinitperiod},
           given={Herman},
           giveni={H\bibinitperiod}}}%
        {{hash=498ba0cb013cf4647b2da42caf91e88f}{%
           family={Teytel},
           familyi={T\bibinitperiod},
           given={Mikhail},
           giveni={M\bibinitperiod}}}%
      }
      \strng{namehash}{c0547afbfd8c7041353961bc5d211c6e}
      \strng{fullhash}{c0547afbfd8c7041353961bc5d211c6e}
      \strng{bibnamehash}{c0547afbfd8c7041353961bc5d211c6e}
      \strng{authorbibnamehash}{c0547afbfd8c7041353961bc5d211c6e}
      \strng{authornamehash}{c0547afbfd8c7041353961bc5d211c6e}
      \strng{authorfullhash}{c0547afbfd8c7041353961bc5d211c6e}
      \field{labelalpha}{Can+00}
      \field{sortinit}{C}
      \field{sortinithash}{4c244ceae61406cdc0cc2ce1cb1ff703}
      \field{labelnamesource}{author}
      \field{labeltitlesource}{title}
      \field{issn}{0022-2488, 1089-7658}
      \field{journaltitle}{Journal of Mathematical Physics}
      \field{langid}{english}
      \field{month}{8}
      \field{number}{8}
      \field{shortjournal}{Journal of Mathematical Physics}
      \field{title}{Isoperimetric problems for the helicity of vector fields and the {{Biot}}–{{Savart}} and curl operators}
      \field{urlday}{13}
      \field{urlmonth}{5}
      \field{urlyear}{2020}
      \field{volume}{41}
      \field{year}{2000}
      \field{dateera}{ce}
      \field{urldateera}{ce}
      \field{pages}{5615\bibrangedash 5641}
      \range{pages}{27}
      \verb{doi}
      \verb 10.1063/1.533429
      \endverb
      \verb{file}
      \verb /home/rrobin/Zotero/storage/3PGF5MYZ/Cantarella et al. - 2000 - Isoperimetric problems for the helicity of vector .pdf;/home/rrobin/Zotero/storage/5FHPPX3L/Cantarella et al. - 2000 - Isoperimetric problems for the helicity of vector .pdf;/home/rrobin/Zotero/storage/FVW6Q6V8/Cantarella et al. - 2000 - Isoperimetric problems for the helicity of vector .pdf;/home/rrobin/Zotero/storage/JDKJZ5NI/Cantarella et al. - 2000 - Isoperimetric problems for the helicity of vector .pdf;/home/rrobin/Zotero/storage/P82NQPBA/Cantarella et al. - 2000 - Isoperimetric problems for the helicity of vector .pdf
      \endverb
      \verb{urlraw}
      \verb http://aip.scitation.org/doi/10.1063/1.533429
      \endverb
      \verb{url}
      \verb http://aip.scitation.org/doi/10.1063/1.533429
      \endverb
    \endentry
    \entry{cantarella_influence_1999}{incollection}{}
      \name{author}{4}{}{%
        {{hash=f33510f68f461a8bfc76644a4b7dc6a4}{%
           family={Cantarella},
           familyi={C\bibinitperiod},
           given={Jason},
           giveni={J\bibinitperiod}}}%
        {{hash=df1bf8e73e4c0cd8ce9dbb34f689ba8f}{%
           family={Deturck},
           familyi={D\bibinitperiod},
           given={Dennis},
           giveni={D\bibinitperiod}}}%
        {{hash=734ee77d7ef9701f6e25a4db78add6b2}{%
           family={Gluck},
           familyi={G\bibinitperiod},
           given={Herman},
           giveni={H\bibinitperiod}}}%
        {{hash=498ba0cb013cf4647b2da42caf91e88f}{%
           family={Teytel},
           familyi={T\bibinitperiod},
           given={Mikhail},
           giveni={M\bibinitperiod}}}%
      }
      \list{publisher}{1}{%
        {American Geophysical Union (AGU)}%
      }
      \strng{namehash}{b80dce2bae5edf46f56f6d6a6f0e2e1a}
      \strng{fullhash}{b80dce2bae5edf46f56f6d6a6f0e2e1a}
      \strng{bibnamehash}{b80dce2bae5edf46f56f6d6a6f0e2e1a}
      \strng{authorbibnamehash}{b80dce2bae5edf46f56f6d6a6f0e2e1a}
      \strng{authornamehash}{b80dce2bae5edf46f56f6d6a6f0e2e1a}
      \strng{authorfullhash}{b80dce2bae5edf46f56f6d6a6f0e2e1a}
      \field{labelalpha}{Can+99}
      \field{sortinit}{C}
      \field{sortinithash}{4c244ceae61406cdc0cc2ce1cb1ff703}
      \field{labelnamesource}{author}
      \field{labeltitlesource}{title}
      \field{abstract}{This chapter contains sections titled: Two Fundamental Problems Hellcity and Writhing Number Relation between Helicity and Writhing Number How the Geometry of the Domain Influences Helicity Magnetic Fields and Helicity A General Point of View The Modified Biot-Savart Operator Spectral Methods Connection with the Curl Operator Explicit Computation of Energy-Minimizing Vector Fields The Isoperimetric Problem First Variation Formulas Constraints on Any Optimal Domain The Search for Optimal Domains Appendix. The Hodge Decomposition Theorem}
      \field{booktitle}{Magnetic {{Helicity}} in {{Space}} and {{Laboratory Plasmas}}}
      \field{isbn}{978-1-118-66447-6}
      \field{langid}{english}
      \field{title}{Influence of {{Geometry}} and {{Topology}} on {{Helicity}}}
      \field{urlday}{1}
      \field{urlmonth}{6}
      \field{urlyear}{2023}
      \field{year}{1999}
      \field{dateera}{ce}
      \field{urldateera}{ce}
      \field{pages}{17\bibrangedash 24}
      \range{pages}{8}
      \verb{doi}
      \verb 10.1029/GM111p0017
      \endverb
      \verb{urlraw}
      \verb https://onlinelibrary.wiley.com/doi/abs/10.1029/GM111p0017
      \endverb
      \verb{url}
      \verb https://onlinelibrary.wiley.com/doi/abs/10.1029/GM111p0017
      \endverb
      \keyw{Magnetic reconnection,Particles (Nuclear physics)-Helicity,Plasma (Ionized gases),Plasma astrophysics}
    \endentry
    \entry{dewarStellaratorSymmetry1998}{article}{}
      \name{author}{2}{}{%
        {{hash=7b0c739a8a7a8ab95c58f6eb4cb104c3}{%
           family={Dewar},
           familyi={D\bibinitperiod},
           given={R.L.},
           giveni={R\bibinitperiod}}}%
        {{hash=81a333de59451f9498d9936ae7ae7021}{%
           family={Hudson},
           familyi={H\bibinitperiod},
           given={S.R.},
           giveni={S\bibinitperiod}}}%
      }
      \strng{namehash}{8d539dee97d722879e41e5690095edcd}
      \strng{fullhash}{8d539dee97d722879e41e5690095edcd}
      \strng{bibnamehash}{8d539dee97d722879e41e5690095edcd}
      \strng{authorbibnamehash}{8d539dee97d722879e41e5690095edcd}
      \strng{authornamehash}{8d539dee97d722879e41e5690095edcd}
      \strng{authorfullhash}{8d539dee97d722879e41e5690095edcd}
      \field{labelalpha}{DH98}
      \field{sortinit}{D}
      \field{sortinithash}{c438b3d5d027251ba63f5ed538d98af5}
      \field{labelnamesource}{author}
      \field{labeltitlesource}{title}
      \field{abstract}{A simple and general definition of stellarator symmetry is presented and its relation to previous definitions discussed. It is shown that the field-line flow in systems possessing stellarator symmetry is time-reversal invariant if the the toroidal angle is regarded as “time”.}
      \field{issn}{01672789}
      \field{journaltitle}{Physica D: Nonlinear Phenomena}
      \field{langid}{english}
      \field{number}{1-2}
      \field{shortjournal}{Physica D: Nonlinear Phenomena}
      \field{title}{Stellarator symmetry}
      \field{urlday}{27}
      \field{urlmonth}{3}
      \field{urlyear}{2021}
      \field{volume}{112}
      \field{year}{1998}
      \field{dateera}{ce}
      \field{urldateera}{ce}
      \field{pages}{275\bibrangedash 280}
      \range{pages}{6}
      \verb{doi}
      \verb 10.1016/S0167-2789(97)00216-9
      \endverb
      \verb{file}
      \verb /home/rrobin/Zotero/storage/D6GMTCRI/Dewar et Hudson - 1998 - Stellarator symmetry.pdf;/home/rrobin/Zotero/storage/DE7B56RN/Dewar et Hudson - 1998 - Stellarator symmetry.pdf;/home/rrobin/Zotero/storage/MKUFWNM6/Dewar et Hudson - 1998 - Stellarator symmetry.pdf;/home/rrobin/Zotero/storage/YXCD25E4/Dewar et Hudson - 1998 - Stellarator symmetry.pdf;/home/rrobin/Zotero/storage/FQ54B3FU/S0167278997002169.html
      \endverb
      \verb{urlraw}
      \verb https://linkinghub.elsevier.com/retrieve/pii/S0167278997002169
      \endverb
      \verb{url}
      \verb https://linkinghub.elsevier.com/retrieve/pii/S0167278997002169
      \endverb
    \endentry
    \entry{encisoNonexistenceAxisymmetricOptimal2020}{article}{}
      \name{author}{2}{}{%
        {{hash=0e8f2c2f78978ae3889dfad261f3cf00}{%
           family={Enciso},
           familyi={E\bibinitperiod},
           given={Alberto},
           giveni={A\bibinitperiod}}}%
        {{hash=009fe385c808a2b29efee58424979b41}{%
           family={Peralta-Salas},
           familyi={P\bibinithyphendelim S\bibinitperiod},
           given={Daniel},
           giveni={D\bibinitperiod}}}%
      }
      \strng{namehash}{9c2ae77efeb2db0d2be214cfba638b9d}
      \strng{fullhash}{9c2ae77efeb2db0d2be214cfba638b9d}
      \strng{bibnamehash}{9c2ae77efeb2db0d2be214cfba638b9d}
      \strng{authorbibnamehash}{9c2ae77efeb2db0d2be214cfba638b9d}
      \strng{authornamehash}{9c2ae77efeb2db0d2be214cfba638b9d}
      \strng{authorfullhash}{9c2ae77efeb2db0d2be214cfba638b9d}
      \field{labelalpha}{EP23}
      \field{sortinit}{E}
      \field{sortinithash}{c554bd1a0b76ea92b9f105fe36d9c7b0}
      \field{labelnamesource}{author}
      \field{labeltitlesource}{title}
      \field{issn}{2036-2145}
      \field{journaltitle}{ANNALI SCUOLA NORMALE SUPERIORE - CLASSE DI SCIENZE}
      \field{title}{Non-existence of axisymmetric optimal domains with smooth boundary for the first curl eigenvalue}
      \field{urlday}{25}
      \field{urlmonth}{1}
      \field{urlyear}{2024}
      \field{year}{2023}
      \field{urldateera}{ce}
      \field{pages}{311\bibrangedash 327}
      \range{pages}{17}
      \verb{doi}
      \verb 10.2422/2036-2145.202010_008
      \endverb
      \verb{file}
      \verb Full Text PDF:C\:\\Users\\rorousse\\Zotero\\storage\\5ANIBSV7\\Enciso et Peralta-Salas - 2023 - Non-existence of axisymmetric optimal domains with.pdf:application/pdf
      \endverb
      \verb{urlraw}
      \verb https://journals.sns.it/index.php/annaliscienze/article/view/5739
      \endverb
      \verb{url}
      \verb https://journals.sns.it/index.php/annaliscienze/article/view/5739
      \endverb
    \endentry
    \entry{ern_analysis_2018}{article}{}
      \name{author}{2}{}{%
        {{hash=e85a4a72db4d523e27330f9eacdfc4a7}{%
           family={Ern},
           familyi={E\bibinitperiod},
           given={Alexandre},
           giveni={A\bibinitperiod}}}%
        {{hash=48c9c127feed47f563af418f0ac404bf}{%
           family={Guermond},
           familyi={G\bibinitperiod},
           given={Jean-Luc},
           giveni={J\bibinithyphendelim L\bibinitperiod}}}%
      }
      \strng{namehash}{eedb8ad8b78016d85df98a4344f910db}
      \strng{fullhash}{eedb8ad8b78016d85df98a4344f910db}
      \strng{bibnamehash}{eedb8ad8b78016d85df98a4344f910db}
      \strng{authorbibnamehash}{eedb8ad8b78016d85df98a4344f910db}
      \strng{authornamehash}{eedb8ad8b78016d85df98a4344f910db}
      \strng{authorfullhash}{eedb8ad8b78016d85df98a4344f910db}
      \field{labelalpha}{EG18}
      \field{sortinit}{E}
      \field{sortinithash}{c554bd1a0b76ea92b9f105fe36d9c7b0}
      \field{labelnamesource}{author}
      \field{labeltitlesource}{title}
      \field{abstract}{We derive H(curl)-error estimates and improved L2-error estimates for the Maxwell equations approximated using edge finite elements. These estimates only invoke the expected regularity pickup of the exact solution in the scale of the Sobolev spaces, which is typically lower than 12 and can be arbitrarily close to 0 when the material properties are heterogeneous. The key tools for the analysis are commuting quasi-interpolation operators in H(curl)- and H(div)-conforming finite element spaces and, most crucially, newly-devised quasi-interpolation operators delivering optimal estimates on the decay rate of the best-approximation error for functions with Sobolev smoothness index arbitrarily close to 0. The proposed analysis entirely bypasses the technique known in the literature as the discrete compactness argument.}
      \field{issn}{0898-1221}
      \field{journaltitle}{Computers \& Mathematics with Applications}
      \field{month}{2}
      \field{number}{3}
      \field{title}{Analysis of the edge finite element approximation of the {{Maxwell}} equations with low regularity solutions}
      \field{urlday}{28}
      \field{urlmonth}{9}
      \field{urlyear}{2023}
      \field{volume}{75}
      \field{year}{2018}
      \field{dateera}{ce}
      \field{urldateera}{ce}
      \field{pages}{918\bibrangedash 932}
      \range{pages}{15}
      \verb{doi}
      \verb 10.1016/j.camwa.2017.10.017
      \endverb
      \verb{urlraw}
      \verb https://www.sciencedirect.com/science/article/pii/S089812211730665X
      \endverb
      \verb{url}
      \verb https://www.sciencedirect.com/science/article/pii/S089812211730665X
      \endverb
      \keyw{Aubin–Nitsche duality argument,Discrete Poincaré inequality,Edge finite elements,Heterogeneous coefficients,Maxwell equations,Quasi-interpolation}
    \endentry
    \entry{falgoutHypreLibraryHigh2002}{inproceedings}{}
      \name{author}{2}{}{%
        {{hash=84b11e4dfa0d4d81f771416a9e88f1e9}{%
           family={Falgout},
           familyi={F\bibinitperiod},
           given={Robert\bibnamedelima D.},
           giveni={R\bibinitperiod\bibinitdelim D\bibinitperiod}}}%
        {{hash=7a0c922f70357d95bb9b30a9cfdab664}{%
           family={Yang},
           familyi={Y\bibinitperiod},
           given={Ulrike\bibnamedelima Meier},
           giveni={U\bibinitperiod\bibinitdelim M\bibinitperiod}}}%
      }
      \name{editor}{4}{}{%
        {{hash=24951e6a6503cb8a1fec551397605657}{%
           family={Sloot},
           familyi={S\bibinitperiod},
           given={Peter\bibnamedelimb M.\bibnamedelimi A.},
           giveni={P\bibinitperiod\bibinitdelim M\bibinitperiod\bibinitdelim A\bibinitperiod}}}%
        {{hash=4b96ae6720f690c55147125b400f2f02}{%
           family={Hoekstra},
           familyi={H\bibinitperiod},
           given={Alfons\bibnamedelima G.},
           giveni={A\bibinitperiod\bibinitdelim G\bibinitperiod}}}%
        {{hash=70f148d7da4524c4a01e947e0f491db0}{%
           family={Tan},
           familyi={T\bibinitperiod},
           given={C.\bibnamedelimi J.\bibnamedelimi Kenneth},
           giveni={C\bibinitperiod\bibinitdelim J\bibinitperiod\bibinitdelim K\bibinitperiod}}}%
        {{hash=0e7c85e7550c9cf876b4ddf487193e8f}{%
           family={Dongarra},
           familyi={D\bibinitperiod},
           given={Jack\bibnamedelima J.},
           giveni={J\bibinitperiod\bibinitdelim J\bibinitperiod}}}%
      }
      \list{location}{1}{%
        {Berlin, Heidelberg}%
      }
      \list{publisher}{1}{%
        {Springer}%
      }
      \strng{namehash}{aa164e9396ab825faa98b4d2989344d5}
      \strng{fullhash}{aa164e9396ab825faa98b4d2989344d5}
      \strng{bibnamehash}{aa164e9396ab825faa98b4d2989344d5}
      \strng{authorbibnamehash}{aa164e9396ab825faa98b4d2989344d5}
      \strng{authornamehash}{aa164e9396ab825faa98b4d2989344d5}
      \strng{authorfullhash}{aa164e9396ab825faa98b4d2989344d5}
      \strng{editorbibnamehash}{54b44d1f1e6235294ec946600cbb2ae1}
      \strng{editornamehash}{54b44d1f1e6235294ec946600cbb2ae1}
      \strng{editorfullhash}{54b44d1f1e6235294ec946600cbb2ae1}
      \field{labelalpha}{FY02}
      \field{sortinit}{F}
      \field{sortinithash}{fb0c0faa89eb6abae8213bf60e6799ea}
      \field{labelnamesource}{author}
      \field{labeltitlesource}{shorttitle}
      \field{abstract}{hypre is a software library for the solution of large, sparse linear systems on massively parallel computers. Its emphasis is on modern powerful and scalable preconditioners. hypre provides various conceptual interfaces to enable application users to access the library in the way they naturally think about their problems. This paper presents the conceptual interfaces in hypre. An overview of the preconditioners that are available in hypre is given, including some numerical results that show the efficiency of the library.}
      \field{booktitle}{Computational {{Science}} — {{ICCS}} 2002}
      \field{isbn}{978-3-540-47789-1}
      \field{langid}{english}
      \field{series}{Lecture {{Notes}} in {{Computer Science}}}
      \field{shorttitle}{hypre}
      \field{title}{hypre: {{A Library}} of {{High Performance Preconditioners}}}
      \field{year}{2002}
      \field{dateera}{ce}
      \field{pages}{632\bibrangedash 641}
      \range{pages}{10}
      \verb{doi}
      \verb 10.1007/3-540-47789-6_66
      \endverb
      \verb{file}
      \verb /home/rrobin/Zotero/storage/JY6HZQWB/Falgout et Yang - 2002 - hypre A Library of High Performance Preconditione.pdf
      \endverb
    \endentry
    \entry{gernerExistenceOptimalDomains2023}{online}{}
      \name{author}{1}{}{%
        {{hash=db2c48966eab2700ec6b716e50936ee3}{%
           family={Gerner},
           familyi={G\bibinitperiod},
           given={Wadim},
           giveni={W\bibinitperiod}}}%
      }
      \strng{namehash}{db2c48966eab2700ec6b716e50936ee3}
      \strng{fullhash}{db2c48966eab2700ec6b716e50936ee3}
      \strng{bibnamehash}{db2c48966eab2700ec6b716e50936ee3}
      \strng{authorbibnamehash}{db2c48966eab2700ec6b716e50936ee3}
      \strng{authornamehash}{db2c48966eab2700ec6b716e50936ee3}
      \strng{authorfullhash}{db2c48966eab2700ec6b716e50936ee3}
      \field{extraname}{1}
      \field{labelalpha}{Ger23}
      \field{sortinit}{G}
      \field{sortinithash}{62eb2aa29549e4fdbd3cb154ec5711cb}
      \field{extraalpha}{1}
      \field{labelnamesource}{author}
      \field{labeltitlesource}{title}
      \field{abstract}{In the present work we present a general framework which guarantees the existence of optimal domains for isoperimetric problems within the class of \$C\^\{1,1\}\$-regular domains satisfying a uniform ball condition as long as the desired objective function satisfies certain properties. We then verify that the helicity isoperimetric problem studied by Cantarella, DeTurck, Gluck and Teytel in \textbackslash cite\{CDGT002\} satisfies the conditions of our framework and hence establish the existence of optimal domains within the given class of domains. We additionally use the same framework to prove the existence of optimal domains among uniform \$C\^\{1,1\}\$-domains for a first curl eigenvalue problem which has been studied recently for other classes of domains in \textbackslash cite\{EGPS23\}.}
      \field{day}{22}
      \field{eprintclass}{math-ph}
      \field{eprinttype}{arxiv}
      \field{month}{5}
      \field{pubstate}{preprint}
      \field{title}{Existence of optimal domains for the helicity maximisation problem among domains satisfying a uniform ball condition}
      \field{year}{2023}
      \field{dateera}{ce}
      \verb{doi}
      \verb 10.48550/arXiv.2305.13642
      \endverb
      \verb{eprint}
      \verb 2305.13642
      \endverb
      \verb{file}
      \verb /home/rrobin/Zotero/storage/QWYFLXNR/Gerner - 2023 - Existence of optimal domains for the helicity maxi.pdf;/home/rrobin/Zotero/storage/HYYNB759/2305.html
      \endverb
      \verb{urlraw}
      \verb http://arxiv.org/abs/2305.13642
      \endverb
      \verb{url}
      \verb http://arxiv.org/abs/2305.13642
      \endverb
    \endentry
    \entry{gernerIsoperimetricProblemFirst2023}{article}{}
      \name{author}{1}{}{%
        {{hash=db2c48966eab2700ec6b716e50936ee3}{%
           family={Gerner},
           familyi={G\bibinitperiod},
           given={Wadim},
           giveni={W\bibinitperiod}}}%
      }
      \strng{namehash}{db2c48966eab2700ec6b716e50936ee3}
      \strng{fullhash}{db2c48966eab2700ec6b716e50936ee3}
      \strng{bibnamehash}{db2c48966eab2700ec6b716e50936ee3}
      \strng{authorbibnamehash}{db2c48966eab2700ec6b716e50936ee3}
      \strng{authornamehash}{db2c48966eab2700ec6b716e50936ee3}
      \strng{authorfullhash}{db2c48966eab2700ec6b716e50936ee3}
      \field{extraname}{2}
      \field{labelalpha}{Ger23}
      \field{sortinit}{G}
      \field{sortinithash}{62eb2aa29549e4fdbd3cb154ec5711cb}
      \field{extraalpha}{2}
      \field{labelnamesource}{author}
      \field{labeltitlesource}{title}
      \field{abstract}{We consider an isoperimetric problem involving the smallest positive and largest negative curl eigenvalues on abstract ambient manifolds, with a focus on the standard model spaces. We in particular show that the corresponding eigenvalues on optimal domains, assuming optimal domains exist, must be simple in the Euclidean and hyperbolic setting. This generalises a recent result by Enciso and Peralta-Salas who showed the simplicity for axisymmetric optimal domains with connected boundary in Euclidean space. We then generalise another recent result by Enciso and Peralta-Salas, namely that the points of any rotationally symmetric optimal domain with connected boundary in Euclidean space which are closest to the symmetry axis must disconnect the boundary, to the hyperbolic setting, as well as strengthen it in the Euclidean case by removing the connected boundary assumption.}
      \field{issn}{0022-247X}
      \field{journaltitle}{Journal of Mathematical Analysis and Applications}
      \field{month}{3}
      \field{number}{2}
      \field{title}{Isoperimetric problem for the first curl eigenvalue}
      \field{urlday}{25}
      \field{urlmonth}{1}
      \field{urlyear}{2024}
      \field{volume}{519}
      \field{year}{2023}
      \field{urldateera}{ce}
      \field{pages}{126808}
      \range{pages}{1}
      \verb{doi}
      \verb 10.1016/j.jmaa.2022.126808
      \endverb
      \verb{file}
      \verb ScienceDirect Snapshot:C\:\\Users\\rorousse\\Zotero\\storage\\7ZTP4JHL\\S0022247X22008228.html:text/html;Version soumise:C\:\\Users\\rorousse\\Zotero\\storage\\PMWGE97L\\Gerner - 2023 - Isoperimetric problem for the first curl eigenvalu.pdf:application/pdf
      \endverb
      \verb{urlraw}
      \verb https://www.sciencedirect.com/science/article/pii/S0022247X22008228
      \endverb
      \verb{url}
      \verb https://www.sciencedirect.com/science/article/pii/S0022247X22008228
      \endverb
      \keyw{Beltrami fields,Curl operator,Isoperimetric problems,Killing fields,Spectral theory}
    \endentry
    \entry{geuzaineGmsh3DFinite2009}{article}{}
      \name{author}{2}{}{%
        {{hash=5aaa16aa8f1da7cc84b7e5a1aa81f7dc}{%
           family={Geuzaine},
           familyi={G\bibinitperiod},
           given={Christophe},
           giveni={C\bibinitperiod}}}%
        {{hash=a929cc4e99da7a492687d2e3f3c6843a}{%
           family={Remacle},
           familyi={R\bibinitperiod},
           given={Jean-François},
           giveni={J\bibinithyphendelim F\bibinitperiod}}}%
      }
      \strng{namehash}{ed8c8a761b707a4d18fa3b6bf744ce90}
      \strng{fullhash}{ed8c8a761b707a4d18fa3b6bf744ce90}
      \strng{bibnamehash}{ed8c8a761b707a4d18fa3b6bf744ce90}
      \strng{authorbibnamehash}{ed8c8a761b707a4d18fa3b6bf744ce90}
      \strng{authornamehash}{ed8c8a761b707a4d18fa3b6bf744ce90}
      \strng{authorfullhash}{ed8c8a761b707a4d18fa3b6bf744ce90}
      \field{labelalpha}{GR09}
      \field{sortinit}{G}
      \field{sortinithash}{62eb2aa29549e4fdbd3cb154ec5711cb}
      \field{labelnamesource}{author}
      \field{labeltitlesource}{shorttitle}
      \field{abstract}{Gmsh is an open-source 3-D finite element grid generator with a build-in CAD engine and post-processor. Its design goal is to provide a fast, light and user-friendly meshing tool with parametric input and advanced visualization capabilities. This paper presents the overall philosophy, the main design choices and some of the original algorithms implemented in Gmsh. Copyright © 2009 John Wiley \& Sons, Ltd.}
      \field{issn}{1097-0207}
      \field{journaltitle}{International Journal for Numerical Methods in Engineering}
      \field{langid}{english}
      \field{number}{11}
      \field{shorttitle}{Gmsh}
      \field{title}{Gmsh: {{A}} 3-{{D}} finite element mesh generator with built-in pre- and post-processing facilities}
      \field{urlday}{19}
      \field{urlmonth}{10}
      \field{urlyear}{2023}
      \field{volume}{79}
      \field{year}{2009}
      \field{dateera}{ce}
      \field{urldateera}{ce}
      \field{pages}{1309\bibrangedash 1331}
      \range{pages}{23}
      \verb{doi}
      \verb 10.1002/nme.2579
      \endverb
      \verb{file}
      \verb /home/rrobin/Zotero/storage/GUFDB3FW/Geuzaine et Remacle - 2009 - Gmsh A 3-D finite element mesh generator with bui.pdf;/home/rrobin/Zotero/storage/PK6JIDZZ/nme.html
      \endverb
      \verb{urlraw}
      \verb https://onlinelibrary.wiley.com/doi/abs/10.1002/nme.2579
      \endverb
      \verb{url}
      \verb https://onlinelibrary.wiley.com/doi/abs/10.1002/nme.2579
      \endverb
    \endentry
    \entry{henrotShapeVariationOptimization2018}{book}{}
      \name{author}{2}{}{%
        {{hash=21468a2c76602c48d253ae1b628f6c98}{%
           family={Henrot},
           familyi={H\bibinitperiod},
           given={Antoine},
           giveni={A\bibinitperiod}}}%
        {{hash=f03f4e52a698d191cb24ae85e3678607}{%
           family={Pierre},
           familyi={P\bibinitperiod},
           given={Michel},
           giveni={M\bibinitperiod}}}%
      }
      \list{location}{1}{%
        {Zürich, Switzerland}%
      }
      \list{publisher}{1}{%
        {European Math. Soc (EMS)}%
      }
      \strng{namehash}{ba059818f98c421b80c05019275a2755}
      \strng{fullhash}{ba059818f98c421b80c05019275a2755}
      \strng{bibnamehash}{ba059818f98c421b80c05019275a2755}
      \strng{authorbibnamehash}{ba059818f98c421b80c05019275a2755}
      \strng{authornamehash}{ba059818f98c421b80c05019275a2755}
      \strng{authorfullhash}{ba059818f98c421b80c05019275a2755}
      \field{labelalpha}{HP18}
      \field{sortinit}{H}
      \field{sortinithash}{6db6145dae8dc9e1271a8d556090b50a}
      \field{labelnamesource}{author}
      \field{labeltitlesource}{shorttitle}
      \field{isbn}{978-3-03719-178-1}
      \field{langid}{english}
      \field{series}{{{EMS Tracts}} in {{Mathematics}}}
      \field{shorttitle}{Shape {{Variation}} and {{Optimization}}}
      \field{title}{Shape variation and optimization: a geometrical analysis}
      \field{urlday}{28}
      \field{urlmonth}{5}
      \field{urlyear}{2020}
      \field{volume}{28}
      \field{year}{2018}
      \field{dateera}{ce}
      \field{urldateera}{ce}
      \verb{file}
      \verb /home/rrobin/Zotero/storage/2EDF45VY/Henrot et Pierre - 2018 - Shape Variation and Optimization A Geometrical An.pdf;/home/rrobin/Zotero/storage/7EH7V7HD/Henrot et Pierre - 2018 - Shape Variation and Optimization A Geometrical An.pdf;/home/rrobin/Zotero/storage/9ARJXEM9/Henrot et Pierre - 2018 - Shape Variation and Optimization A Geometrical An.pdf;/home/rrobin/Zotero/storage/E3JXB5WQ/Henrot et Pierre - 2018 - Shape Variation and Optimization A Geometrical An.pdf;/home/rrobin/Zotero/storage/XMM6SKAB/Henrot et Pierre - 2018 - Shape Variation and Optimization A Geometrical An.pdf
      \endverb
      \verb{urlraw}
      \verb http://www.ems-ph.org/doi/10.4171/178
      \endverb
      \verb{url}
      \verb http://www.ems-ph.org/doi/10.4171/178
      \endverb
    \endentry
    \entry{hiptmairNodalAuxiliarySpace2007}{article}{}
      \name{author}{2}{}{%
        {{hash=1c516e505db7dc2560aa5487e8c6287d}{%
           family={Hiptmair},
           familyi={H\bibinitperiod},
           given={Ralf},
           giveni={R\bibinitperiod}}}%
        {{hash=f5f125e3a6d4ec2c7143063ea1e885c2}{%
           family={Xu},
           familyi={X\bibinitperiod},
           given={Jinchao},
           giveni={J\bibinitperiod}}}%
      }
      \list{publisher}{1}{%
        {Society for Industrial and Applied Mathematics}%
      }
      \strng{namehash}{2c546709503cd7e50326b1f6dcaea0ff}
      \strng{fullhash}{2c546709503cd7e50326b1f6dcaea0ff}
      \strng{bibnamehash}{2c546709503cd7e50326b1f6dcaea0ff}
      \strng{authorbibnamehash}{2c546709503cd7e50326b1f6dcaea0ff}
      \strng{authornamehash}{2c546709503cd7e50326b1f6dcaea0ff}
      \strng{authorfullhash}{2c546709503cd7e50326b1f6dcaea0ff}
      \field{labelalpha}{HX07}
      \field{sortinit}{H}
      \field{sortinithash}{6db6145dae8dc9e1271a8d556090b50a}
      \field{labelnamesource}{author}
      \field{labeltitlesource}{title}
      \field{abstract}{In this paper we present a family of scalable preconditioners for matrices arising in the discretization of \$H(div)\$ problems using the lowest order Raviart--Thomas finite elements. Our approach belongs to the class of “auxiliary space''--based methods and requires only the finite element stiffness matrix plus some minimal additional discretization information about the topology and orientation of mesh entities. We provide a detailed algebraic description of the theory, parallel implementation, and different variants of this parallel auxiliary space divergence solver (ADS) and discuss its relations to the Hiptmair--Xu (HX) auxiliary space decomposition of \$H(div)\$ [SIAM J. Numer. Anal., 45 (2007), pp. 2483--2509] and to the auxiliary space Maxwell solver AMS [J. Comput. Math., 27 (2009), pp. 604--623]. An extensive set of numerical experiments demonstrates the robustness and scalability of our implementation on large-scale \$H(div)\$ problems with large jumps in the material coefficients.}
      \field{issn}{0036-1429}
      \field{journaltitle}{SIAM Journal on Numerical Analysis}
      \field{month}{1}
      \field{number}{6}
      \field{shortjournal}{SIAM J. Numer. Anal.}
      \field{title}{Nodal {{Auxiliary Space Preconditioning}} in {{H}}(curl) and {{H}}(div) {{Spaces}}}
      \field{urlday}{6}
      \field{urlmonth}{12}
      \field{urlyear}{2023}
      \field{volume}{45}
      \field{year}{2007}
      \field{dateera}{ce}
      \field{urldateera}{ce}
      \field{pages}{2483\bibrangedash 2509}
      \range{pages}{27}
      \verb{doi}
      \verb 10.1137/060660588
      \endverb
      \verb{urlraw}
      \verb https://epubs.siam.org/doi/10.1137/060660588
      \endverb
      \verb{url}
      \verb https://epubs.siam.org/doi/10.1137/060660588
      \endverb
    \endentry
    \entry{imbert-gerardIntroductionSymmetriesStellarators2019}{misc}{}
      \name{author}{3}{}{%
        {{hash=4f425bb98a6b48dc0affc51a40b44c0d}{%
           family={Imbert-Gerard},
           familyi={I\bibinithyphendelim G\bibinitperiod},
           given={Lise-Marie},
           giveni={L\bibinithyphendelim M\bibinitperiod}}}%
        {{hash=29ac33c8d0ab582d7470d3ab660298ae}{%
           family={Paul},
           familyi={P\bibinitperiod},
           given={Elizabeth\bibnamedelima J.},
           giveni={E\bibinitperiod\bibinitdelim J\bibinitperiod}}}%
        {{hash=aa562c67a80d4f479a7d91d30220d3c7}{%
           family={Wright},
           familyi={W\bibinitperiod},
           given={Adelle\bibnamedelima M.},
           giveni={A\bibinitperiod\bibinitdelim M\bibinitperiod}}}%
      }
      \list{publisher}{1}{%
        {arXiv}%
      }
      \strng{namehash}{82e46f59e72ef967091b688cc0308515}
      \strng{fullhash}{82e46f59e72ef967091b688cc0308515}
      \strng{bibnamehash}{82e46f59e72ef967091b688cc0308515}
      \strng{authorbibnamehash}{82e46f59e72ef967091b688cc0308515}
      \strng{authornamehash}{82e46f59e72ef967091b688cc0308515}
      \strng{authorfullhash}{82e46f59e72ef967091b688cc0308515}
      \field{labelalpha}{IPW20}
      \field{sortinit}{I}
      \field{sortinithash}{9417e9a1288a9371e2691d999083ed39}
      \field{labelnamesource}{author}
      \field{labeltitlesource}{shorttitle}
      \field{abstract}{In this self-contained document, we aim to present the basic theoretical building blocks to understand modeling of stellarator magnetic fields, some of the challenges associated with modeling, and optimization for designing stellarators. As often as possible, the ideas will be presented using equations and pictures, and references to other relevant introductory material will be included. This document is accessible to those who may not have a physics background but are interested in applications of mathematical and computational tools to stellarator research.}
      \field{month}{8}
      \field{note}{arXiv:1908.05360 [physics]}
      \field{shorttitle}{An {Introduction} to {Stellarators}}
      \field{title}{An {Introduction} to {Stellarators}: {From} magnetic fields to symmetries and optimization}
      \field{urlday}{25}
      \field{urlmonth}{1}
      \field{urlyear}{2024}
      \field{year}{2020}
      \field{urldateera}{ce}
      \verb{doi}
      \verb 10.48550/arXiv.1908.05360
      \endverb
      \verb{file}
      \verb arXiv Fulltext PDF:C\:\\Users\\rorousse\\Zotero\\storage\\BAAXFJJ5\\Imbert-Gerard et al. - 2020 - An Introduction to Stellarators From magnetic fie.pdf:application/pdf;arXiv.org Snapshot:C\:\\Users\\rorousse\\Zotero\\storage\\H8RPKYF7\\1908.html:text/html
      \endverb
      \verb{urlraw}
      \verb http://arxiv.org/abs/1908.05360
      \endverb
      \verb{url}
      \verb http://arxiv.org/abs/1908.05360
      \endverb
      \keyw{Physics - Plasma Physics}
    \endentry
    \entry{laraSpectralApproximationCurl2015}{article}{}
      \name{author}{3}{}{%
        {{hash=04784f621fe2ab65b6caa6726ffeb833}{%
           family={Lara},
           familyi={L\bibinitperiod},
           given={Eduardo},
           giveni={E\bibinitperiod}}}%
        {{hash=05dffe9382dc827a6bcdb38ceeb5fb60}{%
           family={Rodríguez},
           familyi={R\bibinitperiod},
           given={Rodolfo},
           giveni={R\bibinitperiod}}}%
        {{hash=b8b3ddbdbd6a8163c654be84ccf5570a}{%
           family={Venegas},
           familyi={V\bibinitperiod},
           given={Pablo},
           giveni={P\bibinitperiod}}}%
      }
      \list{publisher}{1}{%
        {Discrete and Continuous Dynamical Systems - S}%
      }
      \strng{namehash}{35d605721fa1d41392431e432b08939d}
      \strng{fullhash}{35d605721fa1d41392431e432b08939d}
      \strng{bibnamehash}{35d605721fa1d41392431e432b08939d}
      \strng{authorbibnamehash}{35d605721fa1d41392431e432b08939d}
      \strng{authornamehash}{35d605721fa1d41392431e432b08939d}
      \strng{authorfullhash}{35d605721fa1d41392431e432b08939d}
      \field{labelalpha}{LRV15}
      \field{sortinit}{L}
      \field{sortinithash}{dad3efd0836470093a7b4a7bb756eb8c}
      \field{labelnamesource}{author}
      \field{labeltitlesource}{title}
      \field{abstract}{A numerical scheme based on Nédélec finite elements has beenrecently introduced to solve the eigenvalue problem for the curloperator in simply connected domains. This topological assumption is notjust a technicality, since the eigenvalue problem is ill-posed onmultiply connected domains, in the sense that its spectrum is the wholecomplex plane. However, additional constraints can be added to theeigenvalue problem in order to recover a well-posed problem with adiscrete spectrum. Vanishing circulations on each non-bounding cycle inthe complement of the domain have been chosen as additional constraintsin this paper. A mixed weak formulation including a Lagrange multiplier(that turns out to vanish) is introduced and shown to be well-posed.This formulation is discretized by Nédélec elements, while standardfinite elements are used for the Lagrange multiplier. Spectralconvergence is proved as well as a priori error estimates. It is alsoshown how to implement this finite element discretization taking care ofthese additional constraints. Finally, a numerical test to assess theperformance of the proposed methods is reported.}
      \field{issn}{1937-1632}
      \field{journaltitle}{Discrete and Continuous Dynamical Systems - S}
      \field{langid}{english}
      \field{number}{1}
      \field{shortjournal}{DCDS-S}
      \field{title}{Spectral approximation of the curl operator in multiply connected domains}
      \field{urlday}{4}
      \field{urlmonth}{12}
      \field{urlyear}{2023}
      \field{volume}{9}
      \field{year}{Mon Nov 30 19:00:00 EST 2015}
      \field{urldateera}{ce}
      \field{pages}{235\bibrangedash 253}
      \range{pages}{19}
      \verb{doi}
      \verb 10.3934/dcdss.2016.9.235
      \endverb
      \verb{file}
      \verb /home/rrobin/Zotero/storage/3ZSXULDS/Lara et al. - 2015 - Spectral approximation of the curl operator in mul.pdf
      \endverb
      \verb{urlraw}
      \verb https://www.aimsciences.org/en/article/doi/10.3934/dcdss.2016.9.235
      \endverb
      \verb{url}
      \verb https://www.aimsciences.org/en/article/doi/10.3934/dcdss.2016.9.235
      \endverb
    \endentry
    \entry{mactaggartMagneticHelicityMultiply2019}{article}{}
      \name{author}{2}{}{%
        {{hash=188b25d5b6066b6c526509cd994896b6}{%
           family={MacTaggart},
           familyi={M\bibinitperiod},
           given={D.},
           giveni={D\bibinitperiod}}}%
        {{hash=31171371a77aefab1503a9ea99f99271}{%
           family={Valli},
           familyi={V\bibinitperiod},
           given={A.},
           giveni={A\bibinitperiod}}}%
      }
      \strng{namehash}{60ce05db265051a65458fabee434725a}
      \strng{fullhash}{60ce05db265051a65458fabee434725a}
      \strng{bibnamehash}{60ce05db265051a65458fabee434725a}
      \strng{authorbibnamehash}{60ce05db265051a65458fabee434725a}
      \strng{authornamehash}{60ce05db265051a65458fabee434725a}
      \strng{authorfullhash}{60ce05db265051a65458fabee434725a}
      \field{labelalpha}{MV19}
      \field{sortinit}{M}
      \field{sortinithash}{2e5c2f51f7fa2d957f3206819bf86dc3}
      \field{labelnamesource}{author}
      \field{labeltitlesource}{title}
      \field{abstract}{Magnetic helicity is a fundamental quantity of magnetohydrodynamics that carries topological information about the magnetic field. By ‘topological information’, we usually refer to the linkage of magnetic field lines. For domains that are not simply connected, however, helicity also depends on the topology of the domain. In this paper we expand the standard definition of magnetic helicity in simply connected domains to multiply connected domains in R3R3{\textbackslash}mathbb\{R\}{\textasciicircum}\{3\} of arbitrary topology. We also discuss how using the classic Biot–Savart operator simplifies the expression for helicity and how domain topology affects the physical interpretation of helicity.}
      \field{issn}{0022-3778, 1469-7807}
      \field{journaltitle}{Journal of Plasma Physics}
      \field{month}{10}
      \field{note}{Publisher: Cambridge University Press}
      \field{number}{5}
      \field{title}{Magnetic helicity in multiply connected domains}
      \field{urlday}{25}
      \field{urlmonth}{1}
      \field{urlyear}{2024}
      \field{volume}{85}
      \field{year}{2019}
      \field{urldateera}{ce}
      \field{pages}{775850501}
      \range{pages}{1}
      \verb{doi}
      \verb 10.1017/S0022377819000576
      \endverb
      \verb{file}
      \verb Full Text PDF:C\:\\Users\\rorousse\\Zotero\\storage\\GN5JVYBL\\MacTaggart et Valli - 2019 - Magnetic helicity in multiply connected domains.pdf:application/pdf
      \endverb
      \verb{urlraw}
      \verb https://www.cambridge.org/core/journals/journal-of-plasma-physics/article/magnetic-helicity-in-multiply-connected-domains/E0F6B7669EC03884F0C00F41B61F4356
      \endverb
      \verb{url}
      \verb https://www.cambridge.org/core/journals/journal-of-plasma-physics/article/magnetic-helicity-in-multiply-connected-domains/E0F6B7669EC03884F0C00F41B61F4356
      \endverb
      \keyw{plasma nonlinear phenomena,plasma properties}
    \endentry
    \entry{moffattDegreeKnottednessTangled1969}{article}{}
      \name{author}{1}{}{%
        {{hash=620879a8ac85084f3886e217ca5dff19}{%
           family={Moffatt},
           familyi={M\bibinitperiod},
           given={H.\bibnamedelimi K.},
           giveni={H\bibinitperiod\bibinitdelim K\bibinitperiod}}}%
      }
      \list{publisher}{1}{%
        {Cambridge University Press}%
      }
      \strng{namehash}{620879a8ac85084f3886e217ca5dff19}
      \strng{fullhash}{620879a8ac85084f3886e217ca5dff19}
      \strng{bibnamehash}{620879a8ac85084f3886e217ca5dff19}
      \strng{authorbibnamehash}{620879a8ac85084f3886e217ca5dff19}
      \strng{authornamehash}{620879a8ac85084f3886e217ca5dff19}
      \strng{authorfullhash}{620879a8ac85084f3886e217ca5dff19}
      \field{labelalpha}{Mof69}
      \field{sortinit}{M}
      \field{sortinithash}{2e5c2f51f7fa2d957f3206819bf86dc3}
      \field{labelnamesource}{author}
      \field{labeltitlesource}{title}
      \field{abstract}{Let u(x) be the velocity field in a fluid of infinite extent due to a vorticity distribution w(x) which is zero except in two closed vortex filaments of strengths K1, K2. It is first shown that the integral \textbackslash [ I=\textbackslash int\{\textbackslash bf u\}.\{\textbackslash boldmath \textbackslash omega\}\textbackslash,dV \textbackslash ] is equal to αK1K2 where α is an integer representing the degree of linkage of the two filaments; α = 0 if they are unlinked, ± 1 if they are singly linked. The invariance of I for a continuous localized vorticity distribution is then established for barotropic inviscid flow under conservative body forces. The result is interpreted in terms of the conservation of linkages of vortex lines which move with the fluid.Some examples of steady flows for which I ± 0 are briefly described; in particular, attention is drawn to a family of spherical vortices with swirl (which is closely analogous to a known family of solutions of the equations of magnetostatics); the vortex lines of these flows are both knotted and linked.Two related magnetohydrodynamic invariants discovered by Woltjer (1958a, b) are discussed in ±5.}
      \field{issn}{1469-7645, 0022-1120}
      \field{journaltitle}{Journal of Fluid Mechanics}
      \field{langid}{english}
      \field{month}{1}
      \field{number}{1}
      \field{title}{The degree of knottedness of tangled vortex lines}
      \field{urlday}{21}
      \field{urlmonth}{11}
      \field{urlyear}{2023}
      \field{volume}{35}
      \field{year}{1969}
      \field{dateera}{ce}
      \field{urldateera}{ce}
      \field{pages}{117\bibrangedash 129}
      \range{pages}{13}
      \verb{doi}
      \verb 10.1017/S0022112069000991
      \endverb
      \verb{file}
      \verb /home/rrobin/Zotero/storage/83KXDQSV/Moffatt - 1969 - The degree of knottedness of tangled vortex lines.pdf
      \endverb
      \verb{urlraw}
      \verb https://www.cambridge.org/core/journals/journal-of-fluid-mechanics/article/degree-of-knottedness-of-tangled-vortex-lines/798DFB897A75D97CE3081F4C1DA970F0
      \endverb
      \verb{url}
      \verb https://www.cambridge.org/core/journals/journal-of-fluid-mechanics/article/degree-of-knottedness-of-tangled-vortex-lines/798DFB897A75D97CE3081F4C1DA970F0
      \endverb
    \endentry
    \entry{montielIsoperimetricProblemCurl2023}{online}{}
      \name{author}{1}{}{%
        {{hash=3d2cd26333164d8ed826caafb2cce7ad}{%
           family={Montiel},
           familyi={M\bibinitperiod},
           given={S.},
           giveni={S\bibinitperiod}}}%
      }
      \strng{namehash}{3d2cd26333164d8ed826caafb2cce7ad}
      \strng{fullhash}{3d2cd26333164d8ed826caafb2cce7ad}
      \strng{bibnamehash}{3d2cd26333164d8ed826caafb2cce7ad}
      \strng{authorbibnamehash}{3d2cd26333164d8ed826caafb2cce7ad}
      \strng{authornamehash}{3d2cd26333164d8ed826caafb2cce7ad}
      \strng{authorfullhash}{3d2cd26333164d8ed826caafb2cce7ad}
      \field{labelalpha}{Mon23}
      \field{sortinit}{M}
      \field{sortinithash}{2e5c2f51f7fa2d957f3206819bf86dc3}
      \field{labelnamesource}{author}
      \field{labeltitlesource}{title}
      \field{abstract}{In the last decades, many mathematicians have studied the curl operator in compact three-manifolds , mainly the structure of its spectrum and some isoperimetric problems associated with it. In this paper, we will see that all the compact threemanifolds (both closed and and with non-empty boundary) have always optimal lower bounds for the absolute value of their nonnull eigenvalues of curl. We will also show that these bounds are always attained and compute the optimal domains and the multiplicities of their associated eigenvalues. So, we have solved the isoperimetric problem associated to the curl operator, and, by the way, we have solved and old Cantarella, de Turck, Gluck and Teytel conjecture [CdTGT] in the negative.}
      \field{day}{18}
      \field{eprintclass}{math}
      \field{eprinttype}{arxiv}
      \field{langid}{english}
      \field{month}{7}
      \field{pubstate}{preprint}
      \field{title}{The {{Isoperimetric Problem}} for the {{Curl Operator}}}
      \field{year}{2023}
      \field{dateera}{ce}
      \verb{eprint}
      \verb 2307.09556
      \endverb
      \verb{file}
      \verb /home/rrobin/Zotero/storage/3KQ3B8P6/Montiel - 2023 - The Isoperimetric Problem for the Curl Operator.pdf
      \endverb
      \verb{urlraw}
      \verb http://arxiv.org/abs/2307.09556
      \endverb
      \verb{url}
      \verb http://arxiv.org/abs/2307.09556
      \endverb
    \endentry
    \entry{nedelecMixedFiniteElements1980}{article}{}
      \name{author}{1}{}{%
        {{hash=d2d06e5706bd5e14fc070df88043d97c}{%
           family={Nedelec},
           familyi={N\bibinitperiod},
           given={J.\bibnamedelimi C.},
           giveni={J\bibinitperiod\bibinitdelim C\bibinitperiod}}}%
      }
      \strng{namehash}{d2d06e5706bd5e14fc070df88043d97c}
      \strng{fullhash}{d2d06e5706bd5e14fc070df88043d97c}
      \strng{bibnamehash}{d2d06e5706bd5e14fc070df88043d97c}
      \strng{authorbibnamehash}{d2d06e5706bd5e14fc070df88043d97c}
      \strng{authornamehash}{d2d06e5706bd5e14fc070df88043d97c}
      \strng{authorfullhash}{d2d06e5706bd5e14fc070df88043d97c}
      \field{labelalpha}{Ned80}
      \field{sortinit}{N}
      \field{sortinithash}{98cf339a479c0454fe09153a08675a15}
      \field{labelnamesource}{author}
      \field{labeltitlesource}{title}
      \field{abstract}{We present here some new families of non conforming finite elements in R3. These two families of finite elements, built on tetrahedrons or on cubes are respectively conforming in the spacesH(curl) andH(div). We give some applications of these elements for the approximation of Maxwell's equations and equations of elasticity.}
      \field{day}{1}
      \field{issn}{0945-3245}
      \field{journaltitle}{Numerische Mathematik}
      \field{langid}{english}
      \field{month}{9}
      \field{number}{3}
      \field{shortjournal}{Numer. Math.}
      \field{title}{Mixed finite elements in {{R3}}}
      \field{urlday}{16}
      \field{urlmonth}{11}
      \field{urlyear}{2023}
      \field{volume}{35}
      \field{year}{1980}
      \field{dateera}{ce}
      \field{urldateera}{ce}
      \field{pages}{315\bibrangedash 341}
      \range{pages}{27}
      \verb{doi}
      \verb 10.1007/BF01396415
      \endverb
      \verb{urlraw}
      \verb https://doi.org/10.1007/BF01396415
      \endverb
      \verb{url}
      \verb https://doi.org/10.1007/BF01396415
      \endverb
    \endentry
    \entry{paulAdjointMethodGradientbased2018}{article}{}
      \name{author}{4}{}{%
        {{hash=bf74568e4e07c534fac1a3365c14f80b}{%
           family={Paul},
           familyi={P\bibinitperiod},
           given={E.\bibnamedelimi J.},
           giveni={E\bibinitperiod\bibinitdelim J\bibinitperiod}}}%
        {{hash=2038c21f603629083df10705711967f2}{%
           family={Landreman},
           familyi={L\bibinitperiod},
           given={M.},
           giveni={M\bibinitperiod}}}%
        {{hash=4fd438bc04d4cbb5749e32a89f2c8bb3}{%
           family={Bader},
           familyi={B\bibinitperiod},
           given={A.},
           giveni={A\bibinitperiod}}}%
        {{hash=3760142ae0521badeb659b4359e5b14e}{%
           family={Dorland},
           familyi={D\bibinitperiod},
           given={W.},
           giveni={W\bibinitperiod}}}%
      }
      \list{publisher}{1}{%
        {IOP Publishing}%
      }
      \strng{namehash}{eb43979c66d3268bfbade5ffb348b1cc}
      \strng{fullhash}{eb43979c66d3268bfbade5ffb348b1cc}
      \strng{bibnamehash}{eb43979c66d3268bfbade5ffb348b1cc}
      \strng{authorbibnamehash}{eb43979c66d3268bfbade5ffb348b1cc}
      \strng{authornamehash}{eb43979c66d3268bfbade5ffb348b1cc}
      \strng{authorfullhash}{eb43979c66d3268bfbade5ffb348b1cc}
      \field{labelalpha}{Pau+18}
      \field{sortinit}{P}
      \field{sortinithash}{bb5b15f2db90f7aef79bb9e83defefcb}
      \field{labelnamesource}{author}
      \field{labeltitlesource}{title}
      \field{abstract}{We present a method for stellarator coil design via gradient-based optimization of the coil-winding surface. The REGCOIL (Landreman 2017 Nucl. Fusion 57 046003) approach is used to obtain the coil shapes on the winding surface using a continuous current potential. We apply the adjoint method to calculate derivatives of the objective function, allowing for efficient computation of analytic gradients while eliminating the numerical noise of approximate derivatives. We are able to improve engineering properties of the coils by targeting the root-mean-squared current density in the objective function. We obtain winding surfaces for W7-X and HSX which simultaneously decrease the normal magnetic field on the plasma surface and increase the surface-averaged distance between the coils and the plasma in comparison with the actual winding surfaces. The coils computed on the optimized surfaces feature a smaller toroidal extent and curvature and increased inter-coil spacing. A technique for computation of the local sensitivity of figures of merit to normal displacements of the winding surface is presented, with potential applications for understanding engineering tolerances.}
      \field{issn}{0029-5515}
      \field{journaltitle}{Nuclear Fusion}
      \field{langid}{english}
      \field{month}{5}
      \field{number}{7}
      \field{shortjournal}{Nucl. Fusion}
      \field{title}{An adjoint method for gradient-based optimization of stellarator coil shapes}
      \field{urlday}{10}
      \field{urlmonth}{6}
      \field{urlyear}{2022}
      \field{volume}{58}
      \field{year}{2018}
      \field{dateera}{ce}
      \field{urldateera}{ce}
      \field{pages}{076015}
      \range{pages}{1}
      \verb{doi}
      \verb 10.1088/1741-4326/aac1c7
      \endverb
      \verb{file}
      \verb /home/rrobin/Zotero/storage/F46KC4UA/Paul et al. - 2018 - An adjoint method for gradient-based optimization .pdf
      \endverb
      \verb{urlraw}
      \verb https://doi.org/10.1088/1741-4326/aac1c7
      \endverb
      \verb{url}
      \verb https://doi.org/10.1088/1741-4326/aac1c7
      \endverb
    \endentry
    \entry{privatOptimalShapeStellarators2022}{article}{}
      \name{author}{3}{}{%
        {{hash=5c28b65eed681121f5cd0d4eadb4fd3e}{%
           family={Privat},
           familyi={P\bibinitperiod},
           given={Yannick},
           giveni={Y\bibinitperiod}}}%
        {{hash=9a2f0cf98c97c50a1f7975672b9a15ee}{%
           family={Robin},
           familyi={R\bibinitperiod},
           given={Rémi},
           giveni={R\bibinitperiod}}}%
        {{hash=b4e4ea11e20c124f5c380bd847d20538}{%
           family={Sigalotti},
           familyi={S\bibinitperiod},
           given={Mario},
           giveni={M\bibinitperiod}}}%
      }
      \strng{namehash}{2b182d0d840beb2a065e3ad0c0d9642e}
      \strng{fullhash}{2b182d0d840beb2a065e3ad0c0d9642e}
      \strng{bibnamehash}{2b182d0d840beb2a065e3ad0c0d9642e}
      \strng{authorbibnamehash}{2b182d0d840beb2a065e3ad0c0d9642e}
      \strng{authornamehash}{2b182d0d840beb2a065e3ad0c0d9642e}
      \strng{authorfullhash}{2b182d0d840beb2a065e3ad0c0d9642e}
      \field{extraname}{1}
      \field{labelalpha}{PRS22}
      \field{sortinit}{P}
      \field{sortinithash}{bb5b15f2db90f7aef79bb9e83defefcb}
      \field{labelnamesource}{author}
      \field{labeltitlesource}{title}
      \field{abstract}{We are interested in the design of stellarators, devices for the production of controlled nuclear fusion reactions alternative to tokamaks. The confinement of the plasma is entirely achieved by a helical magnetic field created by the complex arrangement of coils fed by high currents around a toroidal domain. Such coils describe a surface called “coil winding surface” (CWS). In this paper, we model the design of the CWS as a shape optimization problem, so that the cost functional reflects both optimal plasma confinement properties, through a least square discrepancy, and also manufacturability, thanks to geometrical terms involving the lateral surface or the curvature of the CWS. We completely analyze the resulting problem: on the one hand, we establish the existence of an optimal shape, prove the shape differentiability of the criterion, and provide the expression of the differential in a workable form. On the other hand, we propose a numerical method and perform simulations of optimal stellarator shapes. We discuss the efficiency of our approach with respect to the literature in this area. Résumé Nous nous intéressons à la conception de stellarators, qui sont une technologie de réacteurs à fusion thermonucléaire contrôlée alternative aux tokamaks. Le confinement du plasma est entièrement réalisé par un champ magnétique hélicoïdal créé par l'agencement complexe de bobines alimentées par de forts courants autour d'un domaine toroïdal. De telles bobines décrivent une surface appelée “coil winding surface” (CWS). Dans cet article, nous modélisons la conception de la CWS comme un problème d'optimisation de forme, de sorte que la fonctionnelle de coût reflète à la fois les propriétés optimales de confinement du plasma, à travers une pénalisation des moindres carrés, et également la manufacturabilité, grâce à des termes géométriques impliquant la surface latérale ou la courbure de la CWS. Nous analysons complètement le problème résultant : d'une part, nous établissons l'existence d'une forme optimale, prouvons la différentiabilité de forme du critère et fournissons l'expression de la différentielle sous une forme exploitable. D'autre part, nous proposons une méthode numérique et effectuons des simulations de formes optimales de stellarator. Nous discutons de l'efficacité de notre approche par rapport à la littérature dans ce domaine.}
      \field{issn}{0021-7824}
      \field{journaltitle}{Journal de Mathématiques Pures et Appliquées}
      \field{langid}{english}
      \field{month}{7}
      \field{shortjournal}{J. Math. Pures Appl.}
      \field{title}{Optimal shape of stellarators for magnetic confinement fusion}
      \field{urlday}{14}
      \field{urlmonth}{6}
      \field{urlyear}{2022}
      \field{volume}{163}
      \field{year}{2022}
      \field{dateera}{ce}
      \field{urldateera}{ce}
      \field{pages}{231\bibrangedash 264}
      \range{pages}{34}
      \verb{doi}
      \verb 10.1016/j.matpur.2022.05.005
      \endverb
      \verb{file}
      \verb /home/rrobin/Zotero/storage/J7TYRQH3/Privat et al. - 2022 - Optimal shape of stellarators for magnetic confine.pdf;/home/rrobin/Zotero/storage/S3NLXEX5/S0021782422000587.html
      \endverb
      \verb{urlraw}
      \verb https://hal.inria.fr/hal-03472623
      \endverb
      \verb{url}
      \verb https://hal.inria.fr/hal-03472623
      \endverb
    \endentry
    \entry{privatExistenceSurfacesOptimizing2022}{article}{}
      \name{author}{3}{}{%
        {{hash=5c28b65eed681121f5cd0d4eadb4fd3e}{%
           family={Privat},
           familyi={P\bibinitperiod},
           given={Yannick},
           giveni={Y\bibinitperiod}}}%
        {{hash=9a2f0cf98c97c50a1f7975672b9a15ee}{%
           family={Robin},
           familyi={R\bibinitperiod},
           given={Rémi},
           giveni={R\bibinitperiod}}}%
        {{hash=b4e4ea11e20c124f5c380bd847d20538}{%
           family={Sigalotti},
           familyi={S\bibinitperiod},
           given={Mario},
           giveni={M\bibinitperiod}}}%
      }
      \strng{namehash}{2b182d0d840beb2a065e3ad0c0d9642e}
      \strng{fullhash}{2b182d0d840beb2a065e3ad0c0d9642e}
      \strng{bibnamehash}{2b182d0d840beb2a065e3ad0c0d9642e}
      \strng{authorbibnamehash}{2b182d0d840beb2a065e3ad0c0d9642e}
      \strng{authornamehash}{2b182d0d840beb2a065e3ad0c0d9642e}
      \strng{authorfullhash}{2b182d0d840beb2a065e3ad0c0d9642e}
      \field{extraname}{2}
      \field{labelalpha}{PRS24}
      \field{sortinit}{P}
      \field{sortinithash}{bb5b15f2db90f7aef79bb9e83defefcb}
      \field{labelnamesource}{author}
      \field{labeltitlesource}{title}
      \field{abstract}{Yannick Privat, Rémi Robin, Mario Sigalotti}
      \field{issn}{1463-9963}
      \field{journaltitle}{Interfaces and Free Boundaries}
      \field{month}{6}
      \field{title}{Existence of surfaces optimizing geometric and {PDE} shape functionals under reach constraint}
      \field{urlday}{25}
      \field{urlmonth}{7}
      \field{urlyear}{2024}
      \field{year}{2024}
      \field{urldateera}{ce}
      \verb{doi}
      \verb 10.4171/ifb/523
      \endverb
      \verb{urlraw}
      \verb https://ems.press/journals/ifb/articles/14297927
      \endverb
      \verb{url}
      \verb https://ems.press/journals/ifb/articles/14297927
      \endverb
    \endentry
    \entry{raviartMixedFiniteElement1977}{inproceedings}{}
      \name{author}{2}{}{%
        {{hash=c4ebe38117cda0b142689e8a84abb34e}{%
           family={Raviart},
           familyi={R\bibinitperiod},
           given={P.\bibnamedelimi A.},
           giveni={P\bibinitperiod\bibinitdelim A\bibinitperiod}}}%
        {{hash=a618f65af74dcf6ee4ecbef5e5f2edf5}{%
           family={Thomas},
           familyi={T\bibinitperiod},
           given={J.\bibnamedelimi M.},
           giveni={J\bibinitperiod\bibinitdelim M\bibinitperiod}}}%
      }
      \name{editor}{2}{}{%
        {{hash=435b04704c22d97a4b3c640039eee281}{%
           family={Galligani},
           familyi={G\bibinitperiod},
           given={Ilio},
           giveni={I\bibinitperiod}}}%
        {{hash=93e7fefc0e07db67f3bf58f20f48d6ec}{%
           family={Magenes},
           familyi={M\bibinitperiod},
           given={Enrico},
           giveni={E\bibinitperiod}}}%
      }
      \list{location}{1}{%
        {Berlin, Heidelberg}%
      }
      \list{publisher}{1}{%
        {Springer}%
      }
      \strng{namehash}{3c0d5aad48dfee73bd17db311f9c9149}
      \strng{fullhash}{3c0d5aad48dfee73bd17db311f9c9149}
      \strng{bibnamehash}{3c0d5aad48dfee73bd17db311f9c9149}
      \strng{authorbibnamehash}{3c0d5aad48dfee73bd17db311f9c9149}
      \strng{authornamehash}{3c0d5aad48dfee73bd17db311f9c9149}
      \strng{authorfullhash}{3c0d5aad48dfee73bd17db311f9c9149}
      \strng{editorbibnamehash}{701e1f3974203a47a341a4b60a0b8cfe}
      \strng{editornamehash}{701e1f3974203a47a341a4b60a0b8cfe}
      \strng{editorfullhash}{701e1f3974203a47a341a4b60a0b8cfe}
      \field{labelalpha}{RT77}
      \field{sortinit}{R}
      \field{sortinithash}{b9c68a358aea118dfa887b6e902414a7}
      \field{labelnamesource}{author}
      \field{labeltitlesource}{title}
      \field{booktitle}{Mathematical {{Aspects}} of {{Finite Element Methods}}}
      \field{isbn}{978-3-540-37158-8}
      \field{langid}{english}
      \field{series}{Lecture {{Notes}} in {{Mathematics}}}
      \field{title}{A mixed finite element method for 2-nd order elliptic problems}
      \field{year}{1977}
      \field{dateera}{ce}
      \field{pages}{292\bibrangedash 315}
      \range{pages}{24}
      \verb{doi}
      \verb 10.1007/BFb0064470
      \endverb
      \verb{file}
      \verb /home/rrobin/Zotero/storage/83KE3ZBK/Raviart et Thomas - 1977 - A mixed finite element method for 2-nd order ellip.pdf
      \endverb
    \endentry
    \entry{saadGMRESGeneralizedMinimal1986}{article}{}
      \name{author}{2}{}{%
        {{hash=2e9e890ef0ef7a6569e8afb4da0fbf7a}{%
           family={Saad},
           familyi={S\bibinitperiod},
           given={Youcef},
           giveni={Y\bibinitperiod}}}%
        {{hash=eac3074c8353b52f159e50629f58cfe2}{%
           family={Schultz},
           familyi={S\bibinitperiod},
           given={Martin\bibnamedelima H.},
           giveni={M\bibinitperiod\bibinitdelim H\bibinitperiod}}}%
      }
      \strng{namehash}{e503d866277c9f7effb481e6ffb692d6}
      \strng{fullhash}{e503d866277c9f7effb481e6ffb692d6}
      \strng{bibnamehash}{e503d866277c9f7effb481e6ffb692d6}
      \strng{authorbibnamehash}{e503d866277c9f7effb481e6ffb692d6}
      \strng{authornamehash}{e503d866277c9f7effb481e6ffb692d6}
      \strng{authorfullhash}{e503d866277c9f7effb481e6ffb692d6}
      \field{labelalpha}{SS86}
      \field{sortinit}{S}
      \field{sortinithash}{c319cff79d99c853d775f88277d4e45f}
      \field{labelnamesource}{author}
      \field{labeltitlesource}{title}
      \field{abstract}{We present an iterative method for solving linear systems, which has the property of minimizing at every step the norm of the residual vector over a Krylov subspace. The algorithm is derived from the Arnoldi process for constructing an \$l\_2 \$-orthogonal basis of Krylov subspaces. It can be considered as a generalization of Paige and Saunders’ MINRES algorithm and is theoretically equivalent to the Generalized Conjugate Residual (GCR) method and to ORTHODIR. The new algorithm presents several advantages over GCR and ORTHODIR.}
      \field{journaltitle}{SIAM Journal on Scientific and Statistical Computing}
      \field{number}{3}
      \field{title}{{{GMRES}}: {{A}} generalized minimal residual algorithm for solving nonsymmetric linear systems}
      \field{volume}{7}
      \field{year}{1986}
      \field{dateera}{ce}
      \field{pages}{856\bibrangedash 869}
      \range{pages}{14}
      \verb{doi}
      \verb 10.1137/0907058
      \endverb
      \verb{eprint}
      \verb https://doi.org/10.1137/0907058
      \endverb
      \verb{urlraw}
      \verb https://doi.org/10.1137/0907058
      \endverb
      \verb{url}
      \verb https://doi.org/10.1137/0907058
      \endverb
    \endentry
    \entry{scroggsBasixRuntimeFinite2022}{article}{}
      \name{author}{4}{}{%
        {{hash=aadaec59c647f6faf392faf3d082f623}{%
           family={Scroggs},
           familyi={S\bibinitperiod},
           given={Matthew\bibnamedelima W.},
           giveni={M\bibinitperiod\bibinitdelim W\bibinitperiod}}}%
        {{hash=f83eb8370a9602804f9cd11e05bd33ce}{%
           family={Baratta},
           familyi={B\bibinitperiod},
           given={Igor\bibnamedelima A.},
           giveni={I\bibinitperiod\bibinitdelim A\bibinitperiod}}}%
        {{hash=721b7a4219213ef1a2a17d867236ca35}{%
           family={Richardson},
           familyi={R\bibinitperiod},
           given={Chris\bibnamedelima N.},
           giveni={C\bibinitperiod\bibinitdelim N\bibinitperiod}}}%
        {{hash=997effee24e95b11152025bbad8e00ad}{%
           family={Wells},
           familyi={W\bibinitperiod},
           given={Garth\bibnamedelima N.},
           giveni={G\bibinitperiod\bibinitdelim N\bibinitperiod}}}%
      }
      \strng{namehash}{42ccbb65eff01e181f5f09975654da2e}
      \strng{fullhash}{42ccbb65eff01e181f5f09975654da2e}
      \strng{bibnamehash}{42ccbb65eff01e181f5f09975654da2e}
      \strng{authorbibnamehash}{42ccbb65eff01e181f5f09975654da2e}
      \strng{authornamehash}{42ccbb65eff01e181f5f09975654da2e}
      \strng{authorfullhash}{42ccbb65eff01e181f5f09975654da2e}
      \field{labelalpha}{Scr+22}
      \field{sortinit}{S}
      \field{sortinithash}{c319cff79d99c853d775f88277d4e45f}
      \field{extraalpha}{1}
      \field{labelnamesource}{author}
      \field{labeltitlesource}{shorttitle}
      \field{abstract}{Scroggs et al., (2022). Basix: a runtime finite element basis evaluation library. Journal of Open Source Software, 7(73), 3982, https://doi.org/10.21105/joss.03982}
      \field{day}{25}
      \field{issn}{2475-9066}
      \field{journaltitle}{Journal of Open Source Software}
      \field{langid}{english}
      \field{month}{5}
      \field{number}{73}
      \field{shorttitle}{Basix}
      \field{title}{Basix: a runtime finite element basis evaluation library}
      \field{urlday}{24}
      \field{urlmonth}{10}
      \field{urlyear}{2023}
      \field{volume}{7}
      \field{year}{2022}
      \field{dateera}{ce}
      \field{urldateera}{ce}
      \field{pages}{3982}
      \range{pages}{1}
      \verb{doi}
      \verb 10.21105/joss.03982
      \endverb
      \verb{file}
      \verb /home/rrobin/Zotero/storage/WW5D3N8N/Scroggs et al. - 2022 - Basix a runtime finite element basis evaluation l.pdf
      \endverb
      \verb{urlraw}
      \verb https://joss.theoj.org/papers/10.21105/joss.03982
      \endverb
      \verb{url}
      \verb https://joss.theoj.org/papers/10.21105/joss.03982
      \endverb
    \endentry
    \entry{scroggsConstructionArbitraryOrder2022}{article}{}
      \name{author}{4}{}{%
        {{hash=aadaec59c647f6faf392faf3d082f623}{%
           family={Scroggs},
           familyi={S\bibinitperiod},
           given={Matthew\bibnamedelima W.},
           giveni={M\bibinitperiod\bibinitdelim W\bibinitperiod}}}%
        {{hash=ef2c1d776da6bc66d1a876af3804aff6}{%
           family={Dokken},
           familyi={D\bibinitperiod},
           given={Jørgen\bibnamedelima S.},
           giveni={J\bibinitperiod\bibinitdelim S\bibinitperiod}}}%
        {{hash=721b7a4219213ef1a2a17d867236ca35}{%
           family={Richardson},
           familyi={R\bibinitperiod},
           given={Chris\bibnamedelima N.},
           giveni={C\bibinitperiod\bibinitdelim N\bibinitperiod}}}%
        {{hash=997effee24e95b11152025bbad8e00ad}{%
           family={Wells},
           familyi={W\bibinitperiod},
           given={Garth\bibnamedelima N.},
           giveni={G\bibinitperiod\bibinitdelim N\bibinitperiod}}}%
      }
      \strng{namehash}{c83538a7c73a2d00c38229fefe2cfa12}
      \strng{fullhash}{c83538a7c73a2d00c38229fefe2cfa12}
      \strng{bibnamehash}{c83538a7c73a2d00c38229fefe2cfa12}
      \strng{authorbibnamehash}{c83538a7c73a2d00c38229fefe2cfa12}
      \strng{authornamehash}{c83538a7c73a2d00c38229fefe2cfa12}
      \strng{authorfullhash}{c83538a7c73a2d00c38229fefe2cfa12}
      \field{labelalpha}{Scr+22}
      \field{sortinit}{S}
      \field{sortinithash}{c319cff79d99c853d775f88277d4e45f}
      \field{extraalpha}{2}
      \field{labelnamesource}{author}
      \field{labeltitlesource}{title}
      \field{abstract}{We develop a method for generating degree-of-freedom maps for arbitrary order Ciarlet-type finite element spaces for any cell shape. The approach is based on the composition of permutations and transformations by cell sub-entity. Current approaches to generating degree-of-freedom maps for arbitrary order problems typically rely on a consistent orientation of cell entities that permits the definition of a common local coordinate system on shared edges and faces. However, while orientation of a mesh is straightforward for simplex cells and is a local operation, it is not a strictly local operation for quadrilateral cells and, in the case of hexahedral cells, not all meshes are orientable. The permutation and transformation approach is developed for a range of element types, including arbitrary degree Lagrange, serendipity, and divergence- and curl-conforming elements, and for a range of cell shapes. The approach is local and can be applied to cells of any shape, including general polytopes and meshes with mixed cell types. A number of examples are presented and the developed approach has been implemented in open-source libraries.}
      \field{day}{26}
      \field{issn}{0098-3500}
      \field{journaltitle}{ACM Transactions on Mathematical Software}
      \field{month}{5}
      \field{number}{2}
      \field{shortjournal}{ACM Trans. Math. Softw.}
      \field{title}{Construction of {{Arbitrary Order Finite Element Degree-of-Freedom Maps}} on {{Polygonal}} and {{Polyhedral Cell Meshes}}}
      \field{urlday}{24}
      \field{urlmonth}{10}
      \field{urlyear}{2023}
      \field{volume}{48}
      \field{year}{2022}
      \field{dateera}{ce}
      \field{urldateera}{ce}
      \field{pages}{18:1\bibrangedash 18:23}
      \range{pages}{-1}
      \verb{doi}
      \verb 10.1145/3524456
      \endverb
      \verb{file}
      \verb /home/rrobin/Zotero/storage/N7WT8XRK/Scroggs et al. - 2022 - Construction of Arbitrary Order Finite Element Deg.pdf
      \endverb
      \verb{urlraw}
      \verb https://dl.acm.org/doi/10.1145/3524456
      \endverb
      \verb{url}
      \verb https://dl.acm.org/doi/10.1145/3524456
      \endverb
    \endentry
    \entry{valliVariationalInterpretationBiot2019}{article}{}
      \name{author}{1}{}{%
        {{hash=d8936df78ee25f351f67bb5b57854a86}{%
           family={Valli},
           familyi={V\bibinitperiod},
           given={Alberto},
           giveni={A\bibinitperiod}}}%
      }
      \list{publisher}{1}{%
        {American Institute of Physics}%
      }
      \strng{namehash}{d8936df78ee25f351f67bb5b57854a86}
      \strng{fullhash}{d8936df78ee25f351f67bb5b57854a86}
      \strng{bibnamehash}{d8936df78ee25f351f67bb5b57854a86}
      \strng{authorbibnamehash}{d8936df78ee25f351f67bb5b57854a86}
      \strng{authornamehash}{d8936df78ee25f351f67bb5b57854a86}
      \strng{authorfullhash}{d8936df78ee25f351f67bb5b57854a86}
      \field{labelalpha}{Val19}
      \field{sortinit}{V}
      \field{sortinithash}{02432525618c08e2b03cac47c19764af}
      \field{labelnamesource}{author}
      \field{labeltitlesource}{title}
      \field{abstract}{In this note, we show that the projection of the Biot–Savart operator over the space of divergence-free vector fields that are tangential to the boundary is the solution of a suitable saddle-point variational problem. Since this projected Biot–Savart operator is shown to be compact, its spectrum can be completely characterized. In particular, through a suitable finite element discretization, it becomes possible to compute the helicity of a bounded domain of a general topological shape, via the determination of the eigenvalue of the projected Biot–Savart operator that has a maximum absolute value.}
      \field{issn}{0022-2488}
      \field{journaltitle}{Journal of Mathematical Physics}
      \field{month}{2}
      \field{number}{2}
      \field{shortjournal}{J. Math. Phys.}
      \field{title}{A variational interpretation of the {{Biot}}–{{Savart}} operator and the helicity of a bounded domain}
      \field{urlday}{1}
      \field{urlmonth}{7}
      \field{urlyear}{2022}
      \field{volume}{60}
      \field{year}{2019}
      \field{dateera}{ce}
      \field{urldateera}{ce}
      \field{pages}{021503}
      \range{pages}{1}
      \verb{doi}
      \verb 10.1063/1.5024197
      \endverb
      \verb{file}
      \verb /home/rrobin/Zotero/storage/YEMENYWQ/Valli - 2019 - A variational interpretation of the Biot–Savart op.pdf
      \endverb
      \verb{urlraw}
      \verb https://aip.scitation.org/doi/10.1063/1.5024197
      \endverb
      \verb{url}
      \verb https://aip.scitation.org/doi/10.1063/1.5024197
      \endverb
    \endentry
    \entry{virtanenSciPyFundamentalAlgorithms2020}{article}{}
      \name{author}{35}{}{%
        {{hash=18703a2bb6a62484483c193a212da2f8}{%
           family={Virtanen},
           familyi={V\bibinitperiod},
           given={Pauli},
           giveni={P\bibinitperiod}}}%
        {{hash=646fbfe08374cc41c2f9bd971d8c4725}{%
           family={Gommers},
           familyi={G\bibinitperiod},
           given={Ralf},
           giveni={R\bibinitperiod}}}%
        {{hash=d500f4849030f34359cdb3e1513acf83}{%
           family={Oliphant},
           familyi={O\bibinitperiod},
           given={Travis\bibnamedelima E.},
           giveni={T\bibinitperiod\bibinitdelim E\bibinitperiod}}}%
        {{hash=35bb9c71f55048509a3e9018c349ed73}{%
           family={Haberland},
           familyi={H\bibinitperiod},
           given={Matt},
           giveni={M\bibinitperiod}}}%
        {{hash=fbb0c40f5d70be8ce47ce9daafdf5749}{%
           family={Reddy},
           familyi={R\bibinitperiod},
           given={Tyler},
           giveni={T\bibinitperiod}}}%
        {{hash=9fd9ed8466bbb96364ae008f2a665e6e}{%
           family={Cournapeau},
           familyi={C\bibinitperiod},
           given={David},
           giveni={D\bibinitperiod}}}%
        {{hash=09a667aa6a26526bfcccb2676a494e55}{%
           family={Burovski},
           familyi={B\bibinitperiod},
           given={Evgeni},
           giveni={E\bibinitperiod}}}%
        {{hash=3d6efaaa3d9682e20787eb06ff70a3d7}{%
           family={Peterson},
           familyi={P\bibinitperiod},
           given={Pearu},
           giveni={P\bibinitperiod}}}%
        {{hash=4c7e4c94b846fa41e2fc0a88e0dc656d}{%
           family={Weckesser},
           familyi={W\bibinitperiod},
           given={Warren},
           giveni={W\bibinitperiod}}}%
        {{hash=7447cb057596bc2645d3980bb04f5c78}{%
           family={Bright},
           familyi={B\bibinitperiod},
           given={Jonathan},
           giveni={J\bibinitperiod}}}%
        {{useprefix=true,hash=9a5a65a789013a8d1e8035ec28df9b6e}{%
           family={Walt},
           familyi={W\bibinitperiod},
           given={Stéfan\bibnamedelima J.},
           giveni={S\bibinitperiod\bibinitdelim J\bibinitperiod},
           prefix={van\bibnamedelima der},
           prefixi={v\bibinitperiod\bibinitdelim d\bibinitperiod}}}%
        {{hash=626cc151613864abeb653c0d8172d98c}{%
           family={Brett},
           familyi={B\bibinitperiod},
           given={Matthew},
           giveni={M\bibinitperiod}}}%
        {{hash=57849e8550281b202bd611bf6f11e14b}{%
           family={Wilson},
           familyi={W\bibinitperiod},
           given={Joshua},
           giveni={J\bibinitperiod}}}%
        {{hash=b053969d2c6a9ec8689980fb6463cd56}{%
           family={Millman},
           familyi={M\bibinitperiod},
           given={K.\bibnamedelimi Jarrod},
           giveni={K\bibinitperiod\bibinitdelim J\bibinitperiod}}}%
        {{hash=fbaf80580622bd40577f4a6d38021c0a}{%
           family={Mayorov},
           familyi={M\bibinitperiod},
           given={Nikolay},
           giveni={N\bibinitperiod}}}%
        {{hash=7bcf847eaccba039f7a4523540673aea}{%
           family={Nelson},
           familyi={N\bibinitperiod},
           given={Andrew\bibnamedelimb R.\bibnamedelimi J.},
           giveni={A\bibinitperiod\bibinitdelim R\bibinitperiod\bibinitdelim J\bibinitperiod}}}%
        {{hash=4b3d26f886661aa723985bcfd835ba18}{%
           family={Jones},
           familyi={J\bibinitperiod},
           given={Eric},
           giveni={E\bibinitperiod}}}%
        {{hash=9ad1d38817acd2f00cb7f324ec7d37ea}{%
           family={Kern},
           familyi={K\bibinitperiod},
           given={Robert},
           giveni={R\bibinitperiod}}}%
        {{hash=8d336f110675c46226ece1db501ce712}{%
           family={Larson},
           familyi={L\bibinitperiod},
           given={Eric},
           giveni={E\bibinitperiod}}}%
        {{hash=65b1934a87acb0abe09c469aaf11c326}{%
           family={Carey},
           familyi={C\bibinitperiod},
           given={C\bibnamedelima J},
           giveni={C\bibinitperiod\bibinitdelim J\bibinitperiod}}}%
        {{hash=9989e8a18827e34b15112d671f52bd35}{%
           family={Polat},
           familyi={P\bibinitperiod},
           given={İlhan},
           giveni={İ\bibinitperiod}}}%
        {{hash=b8b88d61c79de60e6e1b5d44e03f5dec}{%
           family={Feng},
           familyi={F\bibinitperiod},
           given={Yu},
           giveni={Y\bibinitperiod}}}%
        {{hash=bf4be16325cb4f641345ca394443fd18}{%
           family={Moore},
           familyi={M\bibinitperiod},
           given={Eric\bibnamedelima W.},
           giveni={E\bibinitperiod\bibinitdelim W\bibinitperiod}}}%
        {{hash=0fd9a0e34f1b2adda41357c948d14986}{%
           family={VanderPlas},
           familyi={V\bibinitperiod},
           given={Jake},
           giveni={J\bibinitperiod}}}%
        {{hash=c6a95a8ced3b86b4e7e60a74bc6ebf5a}{%
           family={Laxalde},
           familyi={L\bibinitperiod},
           given={Denis},
           giveni={D\bibinitperiod}}}%
        {{hash=85242652d69220e83cf71ceb8d90a8cb}{%
           family={Perktold},
           familyi={P\bibinitperiod},
           given={Josef},
           giveni={J\bibinitperiod}}}%
        {{hash=5bca159e697db439e23b8947dfa4b614}{%
           family={Cimrman},
           familyi={C\bibinitperiod},
           given={Robert},
           giveni={R\bibinitperiod}}}%
        {{hash=70b659f5067a8a2efbee66f770681598}{%
           family={Henriksen},
           familyi={H\bibinitperiod},
           given={Ian},
           giveni={I\bibinitperiod}}}%
        {{hash=fa5163c76600eb11a4d07a28f0701cb0}{%
           family={Quintero},
           familyi={Q\bibinitperiod},
           given={E.\bibnamedelimi A.},
           giveni={E\bibinitperiod\bibinitdelim A\bibinitperiod}}}%
        {{hash=db2b4761cc46be347b418e68660c9554}{%
           family={Harris},
           familyi={H\bibinitperiod},
           given={Charles\bibnamedelima R.},
           giveni={C\bibinitperiod\bibinitdelim R\bibinitperiod}}}%
        {{hash=7d86aea5ad1f2b4e27f2f014c71712c2}{%
           family={Archibald},
           familyi={A\bibinitperiod},
           given={Anne\bibnamedelima M.},
           giveni={A\bibinitperiod\bibinitdelim M\bibinitperiod}}}%
        {{hash=3876f7c3dbbb1a17823dcd135d07cfc6}{%
           family={Ribeiro},
           familyi={R\bibinitperiod},
           given={Antônio\bibnamedelima H.},
           giveni={A\bibinitperiod\bibinitdelim H\bibinitperiod}}}%
        {{hash=bab4e5caee2d67831e464ce575022b37}{%
           family={Pedregosa},
           familyi={P\bibinitperiod},
           given={Fabian},
           giveni={F\bibinitperiod}}}%
        {{useprefix=true,hash=1af2e1049c0f42049401999babb9f7b2}{%
           family={Mulbregt},
           familyi={M\bibinitperiod},
           given={Paul},
           giveni={P\bibinitperiod},
           prefix={van},
           prefixi={v\bibinitperiod}}}%
        {{hash=aa8bf7a30651c7bc3d20ff02fc843dd9}{%
           family={{SciPy 1.0 Contributors}},
           familyi={S\bibinitperiod}}}%
      }
      \strng{namehash}{f252038edd7f112d75da3bf0c1edecbc}
      \strng{fullhash}{60816f7538af4874a8d3aab64c605ef0}
      \strng{bibnamehash}{f252038edd7f112d75da3bf0c1edecbc}
      \strng{authorbibnamehash}{f252038edd7f112d75da3bf0c1edecbc}
      \strng{authornamehash}{f252038edd7f112d75da3bf0c1edecbc}
      \strng{authorfullhash}{60816f7538af4874a8d3aab64c605ef0}
      \field{labelalpha}{Vir+20}
      \field{sortinit}{V}
      \field{sortinithash}{02432525618c08e2b03cac47c19764af}
      \field{labelnamesource}{author}
      \field{labeltitlesource}{title}
      \field{journaltitle}{Nature Methods}
      \field{title}{{{SciPy}} 1.0: {{Fundamental}} algorithms for scientific computing in python}
      \field{volume}{17}
      \field{year}{2020}
      \field{dateera}{ce}
      \field{pages}{261\bibrangedash 272}
      \range{pages}{12}
      \verb{doi}
      \verb 10.1038/s41592-019-0686-2
      \endverb
    \endentry
    \entry{warmerW7XHELIASFusion2017}{article}{}
      \name{author}{8}{}{%
        {{hash=88abbdb3449cb10854897d2d522940c3}{%
           family={Warmer},
           familyi={W\bibinitperiod},
           given={F.},
           giveni={F\bibinitperiod}}}%
        {{hash=f8fabb85617d7d4898a242ba324b0719}{%
           family={Bykov},
           familyi={B\bibinitperiod},
           given={V.},
           giveni={V\bibinitperiod}}}%
        {{hash=416ca898870adfc6319dba80e38e051c}{%
           family={Drevlak},
           familyi={D\bibinitperiod},
           given={M.},
           giveni={M\bibinitperiod}}}%
        {{hash=eaacd80e918b1e0383435b117ab78afe}{%
           family={Häußler},
           familyi={H\bibinitperiod},
           given={A.},
           giveni={A\bibinitperiod}}}%
        {{hash=feac46feceb5b9f910b0f51fb4d96a90}{%
           family={Fischer},
           familyi={F\bibinitperiod},
           given={U.},
           giveni={U\bibinitperiod}}}%
        {{hash=19d72854ff37e6747bfc201adedbc526}{%
           family={Stange},
           familyi={S\bibinitperiod},
           given={T.},
           giveni={T\bibinitperiod}}}%
        {{hash=c3527144ba8d59b25987e3d50bb90c53}{%
           family={Beidler},
           familyi={B\bibinitperiod},
           given={C.\bibnamedelimi D.},
           giveni={C\bibinitperiod\bibinitdelim D\bibinitperiod}}}%
        {{hash=021e6421201d25be9eee86d28937b13f}{%
           family={Wolf},
           familyi={W\bibinitperiod},
           given={R.\bibnamedelimi C.},
           giveni={R\bibinitperiod\bibinitdelim C\bibinitperiod}}}%
      }
      \strng{namehash}{f6fbba32547696b5063a307869bd1170}
      \strng{fullhash}{2608a3b8f0db6c2ab36e4c9c6ffb2bed}
      \strng{bibnamehash}{f6fbba32547696b5063a307869bd1170}
      \strng{authorbibnamehash}{f6fbba32547696b5063a307869bd1170}
      \strng{authornamehash}{f6fbba32547696b5063a307869bd1170}
      \strng{authorfullhash}{2608a3b8f0db6c2ab36e4c9c6ffb2bed}
      \field{labelalpha}{War+17}
      \field{sortinit}{W}
      \field{sortinithash}{1af34bd8c148ffb32de1494636b49713}
      \field{labelnamesource}{author}
      \field{labeltitlesource}{shorttitle}
      \field{abstract}{With the recent start of operation of Wendelstein 7-X (W7-X), the helical-axis advanced stellarator concept (HELIAS) has gained renewed attention. In particular a discussion has been started about a research strategy leading from W7-X to a commercial HELIAS fusion power plant. In order to bridge the respective gap in physics and technology between those devices, concepts for an intermediate-step burning plasma stellarator are under discussion. Recent studies have focused on physics issues of importance to such devices. Extending this discussion, engineering and technology considerations for next-step HELIAS devices are discussed in this work.}
      \field{day}{1}
      \field{issn}{0920-3796}
      \field{journaltitle}{Fusion Engineering and Design}
      \field{langid}{english}
      \field{month}{11}
      \field{series}{Proceedings of the 29th {{Symposium}} on {{Fusion Technology}} ({{SOFT-29}}) {{Prague}}, {{Czech Republic}}, {{September}} 5-9, 2016}
      \field{shortjournal}{Fusion Engineering and Design}
      \field{shorttitle}{From {{W7-X}} to a {{HELIAS}} fusion power plant}
      \field{title}{From {{W7-X}} to a {{HELIAS}} fusion power plant: {{On}} engineering considerations for next-step stellarator devices}
      \field{urlday}{9}
      \field{urlmonth}{3}
      \field{urlyear}{2020}
      \field{volume}{123}
      \field{year}{2017}
      \field{dateera}{ce}
      \field{urldateera}{ce}
      \field{pages}{47\bibrangedash 53}
      \range{pages}{7}
      \verb{doi}
      \verb 10.1016/j.fusengdes.2017.05.034
      \endverb
      \verb{file}
      \verb /home/rrobin/Zotero/storage/BE45CA5W/Warmer et al. - 2017 - From W7-X to a HELIAS fusion power plant On engin.pdf;/home/rrobin/Zotero/storage/J9BYIIRQ/Warmer et al. - 2017 - From W7-X to a HELIAS fusion power plant On engin.pdf;/home/rrobin/Zotero/storage/S88H7EJ3/Warmer et al. - 2017 - From W7-X to a HELIAS fusion power plant On engin.pdf;/home/rrobin/Zotero/storage/VWT9Y2FM/Warmer et al. - 2017 - From W7-X to a HELIAS fusion power plant On engin.pdf;/home/rrobin/Zotero/storage/FQN3DVIC/S0920379617305781.html;/home/rrobin/Zotero/storage/SDTHVEMB/S0920379617305781.html;/home/rrobin/Zotero/storage/SM3UMP57/S0920379617305781.html;/home/rrobin/Zotero/storage/X6L4V6YR/S0920379617305781.html;/home/rrobin/Zotero/storage/YY4AD5FL/S0920379617305781.html
      \endverb
      \verb{urlraw}
      \verb http://www.sciencedirect.com/science/article/pii/S0920379617305781
      \endverb
      \verb{url}
      \verb http://www.sciencedirect.com/science/article/pii/S0920379617305781
      \endverb
      \keyw{DEMO,Engineering,HELIAS,Intermediate-step stellarator,Technology,W7-X}
    \endentry
    \entry{zarnstorffPhysicsCompactAdvanced2001}{article}{}
      \name{author}{34}{}{%
        {{hash=410025aef8412310f509481130858e2b}{%
           family={Zarnstorff},
           familyi={Z\bibinitperiod},
           given={M.\bibnamedelimi C.},
           giveni={M\bibinitperiod\bibinitdelim C\bibinitperiod}}}%
        {{hash=eb5154d76ea355f8d1e0ae7ce6e7b209}{%
           family={Berry},
           familyi={B\bibinitperiod},
           given={L.\bibnamedelimi A.},
           giveni={L\bibinitperiod\bibinitdelim A\bibinitperiod}}}%
        {{hash=390c37b095037e8875261f10b95b711d}{%
           family={Brooks},
           familyi={B\bibinitperiod},
           given={A.},
           giveni={A\bibinitperiod}}}%
        {{hash=7c2ba6f1596eae22d4675116e4434a5c}{%
           family={Fredrickson},
           familyi={F\bibinitperiod},
           given={E.},
           giveni={E\bibinitperiod}}}%
        {{hash=aede3346dda4d261e18650fcd48f2e90}{%
           family={Fu},
           familyi={F\bibinitperiod},
           given={G.-Y.},
           giveni={G\bibinithyphendelim Y\bibinitperiod}}}%
        {{hash=28d8e3891a8aba921364acccbd9d4f7e}{%
           family={Hirshman},
           familyi={H\bibinitperiod},
           given={S.},
           giveni={S\bibinitperiod}}}%
        {{hash=fc20750649ce459731e392560c73b949}{%
           family={Hudson},
           familyi={H\bibinitperiod},
           given={S.},
           giveni={S\bibinitperiod}}}%
        {{hash=4f0d910aebef31288c885cd55d208f69}{%
           family={Ku},
           familyi={K\bibinitperiod},
           given={L.-P.},
           giveni={L\bibinithyphendelim P\bibinitperiod}}}%
        {{hash=870fb3624a3009d2f0d456bc28ecf1d4}{%
           family={Lazarus},
           familyi={L\bibinitperiod},
           given={E.},
           giveni={E\bibinitperiod}}}%
        {{hash=ff660b3c9f264a1069c45c338c77a1b8}{%
           family={Mikkelsen},
           familyi={M\bibinitperiod},
           given={D.},
           giveni={D\bibinitperiod}}}%
        {{hash=8b1d535a26bb02cff2731cc49209a387}{%
           family={Monticello},
           familyi={M\bibinitperiod},
           given={D.},
           giveni={D\bibinitperiod}}}%
        {{hash=4a7621742bc24639aa505cfdf47a9070}{%
           family={Neilson},
           familyi={N\bibinitperiod},
           given={G.\bibnamedelimi H.},
           giveni={G\bibinitperiod\bibinitdelim H\bibinitperiod}}}%
        {{hash=649fa29460402ef274a46ec2594f344b}{%
           family={Pomphrey},
           familyi={P\bibinitperiod},
           given={N.},
           giveni={N\bibinitperiod}}}%
        {{hash=3c5e1f7db56ae8e354aa42db8be60fbb}{%
           family={Reiman},
           familyi={R\bibinitperiod},
           given={A.},
           giveni={A\bibinitperiod}}}%
        {{hash=72635612c49d8b29e668f4b2ba41d0b7}{%
           family={Spong},
           familyi={S\bibinitperiod},
           given={D.},
           giveni={D\bibinitperiod}}}%
        {{hash=5c19cb3a1fc8caff8c4b5b1340c1f8c3}{%
           family={Strickler},
           familyi={S\bibinitperiod},
           given={D.},
           giveni={D\bibinitperiod}}}%
        {{hash=c9ddb6741c61609edea3fc3ff4899bd5}{%
           family={Boozer},
           familyi={B\bibinitperiod},
           given={A.},
           giveni={A\bibinitperiod}}}%
        {{hash=502aff825efffcf3f61e6da0de924e3d}{%
           family={Cooper},
           familyi={C\bibinitperiod},
           given={W.\bibnamedelimi A.},
           giveni={W\bibinitperiod\bibinitdelim A\bibinitperiod}}}%
        {{hash=9de3026021e4d8d566534b9a31eaf846}{%
           family={Goldston},
           familyi={G\bibinitperiod},
           given={R.},
           giveni={R\bibinitperiod}}}%
        {{hash=66d0ce3d4b67129967765913285da3aa}{%
           family={Hatcher},
           familyi={H\bibinitperiod},
           given={R.},
           giveni={R\bibinitperiod}}}%
        {{hash=f4cc6e7bda56699e700f3055ee0dceef}{%
           family={Isaev},
           familyi={I\bibinitperiod},
           given={M.},
           giveni={M\bibinitperiod}}}%
        {{hash=5aa696b9e2810f0c8198bc38a0e4d282}{%
           family={Kessel},
           familyi={K\bibinitperiod},
           given={C.},
           giveni={C\bibinitperiod}}}%
        {{hash=2a6637281d28108fe05806663f8053f1}{%
           family={Lewandowski},
           familyi={L\bibinitperiod},
           given={J.},
           giveni={J\bibinitperiod}}}%
        {{hash=6386d6dd7bbb2d9cba740c7e7dedc059}{%
           family={Lyon},
           familyi={L\bibinitperiod},
           given={J.\bibnamedelimi F.},
           giveni={J\bibinitperiod\bibinitdelim F\bibinitperiod}}}%
        {{hash=c9c2c43f7d206b4c1da41b727bc15a22}{%
           family={Merkel},
           familyi={M\bibinitperiod},
           given={P.},
           giveni={P\bibinitperiod}}}%
        {{hash=a09aab77f923101fed15a7941bc0dafa}{%
           family={Mynick},
           familyi={M\bibinitperiod},
           given={H.},
           giveni={H\bibinitperiod}}}%
        {{hash=924dacd6dcf72f97f4ef43a8aaac8e60}{%
           family={Nelson},
           familyi={N\bibinitperiod},
           given={B.\bibnamedelimi E.},
           giveni={B\bibinitperiod\bibinitdelim E\bibinitperiod}}}%
        {{hash=6b490063400f1ec760f47c06b8090c17}{%
           family={Nuehrenberg},
           familyi={N\bibinitperiod},
           given={C.},
           giveni={C\bibinitperiod}}}%
        {{hash=6f28c7b6377d5955b5a9bfd44bdec673}{%
           family={Redi},
           familyi={R\bibinitperiod},
           given={M.},
           giveni={M\bibinitperiod}}}%
        {{hash=3adb8eb2a5c47cd0570c48f66e761e14}{%
           family={Reiersen},
           familyi={R\bibinitperiod},
           given={W.},
           giveni={W\bibinitperiod}}}%
        {{hash=666c949efd09fcf954c8e75d2794718f}{%
           family={Rutherford},
           familyi={R\bibinitperiod},
           given={P.},
           giveni={P\bibinitperiod}}}%
        {{hash=22237c308d0c829be2f77a7d4ba8ab1c}{%
           family={Sanchez},
           familyi={S\bibinitperiod},
           given={R.},
           giveni={R\bibinitperiod}}}%
        {{hash=12b671e73ce9fe1c664efeaf7a509ae9}{%
           family={Schmidt},
           familyi={S\bibinitperiod},
           given={J.},
           giveni={J\bibinitperiod}}}%
        {{hash=bc7735e692a37788f4ab6a5c08ff997f}{%
           family={White},
           familyi={W\bibinitperiod},
           given={R.\bibnamedelimi B.},
           giveni={R\bibinitperiod\bibinitdelim B\bibinitperiod}}}%
      }
      \list{publisher}{1}{%
        {IOP Publishing}%
      }
      \strng{namehash}{52cc253d8d6f2b456085500df3443ea2}
      \strng{fullhash}{24cdd857e26f59d4a51d72ac17b29f16}
      \strng{bibnamehash}{52cc253d8d6f2b456085500df3443ea2}
      \strng{authorbibnamehash}{52cc253d8d6f2b456085500df3443ea2}
      \strng{authornamehash}{52cc253d8d6f2b456085500df3443ea2}
      \strng{authorfullhash}{24cdd857e26f59d4a51d72ac17b29f16}
      \field{labelalpha}{Zar+01}
      \field{sortinit}{Z}
      \field{sortinithash}{8f7b480688e809b50b6f6577b16f3db5}
      \field{labelnamesource}{author}
      \field{labeltitlesource}{title}
      \field{issue}{12A}
      \field{journaltitle}{Plasma Physics and Controlled Fusion}
      \field{month}{11}
      \field{shortjournal}{Plasma Phys. Control. Fusion}
      \field{title}{Physics of the compact advanced stellarator {{NCSX}}}
      \field{volume}{43}
      \field{year}{2001}
      \field{dateera}{ce}
      \field{pages}{A237\bibrangedash A249}
      \range{pages}{-1}
      \verb{doi}
      \verb 10.1088/0741-3335/43/12a/318
      \endverb
      \verb{urlraw}
      \verb https://doi.org/10.1088/0741-3335/43/12a/318
      \endverb
      \verb{url}
      \verb https://doi.org/10.1088/0741-3335/43/12a/318
      \endverb
    \endentry
  \enddatalist
\endrefsection

  \blx@bblend
  \endgroup
  \csnumgdef{blx@labelnumber@\the\c@refsection}{0}}
\makeatother

\title{Shape optimization of harmonic helicity in toroidal domains}

\author[1]{Rémi Robin \thanks{remi.robin@inria.fr}}
\author[2]{Robin Roussel \thanks{robin.roussel@sorbonne-universite.fr}}

\affil[1]{Laboratoire de Physique de l’\'Ecole Normale Supérieure, Mines Paris, Inria, CNRS, ENS-PSL, Sorbonne Université, PSL Research University, Paris, France}
\affil[2]{Laboratoire Jacques-Louis Lions, Sorbonne Université, Inria, Paris}

\begin{document}
\maketitle
\begin{abstract}
    In this paper, we introduce a new shape functional defined for toroidal domains that we call harmonic helicity, and study its shape optimization. Given a toroidal domain, we consider its associated harmonic field. The latter is the magnetic field obtained uniquely up to normalization when imposing zero normal trace and zero electrical current inside the domain. We then study the helicity of this field, which is a quantity of interest in magneto-hydrodynamics corresponding to the $L^2$ product of the field with its image by the Biot--Savart operator. To do so, we begin by discussing the appropriate functional framework and an equivalent PDE characterization. We then focus on shape optimization, and we identify the shape gradient of the harmonic helicity. Finally, we study and implement an efficient numerical scheme to compute harmonic helicity and its shape gradient using finite elements exterior calculus.
\end{abstract}
\tableofcontents
\section{Introduction}


For a given vector field $F$ on a three-dimensional domain $\Omega$, we define its helicity (also known as Biot--Savart helicity) by the formula
\begin{align}
    \label{eq:first_eq_helicity}
    \H(F)=\frac{1}{4 \pi}\int_{\Omega \times \Omega} F(y)\cdot \left(F(x) \times \frac{y-x}{|y-x|^3}\right) dx dy.
\end{align}
This quantity plays an important role in plasma physics, fluid dynamics and magnetohydrodynamics (see e.g. \cite{arnoldTopologyThreedimensionalSteady1966,arnoldTopologicalMethodsHydrodynamics2021}). In the context of electromagnetism, helicity of a magnetic field (called magnetic helicity) can be seen as a scalar quantifying the linkage and twist of the magnetic field \cite{moffattDegreeKnottednessTangled1969,Arnold2014}. \Cref{eq:first_eq_helicity} can be interpreted as a volumic version of the writhe of a curve.

When working with magnetic fields, that is divergence free vector fields, tangent to the boundary of $\Omega$, a natural connection with vector potentials appears. First note that if one introduces the Biot--Savart operator of $F$
\begin{equation}
    \BS(F)(y)=\frac{1}{4\pi}\int_{\Omega} \frac{F(x)\times (y-x)}{|y-x|^3} dx,
\end{equation}
one obtains that the helicity of $F$ is the $L^2$ inner product of $F$ with $\BS(F)$. For simply connected domains, the latter remains true if we replace $\BS(F)$ by any vector potential of $F$ (we recall that $\curl \BS (F)=F$). For three-dimensional domains that are not simply connected, the connection between helicity and vector potentials is slightly more involved and have been established by Bevir and Gray in \cite{bevirRelaxationFluxConsumption1980} for toroidal domains\footnote{\label{ft:toroidal_domain} That is bounded open subsets of $\R^3$ homeomorphic to a full torus.} and \cite{mactaggartMagneticHelicityMultiply2019} for more general ones. We recall the definition and main properties of the magnetic helicity in \cref{subsec:VP_BG}.

A classical mathematical problem related to helicity studied in e.g. \cite{cantarellaIsoperimetricProblemsHelicity2000,cantarella_influence_1999,valliVariationalInterpretationBiot2019,montielIsoperimetricProblemCurl2023} is the maximization of the helicity on $H_0(\div^0, \Omega)$\footnote{we refer to \cref{subsec-functional} for the definition of this Hilbert space.} with fixed $L^2$ norm. The critical vector field of this optimization problem are in fact eigenfields of the curl operator which can be seen thanks to a modified Biot--Savart operator denoted here $\BS'$ \cite{cantarellaIsoperimetricProblemsHelicity2000}.

Elements of shape optimization were also discussed in \cite{cantarellaIsoperimetricProblemsHelicity2000} and \cite{cantarella_influence_1999} to characterize the domain with the highest eigenvalue of $\BS'$ for a given volume. Some properties of the maximizing fields on such a domain were given, but it is also unclear whether the optimal shape exists, as some computations suggest that it would have to be a singular sphere, with North and South Pole collapsed to a single point. Recent works on the existence of an optimal shape and its characterization can be found in \cite{encisoNonexistenceAxisymmetricOptimal2020,gernerIsoperimetricProblemFirst2023,gernerExistenceOptimalDomains2023}.

In this paper, we are interested in a slightly different problem. Given a toroidal domain $\Omega$, we consider the set of harmonic fields, that is the set of vector fields that are divergence free, curl free and tangent to the boundary. By classical results of Hodge theory, this set is a one-dimensional vector space. We are then interested in the helicity of a normalized harmonic field of $\Omega$, where the normalization is related to total flux of currents through the central hole of the torus. Thus, for any regular enough toroidal domain, we define a scalar quantity that we call the \textbf{harmonic helicity} of the domain $\Omega$.

Designing a numerical scheme to compute this shape functional is not obvious. A close problem is the spectral approximation of the curl operator in multiply connected domains; this has been tackled in \cite{laraSpectralApproximationCurl2015,alonso-rodriguezFiniteElementApproximation2018} using finite elements methods. For efficiency considerations, it is important to avoid the computation of the double integral in \cref{eq:first_eq_helicity} and use another vector potential than the Biot--Savart. Using classical results on vector potentials characterizations \cite{amroucheVectorPotentialThreedimensional1998} and tools from finite elements exterior calculus \cite{arnoldFiniteElementExterior2010}, we provide efficient numerical approximation schemes and implementation for the harmonic helicity.

Physical motivations for considering harmonic helicity arise from the design of stellarators, advanced nuclear fusion devices that rely on the confinement of intensely hot plasma through a sophisticated magnetic field. A significant challenge arises due to the inherent impossibility of creating a non-zero magnetic field with constant magnitude on an axisymmetric toroidal domain. Additionally, variations in the magnetic field amplitude, typically inversely proportional to the major radius, result in a vertical drift, which can be mitigated through the implementation of a twisted magnetic field. Optimization of the shape of the coils is a very active field \cite{paulAdjointMethodGradientbased2018,privatOptimalShapeStellarators2022}. A measure of the twisted nature of the magnetic field is expressed by the magnetic helicity. We refer to \cite{imbert-gerardIntroductionSymmetriesStellarators2019} for a very nice introduction to the topic.

In contrast to Tokamaks, which are axisymmetric devices inducing a current inside the plasma to generate a twisting magnetic field, stellarators aim for stability without requiring a current within the plasma. Consequently, the magnetic field employed to stabilize the plasma in a stellarator can be reasonably approximated as a harmonic field within the domain representing the plasma. Hence, we believe that the optimization of the shape of the plasma to increase the harmonic helicity  could give rise to interesting new forms of plasma.

This paper is organized as follows:
\begin{itemize}
    \item In \cref{sec:state_of_the_art}, we recall classical notions used throughout the paper. We begin by the definitions of some functional spaces. Then we properly define the harmonic fields, give two equivalent constructions and define their circulations. Next, we give two PDE formulations to characterize vector potentials. Finally, we make the connection between vector potentials and the helicity through the Bevir--Gray formula.
    \item \Cref{sec:HH_shape_gradient} contains the main contribution of the paper. We begin by introducing the precise definition of the harmonic helicity of a toroidal domain. The rest of the section is dedicated to the computation of the shape derivative of the harmonic helicity. This is done using shape differentiation of several PDE problems and a subtle use of Piola transforms. To the best of the authors' knowledge, the employed methodology is original and may hold applicability in addressing new shape differentiation problems involving other Hilbert complexes.
    \item In \cref{sec:approx_fecc}, we recall the framework of finite element exterior calculus. We then use classical results on approximations of Hodge Laplacian problems. The adaptations of these tools to our problem is not straightforward and provides a method to compute both the helicity and the shape gradient.
    \item In \cref{sec:num}, we provide numerical results of the proposed numerical methods on specific shapes motivated by the study of stellarator plasmas. Then we present two numerical experiments to improve the harmonic helicity of a standard plasma shape.
    \item In \cref{sec:appendix_hodge}, we recall classical results of Hodge theory used throughout the paper. In particular the Hodge decomposition and its connection with the De Rham cohomology.
\end{itemize}

\section{Prerequisites}
\label{sec:state_of_the_art}
In this section, $\Omega$ denotes a Lipschitz toroidal domain of $\R^3$, that is $\partial \Omega$ is locally the graph of a Lipschitz function and $\bar{\Omega}$ is homeomorphic to $D^2\times S^1$ with $D^2$ the closed unit disk and $S^1$ the unit circle. 
%
%
%
%

\subsection{Functional spaces}
\label{subsec-functional}
We recall the definitions of the following classical functional spaces:
\begin{gather}\label{eq:functional_spaces}
    \begin{split}
        &H(\curl, \Omega) = \left\{V \in L^2(\Omega)^3 \mid \curl V\in L^2(\Omega)^3\right\},
        \\
        &H(\div, \Omega) = \left\{V \in L^2(\Omega)^3 \mid \div V\in L^2(\Omega)\right\}.
    \end{split}
\end{gather}
On $H(\curl, \Omega)$ and $H(\div, \Omega)$, the tangential and normal traces $V\times n : \partial \Omega \rightarrow \R^3$ and $V \cdot n: \partial \Omega \rightarrow \R$ are defined respectively by
\begin{equation}\label{eq:IPP_curl}
    \int_{\partial \Omega} (V \times n) \cdot \varphi = { \int_{\Omega} V \cdot \curl \varphi - \int_{\Omega} \curl V \cdot \varphi},
\end{equation}
for every $\varphi$ in $H^1(\Omega)^3$, and 
\begin{equation}\label{eq:IPP_div}
    \int_{\partial \Omega} (V \cdot n) \varphi = { \int_{\Omega} V \cdot \nabla \varphi + \int_{\Omega} \div V \varphi},
\end{equation}
for every $\varphi$ in $H^1(\Omega)$. Since the traces of $H^1(\Omega)^3$ and $H^1(\Omega)$ are $H^{1/2}(\partial \Omega)^3$ and $H^{1/2}(\partial \Omega)$ respectively, $V\times n$ can be defined in $H^{-1/2}(\partial \Omega)^3$\footnote{{In fact, as was shown in \cite{buffaTracesCurlLipschitz2002}, the space of tangential traces of $H(\curl, \Omega)$ is $H^{-1/2}(\div_\Gamma, \partial \Omega)$, so that $H^{-1/2}(\partial \Omega)^3$ is a proper subspace.}}, and $V \cdot n$ in $H^{-1/2}(\partial \Omega)$. Then, we can define 
\begin{gather}\label{eq:functional_spaces_traceless}
    \begin{split}
        &H_0(\curl, \Omega) = \left\{V\in L^2(\Omega)^3 \mid \curl V\in L^2(\Omega)^3, V \times n = 0\right\},
        \\
        &H_0(\div, \Omega) = \left\{V \in L^2(\Omega)^3 \mid \div V\in L^2(\Omega), V \cdot n = 0 \right\}.
    \end{split}
\end{gather}
We also denote by $L^2_0(\Omega)$ the set of functions in $L^2(\Omega)$ which have zero average in $\Omega$. Introducing the following spaces where the differential operator vanishes will also prove to be useful
\begin{gather}\label{eq:functional_spaces^0}
    \begin{split}
        &H\left(\curl^0, \Omega\right)=\left\{V \in H(\curl, \Omega) \mid \curl V = 0\right\},\\
        &H\left(\div^0, \Omega\right)=\left\{V \in H(\div, \Omega) \mid \div V = 0\right\},\\
        &H_0\left(\curl^0, \Omega\right)=H\left(\curl^0, \Omega\right)\cap H_0(\curl, \Omega),\\
        &H_0\left(\div^0, \Omega\right)=H\left(\div^0, \Omega\right)\cap H_0(\div, \Omega).
    \end{split}
\end{gather}

\subsection{Harmonic fields}\label{subsec:harmonic_fields}
The set of harmonic fields in $\Omega$ is defined in the following way
\begin{align}\label{eq:harmonic_fields}
    \mathcal{K}(\Omega)&=\left\{V \in L^2(\Omega)^3\mid \div V =0, \, \curl V=0, V\cdot n =0 \right\} \\
    \notag
    & =H\left(\curl^0, \Omega\right) \cap H_0\left(\div^0, \Omega\right),
\end{align}
Through the identifications with differential forms given in \cref{sec:appendix_hodge}, it is classical that this space is isomorphic to the first De Rham cohomology space of $\Omega$, with each harmonic field giving a natural representent of the corresponding cohomology class. Another characterization based on Hodge decomposition is the following
\begin{align*}
    \mathcal{K}(\Omega)
    &=H\left(\curl^0, \Omega\right) \cap \nabla H^1(\Omega)^\perp.
\end{align*}
We refer to \cref{pr:hodge_decomposition} of \cref{sec:appendix_hodge}.

\begin{figure}
    \center
 
\tikzset{
pattern size/.store in=\mcSize, 
pattern size = 5pt,
pattern thickness/.store in=\mcThickness, 
pattern thickness = 0.3pt,
pattern radius/.store in=\mcRadius, 
pattern radius = 1pt}
\makeatletter
\pgfutil@ifundefined{pgf@pattern@name@_bj65zr824}{
\pgfdeclarepatternformonly[\mcThickness,\mcSize]{_bj65zr824}
{\pgfqpoint{0pt}{0pt}}
{\pgfpoint{\mcSize+\mcThickness}{\mcSize+\mcThickness}}
{\pgfpoint{\mcSize}{\mcSize}}
{
\pgfsetcolor{\tikz@pattern@color}
\pgfsetlinewidth{\mcThickness}
\pgfpathmoveto{\pgfqpoint{0pt}{0pt}}
\pgfpathlineto{\pgfpoint{\mcSize+\mcThickness}{\mcSize+\mcThickness}}
\pgfusepath{stroke}
}}
\makeatother
\tikzset{every picture/.style={line width=0.75pt}} 

\begin{tikzpicture}[x=0.75pt,y=0.75pt,yscale=-1,xscale=1]

\draw  [line width=1.5]  (139.47,53.9) .. controls (180.77,21.5) and (287.86,6) .. (387,17.67) .. controls (486.13,29.33) and (572.13,31.33) .. (575.47,121.33) .. controls (578.8,211.33) and (519.18,261.38) .. (393.47,263.33) .. controls (267.76,265.29) and (191.19,260.85) .. (138.8,222.67) .. controls (86.41,184.48) and (98.16,86.29) .. (139.47,53.9) -- cycle ;
\draw [line width=1.5]    (167.07,123.32) .. controls (190.85,154.79) and (435.6,175.6) .. (495.51,127.95) ;
\draw [line width=1.5]    (188.75,135.75) .. controls (210.51,117.01) and (251.98,108.43) .. (296.84,106.46) .. controls (376.59,102.96) and (467.03,120.36) .. (476.31,138.73) ;
\draw  [color={rgb, 255:red, 74; green, 144; blue, 226 }  ,draw opacity=1 ][fill={rgb, 255:red, 0; green, 0; blue, 0 }  ,fill opacity=0.22 ] (396.87,67.5) .. controls (396.75,42.31) and (404.34,21.89) .. (413.84,21.9) .. controls (423.34,21.91) and (431.14,42.33) .. (431.27,67.52) .. controls (431.39,92.71) and (423.8,113.13) .. (414.3,113.12) .. controls (404.8,113.12) and (397,92.69) .. (396.87,67.5) -- cycle ;
\draw [color={rgb, 255:red, 74; green, 144; blue, 226 }  ,draw opacity=1 ]   (424.37,57.84) -- (431.22,63.89) ;
\draw [color={rgb, 255:red, 74; green, 144; blue, 226 }  ,draw opacity=1 ]   (437.39,57.14) -- (431.22,63.89) ;
\draw  [color={rgb, 255:red, 208; green, 2; blue, 27 }  ,draw opacity=1 ][pattern=_bj65zr824,pattern size=6pt,pattern thickness=0.75pt,pattern radius=0pt, pattern color={rgb, 255:red, 155; green, 155; blue, 155}] (153.6,130.6) .. controls (153.6,98.57) and (235.98,72.6) .. (337.6,72.6) .. controls (439.22,72.6) and (521.6,98.57) .. (521.6,130.6) .. controls (521.6,162.63) and (439.22,188.6) .. (337.6,188.6) .. controls (235.98,188.6) and (153.6,162.63) .. (153.6,130.6) -- cycle ;
\draw [color={rgb, 255:red, 208; green, 2; blue, 27 }  ,draw opacity=1 ]   (193.07,86.84) -- (186.06,98.02) ;
\draw [color={rgb, 255:red, 208; green, 2; blue, 27 }  ,draw opacity=1 ]   (196.21,103.01) -- (186.06,98.02) ;

\draw (436.84,31.7) node [anchor=north west][inner sep=0.75pt]  [color={rgb, 255:red, 74; green, 144; blue, 226 }  ,opacity=1 ]  {$\gamma $};
\draw (406.96,51.44) node [anchor=north west][inner sep=0.75pt]    {$\Sigma $};
\draw (149,87.37) node [anchor=north west][inner sep=0.75pt]    {$\textcolor[rgb]{0.82,0.01,0.11}{\gamma ' }$};
\draw (279,120.4) node [anchor=north west][inner sep=0.75pt]    {$\Sigma '$};
\draw (174,193) node [anchor=north west][inner sep=0.75pt]   [align=left] {$\displaystyle \Omega $};

\end{tikzpicture}
    \caption{Illustration of the curves $\gamma$ and $\gamma'$, and the surfaces $\Sigma$ and $\Sigma'$}
    \label{fig:geometry}
\end{figure}
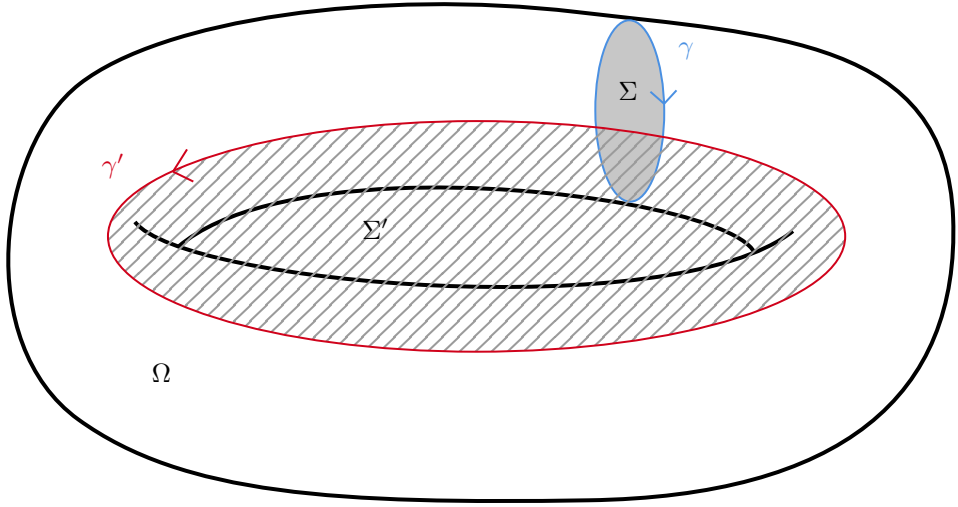

We introduce a poloidal cut $\Sigma$ of $\Omega$, that is, a Lipschitz surface included in $\Omega$, with boundary $\gamma$ contained in $\partial \Omega$ which generates the first homology group of $\Omega^c=\R^3\setminus \Omega$. Similarly, we introduce $\Sigma'$ a Lipschitz surface in $\Omega^c$ with boundary $\gamma'$ that generates the first homology group of $\bar \Omega$. \cref{fig:geometry} illustrates these objects. 
\\
We also define $t$ (resp. $t'$) as unit tangent vector fields on $\gamma$ (resp. $\gamma'$), and $n_\Sigma$ (resp. $n_{\Sigma'}$) as a unit normal vector field on $\Sigma$ (resp. $\Sigma'$). These vector fields define orientations on $\Sigma$, $\Sigma'$, $\gamma$ and $\gamma'$ which are compatible with each other. We refer to \cite[Section 1]{alonso-rodriguezFiniteElementApproximation2018} and \cite[Section 3.a]{amroucheVectorPotentialThreedimensional1998} for the constructions and more precise definitions of these objects.

Formally, we normalize the harmonic field $B$ by fixing its circulation along the toroidal loop $\gamma'$ to be equal to $2\pi$.
%
Hence, as shown in \cite[Proposition 3.14]{amroucheVectorPotentialThreedimensional1998}, we can define the normalized harmonic field $B$ on $\Omega$ by $B = \nabla u$ with $u \in H^1(\Omega \backslash \Sigma)$ defined as follows. $u$ is the unique solution of the variational problem
\begin{equation}\label{eq:B_curl}
    \begin{cases}
        { \int_{\Omega \backslash \Sigma} \nabla u \cdot \nabla v} = 0,
        \\
        [\![u]\!]_\Sigma=2\pi,
    \end{cases}
\end{equation}
for all $v \in H^1(\Omega)$, where $[\![u]\!]_\Sigma$ denote the jump of $u$ across $\Sigma$, that is the difference of traces of $u$ in the direction of $\gamma'$ \cite[Notation 3.9.i]{amroucheVectorPotentialThreedimensional1998}. Note that $\Omega \backslash \Sigma$ is a pseudo-Lipschitz domain, we refer to \cite[Definition 3.1]{amroucheVectorPotentialThreedimensional1998} for further details.

Once we have the normalized harmonic field of $\Omega$, we can define rigorously the circulation of vector fields in $H(\curl, \Omega)$ with curl tangent to the boundary. To do this, we define similarly the harmonic field of $\Omega^{\text{ext}}= \mathcal{B} \cap \bar{\Omega^c}$, denoted $B^{\text{ext}}$, where $\mathcal{B}$ is an open ball containing $\bar{\Omega}$. The normalized harmonic field of $\Omega^{\text{ext}}$ is the gradient of $u^{\text{ext}} \in H^1(\Omega^{\text{ext}} \backslash \Sigma')$, which verifies
\begin{equation}\label{eq:B^ext}
    \begin{cases}
        {\int_{\Omega \backslash \Sigma} \nabla u^{\text{ext}} \cdot \nabla v} = 0,
        \\
        [\![u]\!]_{\Sigma'}=2\pi,
    \end{cases}
\end{equation}
for all $v$ in $H^1\left(\Omega^{\text{ext}}\right)$. Then, as in \cite[Section 2]{alonso-rodriguezFiniteElementApproximation2018}, we are able to define the circulations of a vector field V in $H(\curl, \Omega)$ with curl tangent to the boundary as 
\begin{align}
    \int_{\gamma'} V \cdot t &= \frac{1}{2\pi}\int_{\Omega^{\text{ext}}} \curl \tilde{V} \cdot B^{\text{ext}}, \label{eq:circ_gamma}\\
    \int_{\gamma} V \cdot t' &= \frac{1}{2\pi}\int_{\Omega} \curl V \cdot B, \label{eq:circ_gamma_p}
\end{align}
where $\tilde{V}$ is a continuous extension of $V$ from $H(\curl, \Omega)$ to $H(\curl, \mathcal{B})$. Similarly, we can also write
\begin{align}
    \int_{\gamma'} V \cdot t = -\frac{1}{2\pi}\int_{\partial \Omega} (V \times n) \cdot B^{\text{ext}}, \label{eq:circ_gamma_bound}\\
    \int_{\gamma} V \cdot t' = -\frac{1}{2\pi}\int_{\partial \Omega} (V \times n) \cdot B.
\end{align} 

\begin{remark}\label{rem:sign_B}
    {Although this characterization of the normalized harmonic field makes use of $\Sigma$, choosing another cutting surface may change $B$ only up to a sign. Indeed, $B$ is in $\mathcal{K}(\Omega)$ which is one-dimensional. Therefore, changing $\Sigma$ may change $B$ only by multiplying it by a constant. Furthermore, since $B$ integrates to $2\pi$ along $\gamma'$, we see that this constant is equal to $1$ if we preserve the orientation of $\Sigma$ and $\gamma'$, and equal to $-1$ otherwise.}
\end{remark}

Note that the definition of the normalized harmonic field in \cref{eq:B_curl} is equivalent to the following mixed formulation that will prove to be useful, both for the computation of the shape gradient, and for the numerical scheme.

\begin{proposition}\label{prop:wp_B_div}
    here exists a unique solution to the following problem. Find $(B_\div, u_\div) \in H_0(\div, \Omega) \times L^2_0(\Omega)$ such that, for all $(\tau, v) \in H_0(\div, \Omega) \times L^2(\Omega)$ we have
        \begin{align}
            \begin{split}\label{eq:B_div}
                &{ \int_{\Omega} \left(\div B_\div\right) v} = 0,\\
                &{ \int_{\Omega} B_\div \cdot \tau + \int_{\Omega} u_\div (\div \tau)} = 2\pi\int_{\Sigma}\tau \cdot n_\Sigma.
            \end{split}
        \end{align}
        Furthermore, we have $B_\div = B$ defined in \cref{eq:B_curl}.
\end{proposition}

\begin{proof}
    For the well-posedness, we first notice that this problem is equivalent to the following one. Find $(B_{\div}, u_{\div}, p) \in H_0(\div, \Omega) \times L^2(\Omega) \times \R$ such that 
    \begin{align*}
        &\int_{\Omega} \left(\div B_\div\right) v + p\int_{\Omega} v = 0,
        \\
        &\int_{\Omega} B_\div \cdot \tau + \int_{\Omega} u_\div (\div \tau) = 2\pi \int_{\Sigma}\tau \cdot n_\Sigma,
        \\
        &q\int_{\Omega} u_{\div} = 0,
    \end{align*}
    for all $(\tau, v, q) \in H_0(\div, \Omega) \times L^2(\Omega) \times \R$. Indeed, we find $p=0$ by choosing $v=1$ in the first equation. This problem is then equivalent to the mixed Hodge Laplacian studied in \cite[Section 7]{arnoldFiniteElementExterior2006} with essential boundary conditions and $k=3$. The corresponding bilinear form is then known to satisfy inf-sup conditions \cite[Section 7.5, Remark]{arnoldFiniteElementExterior2006}, and we only need to verify that the right-hand side is continuous with respect to $\tau$ to prove well-posedness. This comes from the continuity of the normal trace from $H(\div, \Omega)$ to $H^{-1/2}(\Sigma)$. As a consequence, we get
    \begin{align*}
        \left| \int_{\Sigma}\tau \cdot n_\Sigma \right|& \leq \|1\|_{H^{1/2}(\Sigma)}\|\tau \cdot n_\Sigma\|_{H^{-1/2}(\Sigma)}\\
        & \leq C \|\tau\|_{H(\div, \Omega)}.
    \end{align*}
    To prove that $B_\div$ and $B$ coincide, we show that $B_\div$ is a harmonic field, and that they are normalized in the same way. Equality then follows from the fact that the set of harmonic fields is one-dimensional. To find that $B_\div$ is a harmonic field, we prove it is orthogonal to $\nabla H^1(\Omega)$ and $\curl H_0(\curl, \Omega)$. 
    \\
    The first fact follows from the first equation of (\ref{eq:B_div}). Indeed, we find that $\div B_\div = 0$, so that $B_\div \in H_0(\div^0, \Omega) = \nabla H^1(\Omega)^{\perp}$.
    \\
    Now, take $\tau = \curl \sigma$ in the second equation of (\ref{eq:B_div}), with $\sigma \in H_0(\curl, \Omega)$. As a consequence, we have $\div \tau = 0$, and
    $$
        \int_{\Sigma}\tau \cdot n_\Sigma = \int_{\gamma} \sigma \cdot t = 0.
    $$
    This gives us ${ \int_{\Omega} B_\div \cdot \tau} = 0$ for all $\tau \in \curl H_0(\curl, \Omega)$. Since $B_\div$ is orthogonal to both $\nabla H^1(\Omega)$ and $\curl H_0(\curl, \Omega)$, we get that it is a harmonic field by Hodge decomposition (see \cref{pr:hodge_decomposition} of \cref{sec:appendix_hodge}).

    Finally, we set $B_\div = \lambda B$, and want to prove that $\lambda=1$. First, by plugging $\tau=B_\div$ in \cref{eq:B_div}, we get
    $$
    \|B_\div\|^2=2\pi\int_{\Sigma}B_\div \cdot n_\Sigma.
    $$
    Now, using the jump condition on $u$ as defined in \cref{eq:B_curl}, and an integration by parts, we get
    \begin{align*}
        \|B\|^2 &= \int_{\Omega \backslash \Sigma} B \cdot \nabla u
        \\
        &=2\pi\int_{\Sigma}B \cdot n_\Sigma.
    \end{align*}
    Now, on the one side we have
    $$
    2\pi\int_\Sigma B_\div \cdot n_\Sigma = 2\pi \lambda\int_\Sigma B \cdot n_\Sigma,
    $$
    and on the other side
    \begin{align*}
        2\pi\int_\Sigma B_\div \cdot n_\Sigma &= \|B_\div\|^2
        \\
        &= \lambda^2 \|B\|^2
        \\
        &= 2\pi \lambda^2 \int_{\Sigma}B \cdot n_\Sigma.
    \end{align*}
    Since $B$ and $B_\div$ are both nonzero, we come to the conclusion that $\lambda=1$.
\end{proof}

\subsection{Vector potentials and Bevir--Gray formula}\label{subsec:VP_BG}

As was mentioned earlier, the numerical computation of the Biot--Savart operator can be very costly. As a consequence, we chose to compute the helicity of the normalized harmonic field by substituting $\mathrm{BS}(B)$ by an appropriate vector potential of $B$. Indeed, since $\curl \mathrm{BS}(B)=B$, we know that any vector potential $A$ of $B$ can only differ from $\mathrm{BS}(B)$ by the sum of a gradient and a harmonic field (see \cref{eq:hodge_decomposition_harm,eq:orth_traceless_curls} in \cref{sec:appendix_hodge}). Since vector fields of $H_0\left(\div^0, \Omega\right)$ are orthogonal to gradient vector fields, we know that ${ \int_{\Omega} B \cdot A}$ can differ from $\mathrm{H}(B)$ only through the harmonic part of the difference between $A$ and $\mathrm{BS}(B)$. However, this difference can be accounted for by modifying the formula of the magnetic helicity, giving a quantity which is invariant under a change of the vector potential. This is given by the well known Bevir--Gray formula \cite{bevirRelaxationFluxConsumption1980} in toroidal domains, which was later generalized to a large class of non-simply connected domains in \cite{mactaggartMagneticHelicityMultiply2019}. In our case, for a vector field $V$ in $H_0\left(\div^0,\Omega\right)$ and $A$ any of its vector potentials, this invariant quantity is given by
\begin{equation}\label{eq:BG}
    \mathrm{H}(V)={ \int_{\Omega} V \cdot A} - \int_{\gamma} A \cdot t \int_{\gamma'} A \cdot t'.
\end{equation}
Note that $\mathrm{BS}(V)$ has zero circulation along $\gamma'$\cite[Section III.A]{cantarellaBiotSavartOperator2001}, so that this invariant quantity does correspond to the usual Biot--Savart helicity.

From this formulation, a natural problem is to find good vector potentials, which are simple to study, both theoretically and numerically. Getting back to the helicity of the normalized harmonic field, we will study two natural vector potential candidates. First, a vector potential given by a classical vector Laplacian problem, which is orthogonal to $B$, therefore canceling the first term of \cref{eq:BG}. Second, a vector potential given in \cite{valliVariationalInterpretationBiot2019}, which is of zero circulation along the toroidal loop $\gamma'$, therefore canceling the second term. As we will see in \cref{sec:num}, using these two vector potentials also allows us to stay in a finite elements exterior calculus setting throughout the numerical computation of the harmonic helicity, and its shape gradient.

The first vector potential, which we denote by $A^1$, is given by the following proposition.
\begin{proposition}\label{pr:wp_vector_potential}
	{Let $B$ be the normalized harmonic field of $\Omega$ as defined in \cref{subsec:harmonic_fields}}. There exists a unique $\left(A^1, u\right) \in H(\curl, \Omega) \times H(\div, \Omega)$ such that, for all $(\tau, v) \in H(\curl, \Omega) \times H(\div, \Omega)$
	\begin{align}
	&{ \int_{\Omega} A^1 \cdot \tau = \int_{\Omega} \curl \tau \cdot u}, \label{eq:A=curlu}
	\\
	&{ \int_{\Omega} \curl A^1 \cdot v + \int_{\Omega} (\div u) (\div v) = \int_{\Omega} B \cdot v}. \label{eq:curlA=B}
	\end{align}
	Furthermore, $A^1$ verifies
	\begin{enumerate}
	\item $A^1$ is in $H_0(\div^0, \Omega)$,\label{item:divA=0}
	\item $\curl A^1 = B$,\label{item:curlA=B}
	\item ${ \int_{\Omega} A^1 \cdot B} = 0$. \label{eq:A_orth_B}
	\end{enumerate}
\end{proposition}

\begin{proof}
    Since the second De Rham cohomology space of $\Omega$ is trivial, the space of harmonic 2-forms of $\Omega$ vanishes, and this variational problem is equivalent to \cite[Equation (7.1)]{arnoldFiniteElementExterior2006} for the case $k=2$. Existence and uniquenesss is then given by \cite[Theorem 7.2]{arnoldFiniteElementExterior2006}. 

    To prove \cref{item:divA=0}, we simply take $\tau$ in $\nabla H^1(\Omega)$ in \cref{eq:A=curlu}. This gives us ${ \int_{\Omega} A^1 \cdot \tau}=0$, so that $A^1$ is orthogonal to $\nabla H^1(\Omega)$, and therefore is in $H_0\left(\div^0, \Omega\right)$ by \cref{eq:orth_gradients} in \cref{sec:appendix_hodge}. To find \cref{item:curlA=B}, we proceed with the same splitting analysis used in \cite{arnoldFiniteElementExterior2006}. We define using the Hodge decomposition \cref{eq:hodge_decomposition_no_harm}, $u=u_\nabla + u_\curl$ with $u_\curl \in \curl H(\curl, \Omega)$ and $u_\nabla \in \nabla H^{1}_{0}(\Omega)$. We now want to prove that $u_\nabla = 0$. Choosing $v \in H\left(\div^0, \Omega\right)^\perp \cap H(\div, \Omega)$ in \cref{eq:curlA=B}, we get
    $$
    { \int_{\Omega} \left(\div u_\nabla\right) (\div v) = \int_{\Omega} B \cdot v} = 0.
    $$
    Using the Poincaré inequality from \cref{pr:continuous_poincare} of \cref{sec:appendix_hodge}, we get $u_\nabla = 0$, so that $\div u = 0$ as claimed.
    \\
    Finally, since $\curl B = 0$, we simply get by choosing $\tau = B$ in \cref{eq:A=curlu} 
    $$
    { \int_{\Omega} A^1 \cdot B} = 0,
    $$
    proving \cref{eq:A_orth_B}.
\end{proof}

For the second vector potential $A^2$, we first need to introduce the following spaces. 
\begin{align*}
	\mathcal{X}(\Omega)&=\left\{V \in H(\curl, \Omega) \mid \curl V \cdot n = 0\right\},\\
	\mathcal{Z}(\Omega)&=\left\{V \in \mathcal{X}(\Omega) \mid \int_{\gamma'} V \cdot t' = 0 \right\}.
\end{align*}
We then introduce the second vector potential $A^2$, given by the following result proven in \cite{valliVariationalInterpretationBiot2019}
\begin{proposition}\label{pr:wp_valli}
    The following problem has a unique solution. Find $(A^2, u) \in \mathcal{Z}(\Omega) \times \nabla H^1(\Omega)$ such that for all $(\tau, v) \in \mathcal{Z}(\Omega) \times \nabla H^1(\Omega)$
    \begin{align}
        \begin{split}\label{eq:A_no_circulation}
            &{ \int_{\Omega} \curl A^2 \cdot \curl \tau + \int_{\Omega} u \cdot \tau = \int_{\Omega} B \cdot \curl \tau},
            \\
            &{ \int_{\Omega} A^2 \cdot v} = 0.
        \end{split}
    \end{align}
    Furthermore, we have $\curl A^2 = B$ and $A^2 \in H_0\left(\div^0, \Omega\right)$.
\end{proposition}

Of course, the two choices of potential vectors are related to each other. Since $A^1$ and $A^2$ are in $H_0\left(\div^0,\Omega\right)$ and have the same curl, we get that $A^1-A^2$ is in $\mathcal{K}(\Omega)$. The fact that $A^2$ has zero circulation along $\gamma'$, and that $B$ has a circulation of $2\pi$ allows us to find the relation 
\begin{align}
    \label{eq_relation_pot_vec_continue}
    A^2 = A^1 - \frac{1}{2\pi}\left(\int_{\gamma'}A^1 \cdot t'\right) B.
\end{align}

\section{Harmonic helicity and its shape derivative}\label{sec:HH_shape_gradient}

As we have seen in the previous section, with each Lipschitz toroidal domain $\Omega$, we are able to associate a normalized harmonic field $B(\Omega)$ and a vector potential of $B(\Omega)$, denoted $A^2(\Omega)$, with zero circulation along $\gamma'$. As a consequence of the Bevir--Gray formula \cref{eq:BG}, the magnetic helicity of $B(\Omega)$ is then given by $\mathrm{H}(B(\Omega))=\int_{\Omega} B(\Omega) \cdot A^2(\Omega)$. In turn, this allows us to define the helicity of the normalized magnetic field of $\Omega$, which we refer to as the harmonic helicity of $\Omega$ for simplicity. As was noted in \cref{rem:sign_B}, $B(\Omega)$ is actually defined as a function of $\Omega$ only up to a sign. However, since the helicity is a quadratic form, this sign indetermination is not relevant, and the harmonic helicity is well-defined as a shape functional.
\begin{definition}\label{def:HH}
    Let $\Omega$ be a Lipschitz toroidal domain. Then, the harmonic helicity of $\Omega$ is defined as 
    $$
    \mathcal{H}(\Omega) = \int_{\Omega} B(\Omega) \cdot A^2(\Omega),
    $$
 where $B(\Omega)$ and $A^2(\Omega)$ are given by the solutions to \cref{eq:B_curl} and \cref{eq:A_no_circulation}.
\end{definition}

The aim of this section is to study how the harmonic helicity varies as a function of the domain $\Omega$. More precisely, we prove that the harmonic helicity is shape Fréchet differentiable under Lipschitz deformation, and we give a formula for its shape derivative.

Before studying the shape differentiability of $\mathcal{H}$, we state the following properties of harmonic helicity. The proof will be given in \cref{subsec:pullbacks_fixed_domain} as we will need to define ways to pullback vector fields beforehand. 

\begin{proposition}\label{prop:vol}
    The following scaling and symmetry properties hold:
    \begin{itemize}
        \item Given $\lambda >0$, we have $\mathcal{H}(\lambda \Omega) =\lambda \mathcal{H}(\Omega)$.
        \item Given a planar reflection $R \in O_3$, we have $\mathcal{H}(R \Omega) =-\mathcal{H}(\Omega)$.
    \end{itemize}
\end{proposition}

\begin{corollary}\label{cor:axisymmetric}
    If $\Omega$ admits a planar symmetry, we have $\mathcal{H}(\Omega) = 0$. In particular, this is the case if $\Omega$ is axisymmetric.
\end{corollary}

{As was mentioned in the introduction, another notion for the helicity of a domain which was studied in \cite{cantarellaIsoperimetricProblemsHelicity2000,cantarella_influence_1999,valliVariationalInterpretationBiot2019,montielIsoperimetricProblemCurl2023} is the maximal value of the helicity on the unit $L^2$ sphere for fields in $H_0\left(\div^0, \Omega\right)$, that is 
$$
\tilde{\mathcal{H}}(\Omega) = \sup \left\{\mathrm{H}(V) \mid V \in H_0\left(\div^0, \Omega\right), \hspace{1ex} \|V\|_{L^2} = 1\right\}.
$$
Since the helicity is quadratic, it is therefore clear that we have the following inequality
$$
\mathcal{H}(\Omega) \leq \|B(\Omega)\|_{L^2}^2 \tilde{\mathcal{H}}(\Omega).
$$
To remove the $L^2$ norm in this inequality, one may normalize the harmonic field using the $L^2$ norm instead of the circulation. Although this makes the comparison between the two quantities more direct, this normalization is less convenient to work with in our context. We note however that one may easily obtain the shape derivative of $\|B(\Omega)\|_{L^2}^2$ using the methods we will introduce in this section, so that this change of normalization may be done without issues for computing the shape derivative of harmonic helicity. We also note that, although $\tilde{\mathcal{H}}$ is known to be bounded when fixing the volume (see \cite[Theorem E]{cantarellaIsoperimetricProblemsHelicity2000}), it is not immediately clear whether this is the case for $\mathcal{H}$. The optimization problems we will consider numerically will be described more precisely in \cref{sec:num}.}

In order to write the shape derivative of the harmonic helicity as a surface integral, we need to assume that $\Omega$ is s-regular for some $s>1/2$ \cite[Section 7.7]{arnoldFiniteElementExterior2006}, that is we have the continuous embedding
\begin{align}
    \label{eq:s_regular}
    H(\curl,\Omega)\cap H_0(\div,\Omega) \xhookrightarrow{} H^s(\Omega).
\end{align}
Indeed, we need $B(\Omega)$ and $A^2(\Omega)$ to have traces in $L^2(\partial \Omega)^3$. For example $\Omega$ being Lipschitz-polyhedral or $\mathcal{C}^{1,1}$ is sufficient~\cite[Prop. 3.7 or Th. 2.9]{amroucheVectorPotentialThreedimensional1998}.

\begin{theorem}\label{th:shape_derivative}
    Let $\Omega$ be a s-regular toroidal domain, with $s> 1/2$, and $\theta$ a vector field in $W^{1, \infty}(\R^3)$ with {$\|D\theta\|_{L^\infty} <1$}. Denoting $\Omega_\theta = (\Id + \theta)(\Omega)$, we have 
    \begin{equation}\label{eq:shape_differentiability}
        \mathcal{H}(\Omega_\theta) = \mathcal{H}(\Omega) + \mathcal{H}'(\Omega;\theta) + o\left(\|\theta\|_{W^{1,\infty}}\right),
    \end{equation}
    where 
    \begin{equation}\label{eq:shape_derivative_formula}
        \mathcal{H}'(\Omega;\theta) = 2\int_{\partial \Omega}{\left( B(\Omega) \cdot A^2(\Omega)\right) \left(\theta \cdot n\right)}.
    \end{equation}
\end{theorem}
In order to prove this theorem, we begin by introducing some transformations in \cref{subsec:pullbacks_fixed_domain}. These transformations, which correspond to usual pullbacks in the language of differential forms, allow us to transform functions and vector fields on the deformed domain $\Omega_\theta$ onto the fixed domain $\Omega$. Since these transformations have good commutation properties with the differential operators, we are then able to use them to differentiate the variational formulations of $B(\Omega_\theta)$ and $A^2(\Omega_\theta)$ in \cref{subsec:differentials_VF}. Once this is done, we can prove \cref{th:shape_derivative} in \cref{subsec:proof_th1} by pulling back the integral defining $\mathcal{H}$ onto a fixed domain, and using the differentiated vector fields from the previous section.

Throughout the rest of this section, $\theta$ denotes a vector field in $W^{1, \infty}(\R^3)$ with {$\|D\theta\|_{L^\infty}<1$}, $\Omega_\theta$ the deformed domain $(\Id + \theta)(\Omega)$, {$D\theta$ denotes the Jacobian matrix of $\theta$ in Cartesian coordinates}, and $J_\theta = \det(\Id + D\theta)$ the determinant of the Jacobian of the transformation. {More generally, given a function (resp. vector field) $u$, we denote by $Du$ its Jacobian matrix in Cartesian coordinates, which is a row vector (resp. square matrix) at each point where it is defined}. To have notations which are less cumbersome, we will also denote $B=B(\Omega)$, $B_\theta = B(\Omega_\theta)$ the solutions to \cref{eq:B_div} in $\Omega$ and $\Omega_\theta$ respectively, and $A=A^2(\Omega)$, $A_\theta=A^2(\Omega_\theta)$ the solutions to \cref{eq:A_no_circulation} in $\Omega$ and $\Omega_\theta$ respectively. We also denote by $n_\theta$ the unit outward pointing vector field on $\partial \Omega_\theta$, and define $\gamma_\theta$, $\gamma'_\theta$, $\Sigma_\theta$ and $\Sigma'_\theta$ in the deformed domain $\Omega_\theta$ as the image by $(\Id+\theta)$ of the corresponding geometrical constructions given in \cref{subsec:harmonic_fields}, with tangent and normal vector fields $t_\theta$, $t'_\theta$, $n_{\Sigma_\theta}$, $n_{\Sigma'_{\theta}}$ respectively.

\subsection{Pullback of the De Rham complex on a fixed domain}\label{subsec:pullbacks_fixed_domain}
In this section, we define maps which allow us to transform elements of our functional spaces on a deformed domain $\Omega_\theta$ back to $\Omega$. These transformations correspond to pullbacks in the language of differential forms, and the corresponding formulas are therefore similar to the Piola mappings used in finite elements exterior calculus. This non-naive treatment is needed as composition by $\Id + \theta$ does not map $H(\curl,\Omega_\theta)$ (resp. $H(\div,\Omega_\theta)$) to $H(\curl,\Omega)$ (resp. $H(\div,\Omega)$).
\begin{definition}\label{def:pullbacks}
    Let $u_0$, $u_1$, $u_2$ and $u_3$ be in $H^1(\Omega_\theta)$, $H(\curl, \Omega_\theta)$, $H(\div, \Omega_\theta)$ and $L^p(\Omega_\theta)$ respectively, with $p$ being 1 or 2. We define 
    \begin{align*}
        \Phi_{\theta}^{0}u_0 &= u_0 \circ (\Id+\theta),
        \\
        \Phi_{\theta}^{1}u_1 &= \left(\Id+D\theta^T\right)u_1 \circ (\Id+\theta),
        \\
        \Phi_{\theta}^{2}u_2 &= J_\theta (\Id + D\theta)^{-1}u_2 \circ (\Id+\theta),
        \\
        \Phi_{\theta}^{3}u_3 &= J_\theta u_3 \circ (\Id+\theta).
    \end{align*}
\end{definition}

Before giving the key properties of these transformations, we start by stating the following useful algebraic identities. These can be found by direct computations.

\begin{lemma}\label{lemma:pullback_identities}
    We have the following identities.
    \begin{align*}
        &\left(\Phi_{\theta}^{0}u\right)\left(\Phi_{\theta}^{3}v\right) = \Phi_{\theta}^{3}(uv) &\forall u \in H^1(\Omega_\theta), v \in L^2(\Omega_\theta),
        \\
        &\left(\Phi_{\theta}^{1}u\right)\cdot\left(\Phi_{\theta}^{2}v\right) = \Phi_{\theta}^{3}(u \cdot v) &\forall u \in H(\curl, \Omega_\theta), v \in H(\div, \Omega_\theta),
        \\
        &(\alpha(\theta)\Phi_{\theta}^{1}u)\cdot\left(\Phi_{\theta}^{1}v\right) = \Phi_{\theta}^{3}(u\cdot v) &\forall u \in H(\curl, \Omega_\theta), v \in H(\curl, \Omega_\theta),
        \\
        &(\alpha(\theta)^{-1}\Phi_{\theta}^{2}u)\cdot\left(\Phi_{\theta}^{2}v\right) = \Phi_{\theta}^{3}(u\cdot v) &\forall u \in H(\div, \Omega_\theta), v \in H(\div, \Omega_\theta),
    \end{align*}
    where $\alpha(\theta)=J_\theta (\Id+D\theta)^{-1}(\Id+D\theta^T)^{-1}$.
\end{lemma}

\begin{proposition}\label{prop:commutativity_pullbacks}
    The diagram
    \[\begin{tikzcd}
        {H^1(\Omega_\theta)} & {H(\curl, \Omega_\theta)} & {H(\div, \Omega_\theta)} & {L^2(\Omega_\theta)} \\
        {H^1(\Omega)} & {H(\curl, \Omega)} & {H(\div, \Omega)} & {L^2(\Omega)}
        \arrow["\nabla", from=1-1, to=1-2]
        \arrow["\curl", from=1-2, to=1-3]
        \arrow["\div", from=1-3, to=1-4]
        \arrow["{\Phi_{\theta}^{0}}"', from=1-1, to=2-1]
        \arrow["{\Phi_{\theta}^{1}}"', from=1-2, to=2-2]
        \arrow["{\Phi_{\theta}^{2}}"', from=1-3, to=2-3]
        \arrow["{\Phi_{\theta}^{3}}"', from=1-4, to=2-4]
        \arrow["\nabla", from=2-1, to=2-2]
        \arrow["\curl", from=2-2, to=2-3]
        \arrow["\div", from=2-3, to=2-4]
    \end{tikzcd}\]
    is commutative. Furthermore, the same diagram can be made with the corresponding traceless spaces.
\end{proposition}

\begin{proof}
    For the proof of the first statement, we choose smooth functions or vector fields for simplicity (both for $u$ and $\theta$). The general cases can then be obtained by density arguments. 
    \\
    We begin by proving $\nabla \circ \Phi_\theta^0 = \Phi_\theta^1 \circ \nabla$. We choose $u$ in $\mathcal{C}^{\infty}(\bar{\Omega}_\theta)$, {and compute using the chain rule
    $$
    D\left(\Phi_{\theta}^{0}u\right) = Du \circ (\Id+\theta) (\Id+D\theta),
    $$
    Since, in Cartesian coordinates, the gradient of a function is the transpose of its Jacobian matrix, we get}
    \begin{align*}
        \nabla \Phi_{\theta}^{0}u &= (\Id+D\theta^T)\nabla u \circ (\Id+\theta)
        \\
        &= \Phi_{\theta}^{1}\nabla u.
    \end{align*}
    Now, we prove that $\curl \circ \Phi_\theta^1 = \Phi_\theta^2 \circ \curl$ by taking $u$ in $\mathcal{C}^{\infty}(\bar{\Omega}_\theta)^3$. We have  
    \begin{align*}
    D\left(\Phi_{\theta}^{1}u\right) &= D\left((\Id+D\theta^T) u \circ (\Id+\theta)\right)
    \\
    &= D\left(D\theta^T\right)u\circ (\Id+\theta)+\left(\Id+D\theta^T\right)Du \circ (\Id+\theta) (\Id+D\theta),
    \end{align*}
    where $D\left(D\theta^T\right)u$ is a symmetric matrix given by
    $$
    \left(D\left(D \theta^T\right)u\right)_{i,j} = \sum_{k=1}^{3}\frac{\partial^2\theta^k}{\partial x^i \partial x^j}u_k,
    $$
    Therefore, we find 
    $$
    D\left(\Phi_{\theta}^{1}u\right)-\left(D\left(\Phi_{\theta}^{1}u\right)\right)^T = \left(\Id+D\theta^T\right)(Du-Du^T) \circ (\Id+\theta) (\Id+D\theta).
    $$
    Now, defining $\mathrm{cr} :\mathfrak{so}(3) \rightarrow \R^3$, where $\mathfrak{so}(3)$ is the space of skew-symmetric matrices, by 
    $$
    \mathrm{cr} 
    \begin{pmatrix}
        0 & -v_3 & v_2 \\
        v_3 & 0 & -v_1 \\
        -v_2 & v_1 & 0
    \end{pmatrix}
    =
    \begin{pmatrix}
        v_1 \\ v_2 \\ v_3
    \end{pmatrix},
    $$
    we can prove that, for any invertible matrix $B$ and skew-symmetric matrix $A$, $\mathrm{cr}(B^T A B) = \det (B) B^{-1} \mathrm{cr}(A)$.
    Since $\curl u = \mathrm{cr}\left(Du - Du^T\right)$, we find 
    \begin{align*}
        \curl \Phi_{\theta}^{1} u &= J_\theta (\Id + D\theta)^{-1}\curl u \circ (\Id + \theta)
        \\
        &= \Phi_{\theta}^{2}\curl u.
    \end{align*}
    For the last commutativity relation $\div \circ \Phi_\theta^2 = \Phi_\theta^3 \circ \div$, we proceed by duality. Take $u$ in $\mathcal{C}^\infty(\bar{\Omega}_\theta)^3$ and $v$ in $\mathcal{C}^{\infty}_{0}(\Omega_\theta)$. First, we notice that the $\Phi_{\theta}^{k}$ are isomorphisms, as $\left(\Phi_{\theta}^{k}\right)^{-1}=\Phi_{(\Id+\theta)^{-1}-\Id}^{k}$. Therefore, 
    \begin{align*}
        \int_{\Omega} \left(\div \Phi_{\theta}^{2}u\right) v &= -\int_{\Omega}\left(\Phi_{\theta}^{2}u\right) \cdot \nabla v
        \\
        &= -\int_{\Omega} \left(\Phi_{\theta}^{2}u\right) \cdot \left[\Phi_{\theta}^{1} \left(\nabla \left(\Phi_{\theta}^{0}\right)^{-1}v\right)\right],
    \end{align*}
    where we used $\Phi_{\theta}^{1}\nabla v = \nabla \Phi_{\theta}^{0} v$. Using \cref{lemma:pullback_identities}, we then find
    \begin{align*}
        \int_{\Omega} \left(\div \Phi_{\theta}^{2}u\right) v &= -\int_{\Omega} \Phi_{\theta}^{3}\left(u \cdot \nabla \left(\Phi_{\theta}^{0}\right)^{-1}v\right)
        \\
        &= -\int_{\Omega_\theta}u \cdot \nabla \left(\Phi_{\theta}^{0}\right)^{-1}v
        \\
        &= \int_{\Omega_\theta} \div u \cdot \left(\Phi_{\theta}^{0}\right)^{-1}v
        \\
        &= \int_{\Omega}\left(\Phi_{\theta}^{3}\div u\right) v.
    \end{align*}
    Now, we prove that traceless functions or vector fields are preserved by the maps $\Phi_{\theta}^{k}$. For $k=0$, we simply take $u \in \mathcal{C}_{0}^{\infty}(\Omega)$, and verify that $u \circ (\Id + \theta)$ is in $H^{1}_{0}(\Omega)$. We then find the desired result by density.
    \\
    The proofs for $k=1$ and $k=2$ being similar, we only write the first one. Take $u$ in $H_{0}(\curl, \Omega_\theta)$ and $v$ in $H(\curl, \Omega)$. Then, we have
    \begin{align*}
        \int_{\partial \Omega} \left(\left(\Phi_{\theta}^{1}u\right) \times n\right) \cdot v &= \int_{\Omega}\curl \left(\Phi_{\theta}^{1} u\right) \cdot v - \int_{\Omega} \left(\Phi_{\theta}^{1}u\right) \cdot \curl v
        \\
        &= \int_{\Omega} \Phi_{\theta}^{2}\left(\curl u\right) \cdot \Phi_{\theta}^{1}\left(\left(\Phi_{\theta}^{1}\right)^{-1}v\right) - \int_{\Omega}\left(\Phi_{\theta}^{1}u\right) \cdot \Phi_{\theta}^{2} \left(\curl \left(\Phi_{\theta}^{1}\right)^{-1}v\right)
        \\
        &=\int_{\Omega_\theta} \curl u \cdot \left(\left(\Phi_{\theta}^{1}\right)^{-1}v\right) - \int_{\Omega_\theta} u \cdot \left(\curl \left(\Phi_{\theta}^{1}\right)^{-1}v\right)
        \\
        &=0,
    \end{align*}
    so that $\Phi_{\theta}^{1}u$ is in $H_0(\curl, \Omega)$ as claimed.
\end{proof}

As the $\Phi_{\theta}^{k}$ transformations correspond to pullbacks of $k$-forms, it is not surprising that $\Phi_{\theta}^{1}$ preserves circulations, and $\Phi_{\theta}^{2}$ preserves fluxes. This result is given by the following lemma.
\begin{lemma}\label{lemma:circ_and_flux_pullback}
    Let $u$ be in $H(\curl, \Omega_\theta)$ with $\curl u \cdot n = 0$, and $v$ be in $H_0(\div, \Omega_\theta)$. We then have the following identities 
    \begin{align*}
        \int_{\gamma'_\theta} u \cdot t'_\theta &= \int_{\gamma'}\Phi_{\theta}^{1} u \cdot t',
        \\
        \int_{\Sigma_\theta} v \cdot n_{\Sigma_\theta} &= \int_{\Sigma} \Phi_{\theta}^{2}v \cdot n_\Sigma.
    \end{align*}
\end{lemma}

\begin{proof}
    For the first equality, we suppose that $u$ is smooth so that we can study the circulation in the usual sense. The general case then follows from a density argument, and the continuity of the circulation as defined in \cite[Section 2]{alonso-rodriguezFiniteElementApproximation2018} with respect to the $H(\curl, \Omega)$ norm. We identify $\gamma'$ with a Lipschitz embedding $\gamma' : S^1 \rightarrow \partial \Omega$, so that the derivative $\dot{\gamma}'(s)$ is defined for a.e. $s \in S^1$. $\gamma'_{\theta}$ is then identified as the embedding $(\Id+\theta)\circ \gamma'$, for which the derivative is given by 
    $$
    \dot{\gamma}'_\theta = \left((\Id + D\theta) \circ \gamma'\right) \dot{\gamma}'.
    $$
    The circulation of $u$ along $\gamma'_\theta$ is then given by 
    $$
    \int_{\gamma'_\theta} u \cdot t'_\theta = \int_{S^1} \left[u \circ \gamma'_\theta\right](s) \cdot \dot{\gamma}'_{\theta}(s)ds.
    $$
    Using the formula for $\gamma_\theta'$ and its derivative as an embedding, we get 
    \begin{align*}
        \int_{\gamma'_\theta} u \cdot t'_\theta &= \int_{S^1} \left[ u \circ (\Id + \theta) \circ \gamma'\right](s) \cdot \left[\left((\Id+D\theta)\circ \gamma'\right) \dot{\gamma}'\right](s)ds
        \\
        &= \int_{S^1} \left[\left(\Id+D\theta^T\right) u \circ(\Id + \theta)\right]\circ \gamma'(s) \cdot \dot{\gamma}'(s) ds
        \\
        &= \int_{\gamma'} \Phi_{\theta}^{1} u \cdot t'.
    \end{align*}
    
    For the second equality, we choose $v \in H_0\left(\div, \Omega_\theta\right)$ and use the following equality given in \cite[Lemma 1]{alonso-rodriguezFiniteElementApproximation2018} for $\psi \in H^1(\Omega_\theta \backslash \Sigma_\theta)$
    \begin{equation}\label{eq:flux_and_jump}
        \int_{\Sigma_\theta} v \cdot n_{\Sigma_\theta} [\![\psi]\!]_{\Sigma_\theta} = \int_{\Omega_\theta \backslash \Sigma_\theta} v \cdot \nabla \psi + \int_{\Omega_\theta \backslash \Sigma_\theta} (\div v) \psi.
    \end{equation}
    Taking $u_\theta$ to be the solution of \cref{eq:B_curl} in $\Omega_\theta$, we get 
    \begin{align*}
        2\pi \int_{\Sigma_\theta} v \cdot n_{\Sigma_\theta} &= \int_{\Omega_\theta \backslash \Sigma_\theta} v \cdot \nabla u_\theta + \int_{\Omega_\theta \backslash \Sigma_\theta} (\div v) u_\theta
        \\
        &= \int_{\Omega \backslash \Sigma} \left(\Phi_{\theta}^{2}v\right)\cdot \nabla \left(\Phi_{\theta}^0 u_\theta\right) + \int_{\Omega \backslash \Sigma} \left(\div \Phi_\theta^2 v\right) \left(\Phi_\theta^0 u_\theta\right)
        \\
        &= \int_{\Sigma} \left(\Phi_\theta^2 v\right) \cdot n_\Sigma \left[\!\left[\Phi_\theta^0 u_\theta\right]\!\right]_\Sigma.
    \end{align*}
    From the definition of $\Phi_\theta^0$, we then see that $\left[\!\left[\Phi_\theta^0 u_\theta\right]\!\right]_\Sigma = [\![u_\theta]\!]_{\Sigma_\theta} \circ (\Id + \theta)=2\pi$, which proves the desired equality.
\end{proof}

We know from \cref{lemma:pullback_identities} that when we pullback products in $H(\curl, \Omega)$ or $H(\div, \Omega)$, there is a non-homogeneous term $\alpha(\theta)$ which appears. We will have to consider such products in the next sections to differentiate the harmonic helicity. The next lemma shows how $\theta \mapsto \alpha(\theta)$ behaves under differentiation when integrated against products of vector fields.

\begin{lemma}\label{lemma:alpha_prime}
    The mapping $\theta \mapsto \alpha(\theta)$ from $W^{1, \infty}\left(\R^3\right)^3$ to $L^\infty\left(\R^3;{\mathcal{M}_3(\R)}\right)$ defined in \cref{lemma:pullback_identities} is $smooth$, and if $u$, $v$ are in $H(\curl, \Omega)\cap H_0(\div, \Omega)$, its differential at $0$ verifies
    $$
    \int_{\Omega}\left(\alpha'(0;\theta)u\right) \cdot v = \int_\Omega \div u \theta \cdot v + \int_\Omega u \cdot \theta \div v+ \int_\Omega u \times \theta \cdot \curl v - \int_\Omega \curl u \times \theta \cdot v + \int_{\partial \Omega} {(u\cdot v) (\theta \cdot n)}.
    $$
\end{lemma}

\begin{proof}
    We recall that $\alpha(\theta)=J_\theta(\Id+D\theta)^{-1}\left(\Id+D\theta^T\right)^{-1}$. To find that $\alpha$ is smooth, we simply notice that the maps defined on the unit ball of $\mathcal{M}_3(\R)$ given by $A \mapsto (\Id+A)^{-1}$, $A \mapsto (\Id+A^T)^{-1}$ and $A \mapsto \det(\Id+A)$ are smooth. We also have that $\theta \mapsto D\theta$ is linear and bounded from $W^{1, \infty}\left(\R^3\right)^3$ to $L^\infty\left(\R^3;{\mathcal{M}_3(\R)}\right)$. We conclude by composition that $\alpha$ is smooth.
    \\
    Using that the differential of the determinant and the inverse at the identity are, respectively, the trace and minus the identity, we find
    $$
    \alpha'(0; \theta)=\div \theta \Id - D\theta - D\theta^T.
    $$
    To prove the final point, we choose $u$ and $v$ in $H^1(\Omega)^3$ such that $u\cdot n = v \cdot n = 0$. Then, $u \cdot \theta$ and $u \times \theta$ are in $H^1$, and we have 
    \begin{align*}
        \nabla(u \cdot \theta) &= Du^T \theta + D\theta^T u,
        \\
        \curl(u \times \theta) &= Du \theta - D\theta u - \div u \theta + \div \theta u.
    \end{align*}
    Combining these formulas, we find
    \begin{align*}
        \nabla(u \cdot \theta)-\curl(u \times \theta) &= -Du \theta + D \theta u + \div u \theta - \div \theta u + Du^T\theta + D\theta^T u
        \\
        &= -\alpha'(0; \theta)u - (Du - Du^T)\theta + \div u \theta
        \\
        &= -\alpha'(0; \theta) - \curl u \times \theta + \div u \theta,
    \end{align*}
    so that 
    $$
    \alpha'(0; \theta)u=\div u \theta + \curl(u \times \theta) - \curl u \times \theta \cdot v - \nabla(u \cdot \theta). 
    $$
    Now, using integration by parts 
    \begin{align*}
        \int_{\Omega}\alpha'(0; \theta) u \cdot v =& \int_\Omega \div u \theta \cdot v + \curl(u \times \theta)\cdot v - \curl u \times \theta \cdot v - \nabla(u\cdot \theta)\cdot v
        \\
        =&\int_\Omega \div u \theta \cdot v + \int_{\Omega}u \times \theta \cdot \curl v - \int_{\partial \Omega} (u\times \theta)\times n \cdot v - \int_\Omega \curl u \times \theta \cdot v + \int_\Omega u \cdot \theta \div v
        \\
        =&\int_\Omega \div u \theta \cdot v + \int_{\Omega}u \times \theta \cdot \curl v + \int_{\partial \Omega} {(u \cdot v) (\theta \cdot n)} - \int_\Omega \curl u \times \theta \cdot v + \int_\Omega u \cdot \theta \div v,
    \end{align*}
    {where we used the formula $(a \times b) \times c = (c \cdot a) b - (c \cdot b) a$ and the fact that $u\in H_0\left(\div, \Omega \right)$ }. Finally, by density, taking any $u$ and $v$ in $H(\curl, \Omega)\cap H_0(\div, \Omega)$, we get the desired formula. Note that we need s-regularity \cref{eq:s_regular} with $s> 1/2$, to ensure that $u,v\in L^2(\partial \Omega)^3$.
\end{proof}
{
Before studying the shape differentiability of harmonic helicity in the next sections, we prove \cref{prop:vol}.

\begin{proof}
    We begin by proving the first point. Let $u(\Omega)$ and $u(\lambda \Omega)$ be the solution to \cref{eq:B_curl} in $\Omega$ and $\lambda \Omega$ respectively. Then, it is clear that $u(\lambda \Omega) = u\left(\cdot/\lambda\right)$. Therefore, we have $B(\lambda \Omega) = \tilde{\nabla} u_\lambda = \frac{1}{\lambda} B(\Omega)(\cdot/\lambda)$. Similarly, we have $A^2(\lambda \Omega) = \frac{1}{\lambda} A(\Omega)(\cdot / \lambda)$. Writing 
    \begin{align*}
        \mathcal{H}(\lambda \Omega) &= \int_{\lambda \Omega} B(\lambda \Omega) \cdot A^2(\lambda \Omega)
        \\
        &= \lambda^3 \int_{\Omega} \left(\frac{1}{\lambda} B(\Omega)\right) \cdot \left(\frac{1}{\lambda} A^2(\Omega)\right) 
        \\
        &= \lambda \mathcal{H}(\Omega)
    \end{align*}
    For the second point, let $R$ be a planar reflection. We introduce as in \cref{def:pullbacks}
    \begin{align*}
        \Phi_R^1 u &= R^T u \circ R,
        \\
        \Phi_R^2 u &= (\det R) R^{-1} u \circ R,
    \end{align*}
    for $u$ a vector field of $R\Omega$. Since $R$ is a planar reflection, we have $\det R = -1$ and $R^{-1} = R^{T}$ so that $\Phi_R^2 u = -\Phi_R^1 u$. Let $u(\Omega)$ and $u(R \Omega)$ be the solutions to \cref{eq:B_curl} in $\Omega$ and $R \Omega$ respectively. Then, since $R$ is an isometry, we have 
    $$
    u(\Omega) = u(R \Omega) \circ R,
    $$
    so that using the definition of the normalized harmonic field and \cref{prop:commutativity_pullbacks}, we have $B(\Omega) = \Phi_R^1 B(R \Omega)$. We now prove that $A^2(\Omega) = \Phi_R^2 A^2(R \Omega)$. First, since $\Phi_R^2 = -\Phi_R^1$, it is clear from \cref{prop:commutativity_pullbacks} that $\Phi_R^2 A^2(R \Omega)$ is in $H_0\left(\div^0, \Omega\right)$. Using once again \cref{prop:commutativity_pullbacks}, we also have
    \begin{align*}
        \curl \Phi_R^2 A^2(R \Omega) &= -\curl \Phi_R^1 A^2(R \Omega)
        \\
        &= -\Phi_R^2 B(R \Omega)
        \\
        &= \Phi_R^1 B(R \Omega)
        \\
        &= B(\Omega).
    \end{align*}
    Moreover, using \cref{lemma:circ_and_flux_pullback}, we know that $\Phi_R^2 A^2(R \Omega) = -\Phi_R^1 A^2(R \Omega)$ is in $\mathcal{Z}(\Omega)$. We then deduce from \cite[Proposition 2]{valliVariationalInterpretationBiot2019} that $\Phi_R^2 A^2(R \Omega) = A^2(\Omega)$. Finally, since $R$ is an isometry, we have 
    \begin{align*}
        \mathcal{H}(R \Omega) &= \int_{R \Omega} B(R \Omega) \cdot A^2 (R \Omega)
        \\
        &= \int_{\Omega} \left(R^T B(R \Omega) \circ R\right) \cdot \left(R^T A^2(R \Omega) \circ R\right),
        \\
        &= -\int_{\Omega} \Phi_R^1 B(R \Omega) \cdot \Phi_R^2 A^2(R \Omega)
        \\
        &= -\mathcal{H}(\Omega).
    \end{align*}
\end{proof}
}

\subsection{Differentiability of the PDEs}\label{subsec:differentials_VF}
In this section we prove the differentiability of $\theta \mapsto \Phi_{\theta}^{2}B_\theta$ and $\theta \mapsto \Phi_{\theta}^{1}A_\theta$. As is classical in shape variation of PDEs (see for example \cite[Chapter 5]{henrotShapeVariationOptimization2018}), this is done using an implicit function argument on a pulled back version of the variational formulations. Although it is quite common in the literature to recover Eulerian derivatives from this, that is, the differential of $\theta \mapsto B_\theta$ and $\theta \mapsto A_\theta$ directly, we chose to skip this step here. Indeed, we will see that this step is not necessary to recover the shape derivative of the harmonic helicity. Furthermore, these Eulerian derivatives often suffer from a loss of regularity compared to the material ones, and satisfy affine variational formulations which are less practical to deal with in our case.
\begin{proposition}\label{prop:B_differentiability}
    Let $(B_\theta, u_\theta)$ be the solution to the following variational problem. 
    \\
    Find $(B_\theta, u_\theta) \in H_0(\div, \Omega_\theta) \times L^2_0(\Omega_\theta)$ such that, for all $(\tau, v) \in H_0(\div, \Omega_\theta) \times L^2_0(\Omega_\theta)$ we have
    \begin{align}
        &\int_{\Omega_\theta} B_\theta \cdot \tau  + \int_{\Omega_\theta} u_\theta \div \tau = 2\pi\int_{\Sigma_\theta}\tau \cdot n_{\Sigma_\theta},\label{eq:vf_b_theta_one}
        \\
        &\int_{\Omega_\theta }\left(\div B_\theta\right) v  = 0.\label{eq:vf_b_theta_two}
    \end{align}
    Then, $\theta \mapsto (\Phi_{\theta}^{2}B_\theta, \Phi_{\theta}^{0}u_\theta)$ is $\mathcal{C}^1$ in a neighborhood of $0$, and its differential at zero $(B', u')$ verifies 
    \begin{gather}\label{eq:B_derivative}
        \begin{split}
            &\int_{\Omega}B' \cdot \tau + \int_\Omega u' \div \tau = \int_\Omega \left(\alpha'(0;\theta)B\right) \cdot \tau,
            \\
            &\int_{\Omega}\left(\div B'\right) v = 0,
        \end{split}
    \end{gather}
    for all $(\tau, v)$ in $H_0(\div, \Omega)\times L^2_0(\Omega)$.
\end{proposition}

\begin{proof}
    We know from \cref{prop:wp_B_div} that $(B_\theta, u_\theta)$ is defined uniquely. Pulling back the integrals of \cref{eq:vf_b_theta_one} onto $\Omega$, and using \cref{lemma:pullback_identities,lemma:circ_and_flux_pullback,lemma:alpha_prime}, we get
    \begin{align*}
        \int_{\Omega}\Phi_{\theta}^{3}(B_\theta \cdot \tau) + \int_{\Omega}\Phi_{\theta}^{3}(u_\theta \div \tau) &= 2\pi \int_{\Sigma}\left(\Phi_{\theta}^{2}\tau\right) \cdot n_{\Sigma},
        \\
        \int_{\Omega} \left(\alpha(\theta)^{-1}\Phi_{\theta}^{2}B_\theta\right)\cdot \left(\Phi_{\theta}^{2}\tau\right) + \int_{\Omega}\left(\Phi_{\theta}^{0}u_\theta\right) \left(\div \Phi_{\theta}^{2}\tau\right) &= 2\pi \int_{\Sigma}\left(\Phi_{\theta}^{2}\tau\right) \cdot n_{\Sigma}.
    \end{align*}
    Similarly, we get from \cref{eq:vf_b_theta_two}
    $$
    \int_{\Omega}\left(\div \Phi_{\theta}^{2}B_\theta\right) \left(\Phi_{\theta}^{0}v\right) = 0.
    $$
    Of course, since the $\Phi_{\theta}^{k}$ define isomorphisms, we can take test functions $(\tau, v)$ in $H_0(\div, \Omega)\times L^2_0(\Omega)$.
    We now define 
    $$
    F: W^{1, \infty}(\R^3)^3 \times \left(H_0(\div, \Omega)\times L^2_0(\Omega)\right) \rightarrow \left(H_0(\div, \Omega)\times L^2_0(\Omega)\right)'
    $$
    by 
    $$
    F(\theta; \sigma, u)(\tau, v) = \int_{\Omega} \left(\alpha(\theta)^{-1}\sigma\right)\cdot \tau + \int_{\Omega}u \div \tau + \int_{\Omega}\left(\div \sigma\right) v - 2\pi \int_{\Sigma}\tau \cdot n_\Sigma,
    $$
    so that $(B_\theta, u_\theta)$ solving \cref{eq:vf_b_theta_one,eq:vf_b_theta_two} is equivalent to $F(\theta; \Phi_{\theta}^{2}B_\theta, \Phi_{\theta}^{0}u_\theta)=0$.

    Now, we know from \cref{lemma:alpha_prime} that $\theta \mapsto \alpha(\theta)^{-1}$ is $\mathcal{C}^{1}$. Since $(\sigma, u) \mapsto F(\theta; \sigma, u)$ is linear and bounded, we deduce that $F$ is $\mathcal{C}^{1}$. Furthermore, denoting by $D_{\sigma, u}F$ the differential of $F$ with respect to the $(\sigma, u)$ variables,
    $$
    D_{\sigma, u}F(0; B_0, u_0)(\sigma',u')(\tau, v) = \int_{\Omega}\sigma' \cdot \tau + \int_{\Omega}u' \div \tau + \int_{\Omega}\left(\div \sigma'\right) v.
    $$
    We therefore find that $D_{\sigma, u}F(0; B_0, u_0)$ is an isomorphism from $H_0(\div, \Omega)\times L^2_0(\Omega)$ to $\left(H_0(\div, \Omega)\times L^2_0(\Omega)\right)'$ by the same inf-sup inequalities used to prove the well-posedness of \cref{prop:wp_B_div}. Using the implicit function theorem, we deduce that for small enough $\theta$, there is a unique $C^{1}$ mapping $\theta \mapsto (\sigma(\theta), u(\theta))$ such that $F(\theta; \sigma(\theta), u(\theta))=0$. From uniqueness, we find that $(\sigma(\theta), u(\theta))=(\Phi_{\theta}^{2}B_\theta, \Phi_{\theta}^{0}u_\theta)$.

    Now, to get a variational formulation, we simply differentiate $F(\theta; \Phi_{\theta}^{2}B_\theta, \Phi_{\theta}^{0}u_\theta)=0$ at $\theta=0$. Denoting $\Phi_{\theta}^{2}B_\theta = B + B' + o(||\theta||_{C^{1, 1}})$ and $\Phi_{\theta}^{0}u_\theta = u + u' + o(||\theta||_{C^{1, 1}})$, this gives us
    $$
    \int_{\Omega}B' \cdot \tau + \int_{\Omega}u' \div \tau + \int_{\Omega}\left(\div B'\right) v = \int_{\Omega}(\alpha'(0;\theta)B_0)\cdot \tau,
    $$
    which concludes the proof.
\end{proof}

\begin{proposition}\label{prop:A_differentiability}
    Let $(A_\theta, u_\theta)$ be the solution to the following variational problem.
    \\
    Find $(A_\theta, u_\theta)$ in $\mathcal{Z}(\Omega_\theta)\times \nabla H^1(\Omega_\theta)$ such that, for all $(\tau, v)$ in $\mathcal{Z}(\Omega_\theta)\times \nabla H^1(\Omega_\theta)$ we have 
    \begin{align}
        &\int_{\Omega_\theta}\curl A_\theta \cdot \curl \tau + \int_{\Omega_\theta}u_\theta \cdot \tau = \int_{\Omega_\theta} B_\theta \cdot \curl \tau,\label{eq:vf_a_theta_one}
        \\
        &\int_{\Omega_\theta}A_\theta \cdot v = 0.\label{eq:vf_a_theta_two}
    \end{align}
    Then, $\theta \mapsto (\Phi_{\theta}^{1}A_\theta, \Phi_{\theta}^{2}u_\theta)$ is $\mathcal{C}^1$ in a neighborhood of $0$.
\end{proposition}

\begin{proof}
    We proceed with an implicit function theorem argument, similar to the one used for \cref{prop:B_differentiability}. First, we note that from \cref{prop:commutativity_pullbacks,lemma:circ_and_flux_pullback}, we get $\Phi_{\theta}^{1}\mathcal{Z}(\Omega_\theta)=\mathcal{Z}(\Omega)$, so that the functional spaces of \cref{eq:vf_a_theta_one,eq:vf_a_theta_two} are preserved by the pullbacks.

    Pulling back \cref{eq:vf_a_theta_one,eq:vf_a_theta_two} onto $\Omega$ and using \cref{lemma:pullback_identities}, we get 
    \begin{align*}
        &\int_{\Omega}\left(\alpha(\theta)^{-1}\curl \Phi_{\theta}^{1}A_\theta\right) \cdot \left(\curl \Phi_{\theta}^{1} \tau\right) + \int_{\Omega} \left(\alpha(\theta)\Phi_{\theta}^{1}u_\theta\right)\cdot \left(\Phi_{\theta}^{1}\tau\right) = \int_{\Omega} \left(\alpha(\theta)^{-1} \Phi_{\theta}^{2}B_\theta\right) \cdot \left(\curl \Phi_{\theta}^{1} \tau\right),
        \\
        &\int_{\Omega}\left(\alpha(\theta)\Phi_{\theta}^{1}A_\theta\right) \cdot \left(\Phi_{\theta}^{1}v\right) = 0.
    \end{align*}
    We therefore define 
    $$
    G : W^{1, \infty}(\R^3)^3 \times \left(\mathcal{Z}(\Omega)\times \nabla H^{1}(\Omega)\right) \rightarrow \left(\mathcal{Z}(\Omega)\times \nabla H^{1}(\Omega)\right)'
    $$
    by 
    $$
    G(\theta; \sigma, u)(\tau, v)=\int_{\Omega}\left(\alpha(\theta)^{-1}\curl \sigma\right) \cdot \curl \tau + \int_{\Omega} \left(\alpha(\theta)u\right) \cdot \tau + \int_{\Omega}\left(\alpha(\theta)\sigma\right)\cdot v - \int_{\Omega}\left(\alpha(\theta)^{-1}\Phi_{\theta}^{2}B_\theta\right) \cdot \curl \tau,
    $$
    so that $(A_\theta, u_\theta)$ solves \cref{eq:vf_a_theta_one,eq:vf_a_theta_two} if and only if $G(\theta; \Phi_{\theta}^{1}A_\theta, \Phi_{\theta}^{1}u_\theta)=0$.

    By \cref{lemma:alpha_prime}, we know that $\theta \mapsto \alpha(\theta)$ and $\theta \mapsto \alpha(\theta)^{-1}$ are $\mathcal{C}^{1}$. Furthermore, we know from \cref{prop:B_differentiability} that $\theta \mapsto \Phi_{\theta}^{2}B_\theta$ is $\mathcal{C}^{1}$. Therefore, by linearity and continuity of $G$ with respect to $(\sigma, u)$, we know that $G$ is $\mathcal{C}^{1}$. We have 
    $$
    D_{\sigma, u}G(0; A_0, v_0)(\sigma',u')(\tau,v)=\int_{\Omega}\curl \sigma' \cdot \curl \tau + \int_{\Omega} u' \cdot \tau + \int_{\Omega}\sigma' \cdot v,
    $$
    so that $D_{\sigma, u}G(0; A_0, v_0)$ is an isomorphism by the inf-sup conditions proven in \cite[Section IV]{valliVariationalInterpretationBiot2019}. This proves, by the implicit function theorem, that for $\theta$ small enough there is a unique mapping $\theta \mapsto (\sigma(\theta), u(\theta))$ such that $G(\theta; \sigma(\theta), u(\theta))=0$. By uniqueness, we get $(\sigma(\theta), u(\theta))=(\Phi_{\theta}^{1}A_\theta, \Phi_{\theta}^{1}u_\theta)$.
\end{proof}

\subsection{Proof of \cref{th:shape_derivative}}\label{subsec:proof_th1}
Now that we have introduced ways to pullback functions and vector fields onto the fixed domain $\Omega$, and that we have derived differentiability of $\Phi_\theta^2 B_\theta$ and $\Phi_\theta^1 A_\theta$, we can conclude the proof of \cref{th:shape_derivative}. To do so, we simply pullback the integral of $B_\theta$ against $A_\theta$ onto $\Omega$ using \cref{lemma:pullback_identities}, and use the differentiability results from the last section.
\begin{proof}
    We have
    \begin{align*}
        \mathcal{H}(\Omega_\theta)&=\int_{\Omega_\theta} B_\theta \cdot A_\theta
        \\
        &=\int_{\Omega}\Phi_{\theta}^{3}(B_\theta \cdot A_\theta)
        \\
        &=\int_{\Omega}\left(\Phi_{\theta}^{2}B_\theta\right) \cdot \left(\Phi_{\theta}^{1}A_\theta\right).
    \end{align*}
    From \cref{prop:B_differentiability,prop:A_differentiability}, we know that $\Phi_{\theta}^{2}B_\theta$ and $\Phi_{\theta}^{1}A_\theta$ are differentiable at zero, so that $\mathcal{H}$ is differentiable at $\Omega$ and, denoting their differentials in the direction $\theta$ as $B'$ and $A'$, we have
    $$
    \mathcal{H}'(\Omega;\theta)=\int_{\Omega}B' \cdot A + \int_{\Omega}B \cdot A'.
    $$
    Since $B=\curl A$, and $A$, $A'$ are both in $\mathcal{Z}(\Omega)$, we know from \cite[Lemma 7]{alonso-rodriguezFiniteElementApproximation2018} that 
    $$
    \int_{\Omega}B \cdot A' = \int_{\Omega} \curl A \cdot A' = \int_{\Omega} \curl A' \cdot A.
    $$
    Furthermore, since $\curl A_\theta = B_\theta$, we get $\curl \Phi_{\theta}^{1}A_\theta = \Phi_{\theta}^{2}B_\theta$. By differentiating, we get $\curl A' = B'$, so that 
    $$
    \mathcal{H}'(\Omega;\theta)=2\int_{\Omega}B' \cdot A.
    $$
    Using the first equation in \cref{eq:B_derivative}, we get 
    \begin{align*}
        \int_{\Omega}B' \cdot A = \int_{\Omega}\left(\alpha'(0;\theta)B\right) \cdot A.
    \end{align*}
    Finally, since $B$ and $A$ are in $H(\curl, \Omega)\cap H_0(\div, \Omega)$, we can use \cref{lemma:alpha_prime} to get 
    \begin{align*}
        \int_{\Omega}\left(\alpha'(0;\theta)B'\right) \cdot A &= \int_\Omega \div B \cdot A + \int_\Omega B \times \theta \cdot \curl A - \int_\Omega \curl B \times \theta \cdot A + \int_\Omega B \cdot \theta \div A + \int_{\partial \Omega}(B \cdot A) (\theta \cdot n)
        \\
        &= \int_{\partial \Omega} (B \cdot A) (\theta \cdot n).
    \end{align*}
\end{proof}

\section{Approximation by finite element exterior calculus}
\label{sec:approx_fecc}
In this section we assume that $\Omega$ is a toroidal domain with  polyhedral boundary and that $(\mathcal{T}_h)_{h>0}$ is a quasi-uniform family of tetrahedron meshes of $\bar{\Omega}$, with $h$ the largest diameter of the cells. From \cite[Proposition 3.7]{amroucheVectorPotentialThreedimensional1998}, we know that we can choose $s>1/2$ so that $\Omega$ is $s$-regular. We denote by $\Delta^0(\mathcal{T}_h)$, $\Delta^1(\mathcal{T}_h)$, $\Delta^2(\mathcal{T}_h)$ and $\Delta^3(\mathcal{T}_h)$ the sets of points, edges, faces and cells of the mesh $\mathcal{T}_h$, respectively. For a $k$-simplex $S$ of the mesh $\mathcal{T}_h$, we define $\Delta^i(S)$ as the set of $i$-simplices which have non-empty intersection with $\bar{S}$ for $i < k$. {Moreover, we assume that the cutting surface $\Sigma$ and toroidal curve $\gamma'$ introduced in \cref{subsec:harmonic_fields} are compatible with the mesh, that is, are given by unions of elements of $\Delta^2(\mathcal{T}_h)$ and $\Delta^1(\mathcal{T}_h)$ respectively.} Throughout this section, inequality constants denoted by $C$ are independent of $h$, but may depend on the domain $\Omega$.

Our goal is to provide a finite elements scheme to compute both the harmonic helicity of the domain $\Omega$ and its shape gradient. In \cref{sub-sec:recall-FEM}, we recall the main properties of classical families of elements coming from finite elements exterior calculus. Then, we recall in \cref{sub-sec:Discrete_De_Rham} some notions related to the discrete version of the De Rham complex, and the related discrete harmonic fields. Finally, in \cref{sub-sec:convergence_hh} we prove the convergence of the numerical harmonic helicity.

\subsection{Classical finite elements exterior calculus families}
\label{sub-sec:recall-FEM}
Here, we introduce some classical families in finite elements exterior calculus. The main idea is to define a discretized version of each functional space in the De Rham complex, and stable projections onto these spaces which commute with the differential operators. All these notions were first introduced in \cite{raviartMixedFiniteElement1977,nedelecMixedFiniteElements1980}, and later generalized for differential forms in \cite{arnoldFiniteElementExterior2006}.
\\
For $r$ an integer and $T$ a tetrahedral domain in $\R^3$, we define $\mathcal{P}_r(T)$ as the set of polynomials on $T$ with degree at most $r$, and $\tilde{\mathcal{P}}_r(T)$ the set of homogeneous polynomials on $T$ of degree $r$. {Then, $\mathcal{P}_{r}^{-, \curl}(T)$ and $\mathcal{P}_r^{-, \div}$ are defined as 
\begin{align*}
    \mathcal{P}_{r}^{-, \curl}(T) &= \mathcal{P}_{r-1}(T)^3 \oplus \left\{p \in\tilde{\mathcal{P}}_r(T)^3 \mid p \cdot \mathbf{x} = 0 \right\},
    \\
    \mathcal{P}_{r}^{-, \div}(T) &= \mathcal{P}_{r-1}(T)^3 \oplus \left\{\mathbf{x}p \mid p \in \tilde{\mathcal{P}}_{r-1}(T)\right\}.
\end{align*}}

For $r$ positive, we then define the discretizations
\begin{align*}
    V_{h}^{0}(\Omega)&=\left\{ u \in H^{1}(\Omega) \mid u_{\mid K} \in \mathcal{P}_r(T) \hspace{1ex} \forall T \in \mathcal{T}_h\right\},
    \\
    V_{h}^{1}(\Omega)&=\left\{ u \in H(\curl, \Omega) \mid u_{\mid K} \in {\mathcal{P}_{r}^{-, \curl}(T)} \hspace{1ex} \forall T \in \mathcal{T}_h\right\},
    \\
    V_{h}^{2}(\Omega)&=\left\{ u \in H(\div, \Omega) \mid u_{\mid K} \in {\mathcal{P}_{r}^{-, \div}(T)} \hspace{1ex} \forall T \in \mathcal{T}_h\right\},
    \\
    V_{h}^{3}(\Omega)&=\left\{ u \in L^2(\Omega) \mid u_{\mid K} \in \mathcal{P}_{r-1}(T) \hspace{1ex} \forall T \in \mathcal{T}_h\right\}.
\end{align*}
We then obtain the following sequence
\[\begin{tikzcd}
    {V_{h}^{0}(\Omega)} & {V_{h}^{1}(\Omega)} & {V_{h}^{2}(\Omega)} & {V_{h}^{3}(\Omega)}
    \arrow["\nabla", from=1-1, to=1-2]
    \arrow["\curl", from=1-2, to=1-3]
    \arrow["\div", from=1-3, to=1-4]
\end{tikzcd}.\]
These finite element spaces correspond respectively to the Lagrange $P^r$ elements, the Nedelec first kind elements of order $r$, the Raviart Thomas elements of order $r$, and the {discontinuous $P^{r-1}$ elements}.

One can check easily, using integration by parts, that $u$ is in $V_{h}^{1}(\Omega)$ if and only if the tangential trace of $u$ is continuous along all shared faces of $\Omega$. More precisely, if $T_1$ and $T_2$ are in $\Delta^3(\mathcal{T}_h)$ and $S_1$, $S_2$ are in $\Delta^2(T_1)$ and $\Delta^2(T_2)$ respectively with $T_1 \cap T_2 = S_1 = S_2$, we have 
$$
u \times n_{S_1} + u \times n_{S_2} = 0.
$$
Similarly, we have that $u$ is in $V_{h}^{2}(\Omega)$ if and only the normal trace of $u$ is continuous along all shared faces of $\Omega$.

We also introduce the discrete affine space $V_h^\text{aff}$ to solve \cref{eq:B_curl} numerically. It is defined by 
$$
V_h^\text{aff}(\Omega)= \left\{ v \in H^1(\Omega \backslash \Sigma) \mid v_{|T} \in \mathcal{P}_r(T) \hspace{1ex} \forall T \in \mathcal{T}_h, \text{ and } [\![v]\!]_\Sigma = 2\pi \right\},
$$
where $[\![v]\!]_\Sigma$ denotes the jump of $v$ across $\Sigma$. The corresponding linear space is $V_h^0(\Omega)$.

We also use the smoothed quasi interpolators built in \cite[Section 5.4.]{arnoldFiniteElementExterior2006} for differential forms, denoted by $\Pi_{h}^{k}$. It is known that these quasi interpolators are stable, and that they make the following diagram commute
\[\begin{tikzcd}
    {H^1(\Omega)} & {H(\curl, \Omega)} & {H(\div, \Omega)} & {L^2(\Omega)} \\
    {V_{h}^{0}(\Omega)} & {V_{h}^{1}(\Omega)} & {V_{h}^{2}(\Omega)} & {V_{h}^{3}(\Omega)}
    \arrow["\nabla", from=1-1, to=1-2]
    \arrow["\curl", from=1-2, to=1-3]
    \arrow["\div", from=1-3, to=1-4]
    \arrow["{\Pi_{h}^{0}}"', from=1-1, to=2-1]
    \arrow["{\Pi_{h}^{1}}"', from=1-2, to=2-2]
    \arrow["{\Pi_{h}^{2}}"', from=1-3, to=2-3]
    \arrow["{\Pi_{h}^{3}}"', from=1-4, to=2-4]
    \arrow["\nabla", from=2-1, to=2-2]
    \arrow["\curl", from=2-2, to=2-3]
    \arrow["\div", from=2-3, to=2-4]
\end{tikzcd}.\]
We then have the following approximation estimates on the quasi interpolators (see for example \cite[Theorem 2.2. and 2.3.]{ern_analysis_2018}).
\begin{proposition}\label{prop:approx_quasi_interpolator}
    We have, for all $k$ between $0$ and $3$ and $u$ $H^s$ regular,
    $$
    \left\|u-{\Pi}_{h}^{k}u\right\|_{L^2} \leq C h^s \left\|u\right\|_{H^s}.
    $$
\end{proposition}

\subsection{Discretization of the De Rham complex}
\label{sub-sec:Discrete_De_Rham}
In this section we state some results about the discretization of the De Rham complex. We borrow here most of our notations and lemmas from \cite{arnoldFiniteElementExterior2006}. We begin by defining discrete equivalents of closed and exact fields, which allow us to define the discrete harmonic fields. This then allows us to derive Hodge decompositions in the discrete setting, and uniform Poincaré inequalities. We then state some lemmas which will be useful for the coming convergence results.

First, we use some notations from the differential forms setting to unify some definitions and results. We denote 
$$
H\Lambda^0(\Omega)=H^1(\Omega), \hspace{1ex}H\Lambda^1(\Omega)=H(\curl, \Omega), \hspace{1ex} H\Lambda^2(\Omega)=H(\div, \Omega), \hspace{1ex}\text{and } H\Lambda^3(\Omega)=L^2(\Omega),
$$
as well as
$$
\d^0 = \nabla, \hspace{1ex} \d^1 = \curl, \hspace{1ex} \text{and } \d^2 = \div.
$$
We also denote $d^k=0$ when $k$ is negative or larger than 3. The corresponding traceless spaces are denoted $\mathring{H}\Lambda^k(\Omega)$, and the discrete traceless spaces are given by $\mathring{V}_{h}^{k}(\Omega)=V_{h}^{k}(\Omega) \cap \mathring{H}\Lambda^k(\Omega)$.

We are now able to define the discrete harmonic fields. We denote 
$$
\mathcal{B}_{h}^{k}(\Omega)=d^{k-1}V_{h}^{k-1}(\Omega), \hspace{1ex} \mathcal{Z}_{h}^{k}(\Omega)=\left\{ u \in V_{h}^{k}(\Omega) \mid \d^k u = 0 \right\},
$$
and 
$$
\mathcal{K}_{h}^{k}(\Omega) = \mathcal{Z}_{h}^{k}(\Omega) \cap \mathcal{B}_{h}^{k}(\Omega)^\perp.
$$
Similarly, we define for traceless spaces 
$$
\mathring{\mathcal{B}}_{h}^{k}(\Omega)=d^{k-1}\mathring{V}_{h}^{k-1}(\Omega), \hspace{1ex} \mathring{\mathcal{Z}}_{h}^{k}(\Omega)=\left\{ u \in \mathring{V}_{h}^{k}(\Omega) \mid \d^k u = 0 \right\},
$$
and 
$$
\mathring{\mathcal{K}}_{h}^{k}(\Omega) = \mathring{\mathcal{Z}}_{h}^{k}(\Omega) \cap \mathring{\mathcal{B}}_{h}^{k}(\Omega)^\perp.
$$
From the equivalence of discrete and continuous De Rham cohomology (see e.g. \cite[Section 5.6]{arnoldFiniteElementExterior2010}), we get 
$$
\mathcal{K}_{h}^{0}(\Omega)\cong \R, \hspace{1ex}\mathcal{K}_{h}^{1}(\Omega)\cong \R, \hspace{1ex}\mathcal{K}_{h}^{2}(\Omega)\cong 0 \hspace{1ex} \text{and }\mathcal{K}_{h}^{3}(\Omega)\cong 0.
$$
Similarly, from the equivalence of discrete and continuous De Rham cohomology with boundary condition, and Poincaré duality, we get 
$$
\mathring{\mathcal{K}}_{h}^{0}(\Omega)\cong 0, \hspace{1ex}\mathring{\mathcal{K}}_{h}^{1}(\Omega)\cong 0, \hspace{1ex}\mathring{\mathcal{K}}_{h}^{2}(\Omega)\cong \R \hspace{1ex} \text{and } \mathring{\mathcal{K}}_{h}^{3}(\Omega)\cong \R.
$$
\begin{remark}
    In the continuous case, we defined the harmonic fields as
    $
    \mathcal{K}(\Omega)=H\left(\curl^0, \Omega\right) \cap H_0\left(\div^0, \Omega\right).
    $
    By using $H_0\left(\div^0, \Omega\right) = \nabla H^1(\Omega)^\perp$ and $H\left(\curl^0, \Omega\right) = \curl H_0\left(\curl, \Omega\right)^\perp$ from \cref{pr:orthogonality_relations} of \cref{sec:appendix_hodge}, we recover two expressions for $\mathcal{K}(\Omega)$ which are similar to the discrete ones. However, in the discrete case, the two spaces $\mathcal{K}_{h}^{1}(\Omega)$ and $\mathring{\mathcal{K}}_{h}^{2}(\Omega)$ do not coincide. {Indeed, elements of $\mathcal{K}_h^1(\Omega)$ are curl free but only weakly divergence free and tangent to the boundary meaning that they are not orthogonal to $\nabla H^1(\Omega)$ but only to the finite dimensional space $\nabla V_h^0(\Omega)$. On the other hand, elements of $\mathring{\mathcal{K}}_h^2(\Omega)$ are divergence free and tangent to the boundary, but only weakly curl free, that is, not orthogonal to $\curl H_0\left(\curl^0, \Omega\right)$ but to $\curl \mathring{V}_h^1(\Omega)$.}
\end{remark}
From these definitions, it is now straightforward to find the discrete Hodge decompositions.
\begin{align*}
    V_{h}^{k}(\Omega)&=\mathcal{Z}_{h}^{k}(\Omega)^\perp \oplus\mathcal{K}_{h}^{k}(\Omega)\oplus\mathcal{B}_{h}^{k}(\Omega),
    \\
    \mathring{V}_{h}^{k}(\Omega)&=\mathring{\mathcal{Z}}_{h}^{k}(\Omega)^\perp \oplus\mathring{\mathcal{K}}_{h}^{k}(\Omega)\oplus\mathring{\mathcal{B}}_{h}^{k}(\Omega).
\end{align*}

There is also an equivalent of the Poincaré inequality in the discrete case given in \cite{arnoldFiniteElementExterior2006}[Lemma 5.11].

\begin{proposition}\label{prop:discrete_poincare}
    There exists a constant $C$, independent of $h$, such that for all $u$ in $\mathcal{Z}_{h}^{k}(\Omega)^\perp$
    $$
    \|u\|_{L^2} \leq C\left\|d^k u\right\|_{L^2}.
    $$  
\end{proposition}

We will also be using the following lemma for the proof of convergence of the harmonic field. This lemma is proven in \cite{arnoldFiniteElementExterior2006}[Lemma 5.9] for harmonic forms.
\begin{lemma}\label{lemma:discrete_harmonic_fields}
    For all $u_h$ in $\mathcal{K}_{h}^{1}(\Omega)$ there exists $u$ in $\mathcal{K}(\Omega)$ such that $\left\| u \right\| \leq \left\| u_h \right\|$ and 
    $$
    \left\| u_h - u \right\|_{L^2} \leq \left\| u - \Pi_{h}^{1}u \right\|_{L^2}.
    $$
\end{lemma}

\subsection{Numerical convergence of the harmonic helicity}
\label{sub-sec:convergence_hh}

We begin by studying the approximation of the harmonic field. As we have seen in \cref{subsec:harmonic_fields}, there are two different variational formulations for the harmonic fields. Although they give the same fields in the continuous case, this will not be true at the discrete level. As we will see, the classical Poisson formulation will give a discrete harmonic field in $\mathcal{K}_{h}^{1}(\Omega)$, and the mixed formulation will give a discrete harmonic field in $\mathring{\mathcal{K}}_h^{2}(\Omega)$.

Before studying the convergence of the numerical solutions, we state their well-posedness.
\begin{proposition}\label{prop:WP_B_h_curl}
    There exists a unique solution to the following variational formulation. Find $u_h \in V_{h}^{\rm{aff}}(\Omega)$ such that, for all $v_h \in V_{h}^{0}(\Omega)$,
    \begin{equation}\label{eq:VF_B_h_curl}     
        { \int_{\Omega} u_h \cdot \nabla v_h} = 0.
    \end{equation}
    Furthermore, $B_{h}^{\curl}=\nabla u_h$ is in $\mathcal{K}_{h}^{1}(\Omega)$.
\end{proposition}

\begin{proof}
    We denote by $u$ the solution to \cref{eq:B_curl}, and $B=\nabla u$ the normalized harmonic field of $\Omega$. We begin by noticing that there is a unique solution to the following variational problem. Find $\tilde{u}_h \in V_{h}^{0}(\Omega)$ such that for all $v_h \in V_{h}^{0}(\Omega)$
    \begin{equation}\label{eq:modified_VF_B_h_curl}
        { \int_{\Omega} \nabla \tilde{u}_h \cdot \nabla v_h = - \int_{\Omega} \Pi_{h}^{1}B \cdot v_h}.
    \end{equation}
    Indeed, this problem is a classical discretization of a Poisson equation, so  the Lax--Milgam theorem applies.
    \\
    Now, define $u_h = \tilde{u}_h + \Pi_{h}^{0} u$. Since $\nabla u_h = \nabla \tilde{u}_h + \Pi_{h}^{1} B$, and $\Pi_{h}^{0}u$ is in $V_{h}^{\text{aff}}(\Omega)$, we notice that $u_h$ is solution to \cref{eq:VF_B_h_curl} if and only if $\tilde{u}_h$ is solution to \cref{eq:modified_VF_B_h_curl}. As a consequence, \cref{eq:VF_B_h_curl} is also well-posed.

    We now prove that $\nabla u_h$ is in $V_{h}^{1}(\Omega)$. It is straightforward that interelement continuity will be verified on all surfaces which are not included in $\Sigma_h$. Therefore, we take $S \in \Delta^{2}(\mathcal{T}_h)$ included in $\Sigma_h$. We know that there exist two cells $T_1$ and $T_2$ in $\Delta^3(\mathcal{T}_h)$ such that $S \in \Delta^2(T_1) \cap \Delta^2(T_2)$. We then denote $S=S_1$ when seen as an element of $\Delta^2(T_1)$, and $S=S_2$ when seen as an element of $\Delta^2(T_2)$, and order $T_1$ and $T_2$ so that the exterior normal of $T_1$ on $S_1$ is $n_{\Sigma_h}$. Since $u_h$ is in $V_{h}^{\text{aff}}(\Omega)$, we have $\left.{u_{h}}\right|_{S_2}-\left.{u_{h}}\right|_{S_1}=2 \pi$, so that $\left.\nabla u_h \times n_{\Sigma_h}\right|_{S_2}=\left.\nabla u_h \times n_{\Sigma_h}\right|_{S_1}$. As a consequence, interelement continuity is verified across $S$, so that $B_{h}^{\curl}=\nabla u_h$ is in $V_{h}^{1}(\Omega)$.

    Finally, we prove that $B_{h}^{\curl}$ is in $\mathcal{K}_{h}^{1}(\Omega)$. The fact that $B_{h}^{\curl}$ is in $\mathcal{B}_{h}^{1}(\Omega)^\perp$ is given directly by \cref{eq:VF_B_h_curl}. Now, we take $T$ in $\Delta^3(\mathcal{T}_h)$. Since $\left.{B_{h}^{\curl}}\right|_{T}=\left. {\nabla u_h}\right|_{T}$, we have $\left.{B_{h}^{\curl}}\right|_{T} \in \mathcal{B}_{h}^{1}(T)$. 
    \\
    Now, using that $\mathcal{B}_{h}^{1}(T)$ is a subset of $\mathcal{Z}_{h}^{1}(T)$, and that $B_{h}^{\curl}$ is in $\mathcal{Z}_{h}^{1}(\Omega)$ if and only if $\left.{B_{h}^{\curl}}\right|_{T}$ is in $\mathcal{Z}_{h}^{1} (T)$ for all $T$, {we obtain the desired result}.
\end{proof}

\begin{proposition}\label{prop:WP_B_h_div}
    There exists a unique solution to the following variational formulation. Find $(B_{h}^{\div}, u_h) \in  \mathring{V}_{h}^{2}(\Omega) \times \left(V_{h}^{3}(\Omega) \cap L^2_0(\Omega)\right)$ such that, for all $(\tau_h, v_h) \in \mathring{V}_{h}^{2}(\Omega) \times V_{h}^{3}(\Omega)$
    \begin{align}
        \begin{split}\label{eq:VF_B_h_div}
            &{ \int_{\Omega} \left(\div B_{h}^{\div}\right) v_h} =0,
            \\
            &{ \int_{\Omega} B_{h}^{\div} \cdot \tau_h + \int_{\Omega} u_h \left(\div \tau_h\right)} = 2\pi\int_{\Sigma} \tau_h \cdot n_\Sigma.
        \end{split}
    \end{align}
    Furthermore, $B_{h}^{\div}$ is in $\mathring{\mathcal{K}}_{h}^{2}(\Omega)$.
\end{proposition}

\begin{proof}
    The well-posedness comes from the exact same arguments as in the continuous case, by replacing the Poincaré inequalities and Hodge decompositions by their discrete counterparts.

    We now prove that $B_{h}^{\div}$ is in $\mathring{\mathcal{K}}_{h}^{2}(\Omega)$. The fact that $B_{h}^{\div}$ is in $\mathring{\mathcal{Z}}_{h}^{2}(\Omega)$ comes directly from the first equation of (\ref{eq:VF_B_h_div}). Now, taking $\tau_h=\curl \rho_h$ in the second equation for $\rho_h \in \mathring{V}_h^1(\Omega)$, we get
    \begin{align*}
        { \int_{\Omega} B_{h}^{\div} \cdot \curl \rho_h} &= 2\pi \int_{\Sigma}\curl \rho_h \cdot n_{\Sigma}
        \\
        &= 2\pi \int_{\gamma'} \rho_h \cdot t'
        \\
        &= 0,
    \end{align*}
    so that $B_{h}^{\div}$ is in $\mathring{\mathcal{B}}_{h}^{2}(\Omega)^\perp$.
\end{proof}

We prove the two following approximation results.
\begin{proposition}\label{prop:approx_B}
    There exists a constant $C$ independent of $h$ such that
    \begin{align*}
        \left\|B_{h}^{\curl}-B\right\|_{L^2} &\leq C h^s \left\|B\right\|_{L^2},
        \\
        \left\|B_{h}^{\div}-B\right\|_{L^2} &\leq C h^s \left\|B\right\|_{L^2}.
    \end{align*}
\end{proposition}
\begin{proof}
    We begin by proving the convergence for $B_{h}^{\curl}$. From \cref{lemma:discrete_harmonic_fields}, we know there exists $\tilde{B}_h$ in $\mathcal{K}(\Omega)$ such that 
    $$
    \left\| \tilde{B}_h - B_{h}^{\curl} \right\|_{L^2} \leq \left\| \tilde{B}_h - \Pi_{h}^{1}\tilde{B}_h \right\|_{L^2}.
    $$
    From the continuity of the circulation along $\gamma$ with respect to the $H\curl$ norm, we have 
    \begin{align*}
        \left|\int_{\gamma} \left( \tilde{B}_h - B \right) \cdot t\right| &= \left|\int_{\gamma} \left(\tilde{B}_h - B_{h}^{\curl}\right) \cdot t \right|
        \\
        &\leq C \left\| \tilde{B}_{h} - B_{h}^{\curl} \right\|_{H\curl}
        \\
        &\leq C \left\| \tilde{B}_{h} - B_{h}^{\curl} \right\|_{L^2}
        \\
        &\leq C \left\| \tilde{B}_{h} - \Pi_{h}^{1} \tilde{B}_{h}\right\|_{L^2}.
    \end{align*}
    Since $B$ and $\tilde{B}_{h}$ are in $\mathcal{K}(\Omega)$, we get from \cite[Lemma 5]{valliVariationalInterpretationBiot2019} that
    \begin{align*}
        \left\| \tilde{B}_{h}-B \right\|_{L^2} &\leq C\left( \left\| \div \tilde{B}_{h}- \div B \right\|_{L^2} + \left\| \curl \tilde{B}_{h}- \curl B \right\|_{L^2} + \left|\int_{\gamma} \left( \tilde{B}_h - B \right) \cdot t\right|\right)
        \\
        &\leq C \left\| \tilde{B}_{h} - \Pi_{h}^{1} \tilde{B}_{h}\right\|_{L^2}.
    \end{align*}
    We then get by triangle inequality 
    \begin{align*}
        \left\| B_{h}^{\curl} - B \right\|_{L^2} &\leq \left\| B_{h}^{\curl} - \tilde{B}_{h} \right\|_{L^2} + \left\| \tilde{B}_{h} - B \right\|_{L^2}
        \\
        &\leq C \left\| \tilde{B}_{h} - \Pi_{h}^{1} \tilde{B}_{h}\right\|_{L^2}.
    \end{align*}
    Finally, using \cref{prop:approx_quasi_interpolator} and the continuous injection of $H(\curl, \Omega) \cap H_0(\div, \Omega)$ into $H^s(\Omega)$ for some $1/2 < s \leq 1$, we get 
    \begin{align*}
        \left\| B_{h}^{\curl} - B \right\|_{L^2} &\leq C h^s \left\| \tilde{B}_{h} \right\|_{H^s}
        \\
        &\leq C h^s \left\| \tilde{B}_{h} \right\|_{L^2}
        \\
        &\leq C h^s \left\| B \right\|_{L^2}.
    \end{align*}
    We now prove the convergence of $B_{h}^{\div}$ to $B$. Using \cref{eq:B_div,eq:VF_B_h_div}, we have 
    \begin{align*}
        { \int_{\Omega} B \cdot \tau + \int_{\Omega} u \left(\div \tau_h\right)} &= \int_{\Sigma} \tau_h \cdot n_{\Sigma},
        \\
        { \int_{\Omega} B_{h}^{\div} \cdot \tau + \int_{\Omega} u_h \left(\div \tau_h\right)} &= \int_{\Sigma} \tau_h \cdot n_{\Sigma},
    \end{align*}
    for all $\tau_h$ in $\mathring{V}_{h}^{2}(\Omega)$. As a consequence, for all $\tau_h$ in $\mathring{\mathcal{Z}}_{h}^{2}(\Omega)$, we get
    $$
    { \int_{\Omega} \left(B_{h}^{\div}-B\right) \cdot \tau_h} = 0.
    $$
    Since both $B_{h}^{\div}$ and $\Pi_{h}^{2}B$ are in $\mathring{\mathcal{Z}}_{h}^{2}(\Omega)$, we have
    $$
    { \int_{\Omega} \left(B_{h}^{\div}-B\right) \cdot \left(B_{h}^{\div} - B\right) = \int_{\Omega} \left(B_{h}^{\div} - B\right) \cdot \left(\Pi_{h}^{2}B - B\right)},
    $$
    so that
    $$
    \left\| B_{h}^{\div} - B \right\|_{L^2} \leq \left\| \Pi_{h}^{2}B - B \right\|_{L^2}.
    $$
    Using again \cref{prop:approx_quasi_interpolator} and the fact that $\Omega$ is $s$-regular, we get the desired result
    $$
    \left\| B_{h}^{\div} - B \right\|_{L^2} \leq C h^s \|B\|_{L^2}.
    $$
\end{proof}

We now study the well-posedness and the convergence of the vector potentials of $B$. A first remark is that the vector potential orthogonal to $B$ is well studied in the literature via the Hodge Laplacian (see for example \cite{arnoldFiniteElementExterior2006}). As such, we study the approximation of such vector potentials which we denote $A^1_{h}$, and recover $A_h^2$, the one which is circulation free, by the discrete counterpart of \cref{eq_relation_pot_vec_continue}
\begin{align}
    \label{eq_relation_pot_vec_discrete}
    A^2_{h}=A^1_{h}-\frac{1}{2 \pi}\left(\int_{\gamma'}A^1_{h} \cdot t'\right) B_{h}^{\curl}.
\end{align}

\begin{proposition}\label{prop:WP_A_h}
    There exists a unique solution to the following variational formulation. Find $(A^1_{h}, u_h) \in V_{h}^{1}(\Omega) \times V_{h}^{2}(\Omega)$ such that, for all $(\tau_h, v_h) \in V_{h}^{1}(\Omega) \times V_{h}^{2}(\Omega)$,
    
    \begin{align}
        \begin{split}\label{eq:WP_A_h}
        &{ \int_{\Omega} A^1_h \cdot \tau_h = \int_{\Omega} \curl \tau_h \cdot u_h},
        \\
        &{ \int_{\Omega} \curl A^1_h \cdot v_h + \int_{\Omega} \left(\div u_h\right) \left(\div v_h\right) = \int_{\Omega} B_{h}^{\div} \cdot v_h}. 
        \end{split}
    \end{align}
    Furthermore, we have $\curl A_h^1 = B_h^\div$.
\end{proposition}

\begin{proof}
    The well posedness is proven in the exact same way as in the continuous case \cref{pr:wp_vector_potential}, by replacing the Hodge decomposition and Poincaré inequality by their discrete counterparts. The fact that $\curl A_h^1 = B_h^\div$ is also proven in a similar way as in \cref{pr:wp_vector_potential}. Using $\mathcal{K}_h^2(\Omega) = 0$, and the discrete Hodge decomposition (see \cref{sub-sec:Discrete_De_Rham}), we have
    $$
    V_h^2(\Omega) = \mathcal{Z}_h^2(\Omega)^\perp \oplus \mathcal{B}_h^2(\Omega),
    $$
    so that $u_h = u_h^\nabla + u_h^\curl$ with ${ \int_{\Omega} u_h^\nabla v} = 0$ for all $v_h$ in $V_h^2(\Omega)$ with $\div v_h=0$, and $\div u_h^\curl = 0$. Therefore, taking $v_h=u_h^\curl$ in the second equation of \cref{eq:WP_A_h}, we get 
    $$
    { \int_{\Omega} \left(\div u_h\right) \left(\div u_h^\nabla\right)} = 0,
    $$
    so that, since $\div u_h^\curl = 0$, we obtain $\div u_h^\nabla = \div u_h = 0$. From this, we get ${ \int_{\Omega} A_h^1 \cdot v_h = \int_{\Omega} B_h^\div \cdot v_h}$ for all $v_h \in V_h^2(\Omega)$, and $B_h^\div \in V_h^2(\Omega)$ implies the desired equality.
\end{proof}

\begin{proposition}\label{prop:approx_A}
    There exists a constant $C$ independent of $h$ such that
    $$
    \left\|A_h^1 - A^1\right\|_{L^2} \leq C h^s \left\|B\right\|_{L^2}.
    $$
\end{proposition}

\begin{proof}
    Define $\left(\tilde{A}_h^1, \tilde{u}_h\right)$ as the solution to the following variational problem. Find $\left(\tilde{A}_h^1, \tilde{u}_h\right) \in V_{h}^{1}(\Omega) \times V_{h}^{2}(\Omega)$ such that, for all $(\tau_h, v_h) \in V_{h}^{1}(\Omega) \times V_{h}^{2}(\Omega)$,
    \begin{align*}
        &{ \int_{\Omega} \tilde{A}_h^1 \cdot \tau_h  = \int_{\Omega} \curl \tau_h \cdot \tilde{u}_h},
        \\
        &{ \int_{\Omega} \curl \tilde{A}_h^1 \cdot v_h + \int_{\Omega} \left(\div \tilde{u}_h\right) \left(\div v_h\right) = \int_{\Omega} B \cdot v_h}. 
    \end{align*}
    The well-posedness of this equation is obtained in the same way as for $A_h^1$. Furthermore, by continuity of the resolvant, and independence of the discrete Poincaré inequality constants in $h$, we get 
    $$
    \left\| \tilde{A}_h^1 - A_h^1 \right\|_{L^2} \leq C \left\| B - B_{h}^{\div} \right\|_{L^2}.
    $$
    Furthermore, we obtain from of \cite[Theorem 7.9]{arnoldFiniteElementExterior2006} and $\| B \|_{H^s} \leq C \|B\|_{L^2}$ that 
    $$
    \left\| \tilde{A}^1_h - A^1 \right\|_{L^2} \leq C h^s \| B \|_{L^2}.
    $$
    We then get the desired result from a triangle inequality and \cref{prop:approx_B}.
\end{proof}

As we have seen in \cref{subsec:VP_BG}, there are two ways of computing the harmonic helicity of $\Omega$. The first one is done by computing circulations of $A^1$, and the second by taking the $L^2$ product of $B$ and $A^2$. As we shall see in the following results, we can also recover the numerical harmonic helicity in two similar ways. To do this, we first need to establish convergence of the circulation of $A_h^1$, which is given by the following lemma. 

\begin{lemma}\label{lemma:convergence_circulation_A^1}
    We have 
    $$
    \left| \int_{\gamma'} (A_h^1 - A^1) \cdot t' \right| \leq Ch^s\|B\|_{L^2}.
    $$
\end{lemma}

\begin{proof}
    From \cref{prop:WP_A_h}, we have 
    \begin{align*}
        \left\|A_h^1 - A^1\right\|_{H_\curl}^2 &= \left\|A_h^1 - A^1\right\|_{L^2}^2 + \left\|\curl A_h^1 - \curl A^1\right\|_{L^2}^2
        \\
        &= \left\|A_h^1 - A^1\right\|_{L^2}^2 + \left\|B_h^\div - B\right\|_{L^2}^2,
    \end{align*}
    so that by \cref{prop:approx_B,prop:approx_A}, we get 
    $$
    \left\|A_h^1 - A^1\right\|_{H_\curl} \leq C h^s \|B\|_{L^2}.
    $$
    Now, since $\curl A_h^1 = B_h^\div$, which is in $H_0(\div, \Omega)$, the circulation of $A_h^1$ along $\gamma'$ is well-defined by \cref{eq:circ_gamma_bound}, and this circulation is bounded by the $H(\curl,\Omega)$ norm. We therefore have 
    $$
    \left| \int_{\gamma'}(A_h^1 - A^1) \cdot t \right| \leq Ch^s\|B\|_{L^2}.
    $$
\end{proof}

\begin{corollary}\label{cor:convergence_A^2}
    Defining $A_h^2$ as in \cref{eq_relation_pot_vec_discrete}, we get 
    $$
    \left\|A_h^2 - A^2\right\|_{L^2} \leq Ch^s\|B\|_{L^2}.
    $$
\end{corollary}

\begin{proof}
    We get this result by writing 
    $$
    A_h^2 - A^2 = A_h^1 - A^1 - \frac{1}{2\pi}\left(\int_{\gamma'}\left(A_h^1 - A^1\right)\cdot t \right) B_h^\curl + \frac{1}{2\pi}\left(\int_{\gamma'}A^1 \cdot t \right)\left(B-B_h^\curl\right).
    $$
    {We then obtain the desired result using \cref{prop:approx_B,prop:approx_A,lemma:convergence_circulation_A^1}.}
\end{proof}

\begin{theorem}[Convergence of harmonic helicity]\label{th:hh_convergence}
    We have 
    $$
    \left| \mathcal{H}(\Omega) - \mathrm{H}(B_{h}^\div) \right| \leq Ch^s\|B\|_{L^2}^2,
    $$
    where $\mathrm{H}(B_h^\div)$ can be either computed as 
    $$
    \mathrm{H}(B_h^\div) = -\left(\int_{\gamma'} A_h^1 \cdot t'\right) \left(\int_{\gamma} A_h^1 \cdot t\right),
    $$
    or as
    $$
    \mathrm{H}(B_h^\div) = \int_{\Omega} B_h^\div \cdot A_h^2.
    $$
\end{theorem}

\begin{proof}
    We use the formulation of the helicity through $A_h^2$, the other one being a simple consequence of the Bevir--Gray formula. To obtain the desired estimate, we write 
    $$
    \int_{\Omega} B_{h}^{\div} \cdot A_h^2 - \mathcal{H}(\Omega) = \int_{\Omega} (B_{h}^{\div}-B) \cdot A_h^2 + \int_{\Omega} B \cdot (A_h^2 - A^2).
    $$
    The conclusion then follows from \cref{prop:approx_B,cor:convergence_A^2}.
\end{proof}

\begin{remark}
    \label{rmk:pb_cvg_shape_grad}
    Although we have proven convergence of the harmonic helicity, convergence of the shape gradient is a harder task. Indeed, such a result is related to traces of the numerical solutions. From the estimates obtained in this section, we only have convergence of $B_h^\curl$ and $A_h^{2}$ in $H(\curl, \Omega)$, and of $B_h^\div$ in $H(\div, \Omega)$.
    Hence, the tangential traces of both $A_h^2$ and $B_h^\curl$ converges in $H^{-1/2}(\partial \Omega)$. Nevertheless, as we do not expect the numerical solutions to have any higher regularity than the ones prescribed by the discretized functional spaces\footnote{In the continuous setting, the reasoning works because $A^2$ and $B$ are in $H(\curl, \Omega) \cap H_0(\div, \Omega)$. Nevertheless, the divergence of $A_h^{2}$ and $B_h^\curl$ is not in $L^2(\Omega)$ for Nedelec elements}, this is not sufficient to ensure convergence toward the shape gradient \cref{eq:shape_derivative_formula}. This means that we cannot simply use Sobolev estimates inside $\Omega$, and need to work on the boundary directly in order to obtain, for example, convergence in $L^2(\partial \Omega)^3$ (or at least on the tangential part). Such a convergence result is therefore non-trivial with usual techniques, and this problem remains open.
\end{remark}

\section{Numerical implementation and results}
\label{sec:num}
\subsection{Specificity of simulations for stellarators}
\label{subsec-specifi-stella}
For both historical and practical considerations, it is frequently advantageous for the surfaces under examination to exhibit specific symmetries. In particular, a majority of stellarators are invariant under discrete rotations along the Oz axis, with an angle of $2\pi/N_p$, where $N_p$ takes values of 3 (as in the case of NCSX~\cite{zarnstorffPhysicsCompactAdvanced2001}) or 5 (as observed in W7X~\cite{warmerW7XHELIASFusion2017}). Additionally, stellarators commonly exhibit invariance under the continuous symmetry known as stellarator symmetry, as discussed in detail in \cite[Section 12.3]{imbert-gerardIntroductionSymmetriesStellarators2019} and \cite{dewarStellaratorSymmetry1998}. In practice, these surfaces are represented by a set of coefficients $(R_{m,n}),(Z_{m,n})$ for $m\in \N$ and $n\in \Z$ which define the functions
\begin{align}
    R(u,v)&=\sum_{m \geq 0} \sum_{n\in \Z} R_{m,n}\cos(2 \pi (m u+n v)) ,\label{eq:R_in_fourier}\\
    Z(u,v)&=\sum_{m \geq 0} \sum_{n\in \Z} Z_{m,n}\sin(2 \pi (m u+n v)) .
    \label{eq:Z_in_fourier}
\end{align}
Note the absence of $\sin$ terms for $R$ and $\cos$ terms for $Z$ due to the imposition of stellarator symmetry. The surface is then parametrized in cylindrical coordinates by $(R(u,v),\frac{2\pi v}{N_p},Z(u,v))$ where $N_p$ stands for the discrete symmetry imposed.

For the numerical simulation, we truncate the number of Fourier components in \eqref{eq:R_in_fourier} and \eqref{eq:Z_in_fourier}.

\begin{figure}
\center
\includegraphics[width=0.85\textwidth]{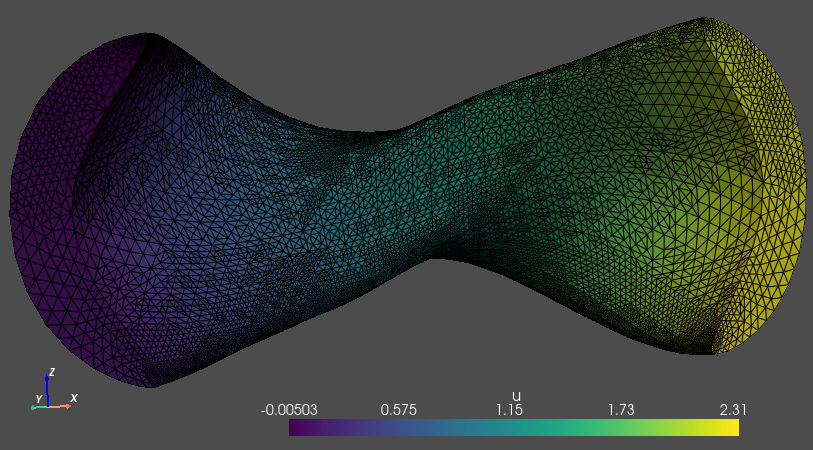}
\includegraphics[width=0.85\textwidth]{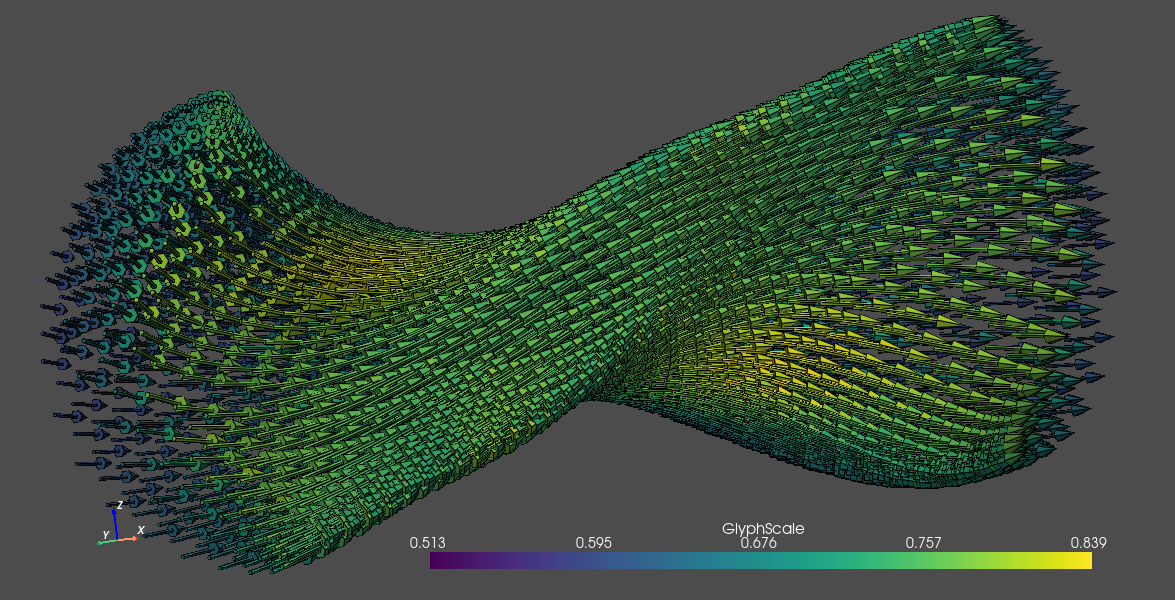}
\includegraphics[width=0.85\textwidth]{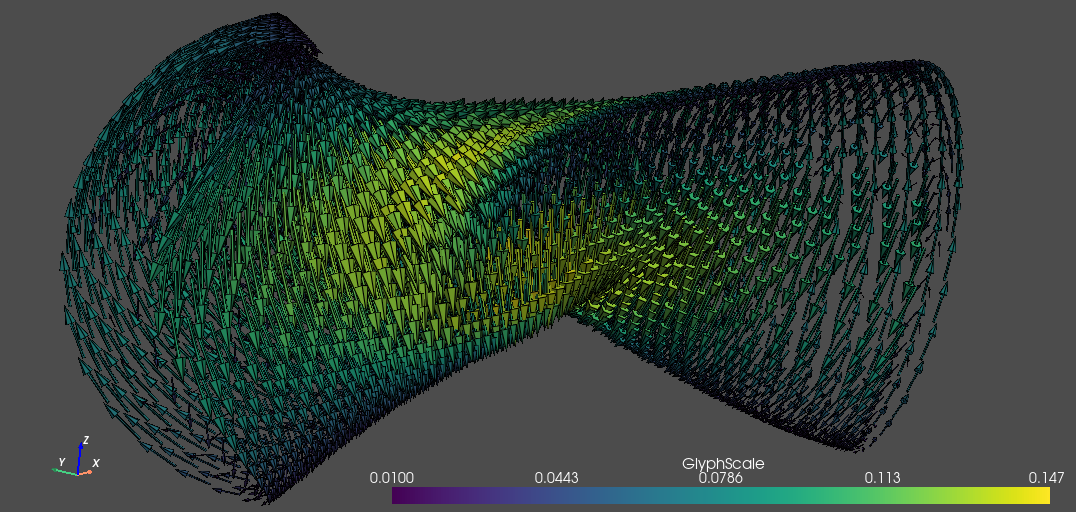}
\caption{One section (i.e. one third) of NCSX plasma. In the upper plot, we represent the function $u_h$ of \cref{eq:VF_B_h_curl}. Its gradient $B_h^\curl$ is shown in the middle figure. Note that the boundary conditions on the left and right cuts are a jump $2\pi/3$, whereas we have Neumann boundary condition on the plasma surface. The bottom figure is a reprentation of $A^2_h$.}
\label{fig:multiple_repr}
\end{figure}
\subsection{Implementation}

Using the symmetries of the system, we only work with one section of the stellarator, characterized by the coefficients $(R_{m,n},Z_{m,n})$. Our first step is to use Gmsh \cite{geuzaineGmsh3DFinite2009} to mesh the interior of a polyhedral approximation of the surface.

Then we use the finite element library FEniCSx (\cite{scroggsConstructionArbitraryOrder2022,scroggsBasixRuntimeFinite2022,alnaesUnifiedFormLanguage2014}) to assemble the finite element problems \cref{eq:VF_B_h_div,eq:VF_B_h_curl,eq:WP_A_h} to obtain $B_{h}^{\div},B_{h}^{\curl}$ and $A_h^1$. We implemented both first or second order FEEC elements defined in \cref{sub-sec:recall-FEM} with adequate periodic conditions at the cuts of our section. To solve the linear system required to get $B_h^{\curl}$, that is a Poisson equation, we use the solver MUMPS: MUltifrontal Massively Parallel sparse direct Solver \cite{amestoyFullyAsynchronousMultifrontal2001,amestoyHybridSchedulingParallel2006}. $B_{h}^{\div}$ and $A_h$ are more expensive and complex to solve as both are defined in mixed formulations. Both are solved using the iterative Krylov method GMRES \cite{saadGMRESGeneralizedMinimal1986} applied after a careful preconditioning. More precisely, we use the block diagonal preconditioner proposed by Arnold et al. in \cite[Section 10.2]{arnoldFiniteElementExterior2006} for Hodge-Laplacian problems. The preconditioning problems are solved using MUMPS. We believe that using an Auxiliary-space Maxwell (AMS) Solver \cite{hiptmairNodalAuxiliarySpace2007} instead\footnote{We had issues porting AMS hypre~\cite{falgoutHypreLibraryHigh2002} solvers to the background sparse linear library Petsc~\cite{balay1998petsc}} of a direct solver for the preconditioner would improve the efficiency and scalability of our method when using more than a few million degrees of freedom.

\subsection{Numerical tests}

\begin{figure}
    \includegraphics[width=\textwidth]{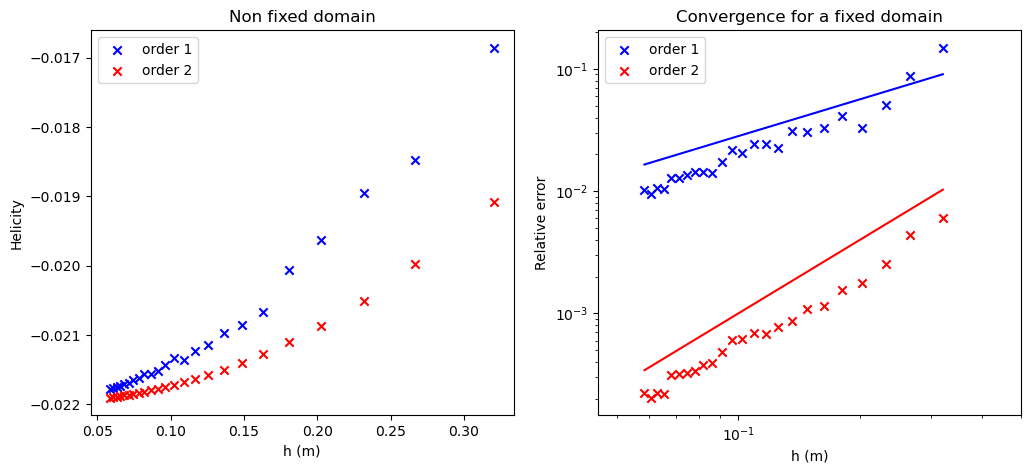}
    \caption{On the left, we provide better and better approximations of the surface to the mesher. We then plot the characteristic size of an element versus the harmonic helicity. On the right, we mesh with different characteristic size $h$ the same polyhedral approximation of the continuous surface { and compare it to a reference solution obtained with second order elements and $h=0.025$. The two lines are merely indicative and their slopes correspond to theoretical first order (blue) and second order (red) convergence rates.}
    }
    \label{fig:helicity_vs_h}
\end{figure}

Once we have computed $B_{h}^{\div},B_{h}^{\curl}$ and $A_h$ that are represented in \cref{fig:multiple_repr}, we compute the harmonic helicity and its shape gradient. {As was stated in \cref{cor:axisymmetric}, we know that axisymmetric tori have zero harmonic helicity. Numerically, we have found in this case a value of $10^{-9}$. This was computed with a major and minor radius of 1 and 0.1 respectively, second order elements, and $h$ equal to 0.025.} For more interesting shapes, as NCSX's plasma, we do not have a reference. As mentioned in \cref{subsec-specifi-stella}, we have a continuous description of the shape (i.e. not a polyhedral one), thus we perform two tests.
\begin{enumerate}
    \item On the left side of \cref{fig:helicity_vs_h}, we perform a better and better approximation of the continuous shape. To this aim we provide finer and finer grid description of the surface to the mesher. Hence, variation of the obtained magnetic helicity are due both to variations of the domain and to the finite element approximation error.
    \item On the right side of \cref{fig:helicity_vs_h}, we fix at the beginning a polyhedral shape that is an approximation of the continuous surface described by the set of coefficients $(R_{m,n},Z_{m,n})$ and use finer and finer meshes of this fixed polyhedral shape as $h$ goes to $0$. Hence, variations of the magnetic helicity are only due to the numerical approximation of the fields. This is the situation described by the theoretical analysis in \cref{sec:approx_fecc}. { The reference solution was obtained on a high-performance computing cluster using a mesh size of $h=0.025$ with second-order finite elements. This resulted in linear systems with approximately 5 million degrees of freedom, which were solved in a matter of minutes.}
\end{enumerate}
{ Except for the reference solution, all simulations can be conveniently performed on a laptop with 32 GB of RAM.} 

\begin{figure}
    \center
    \includegraphics[width=0.85\textwidth]{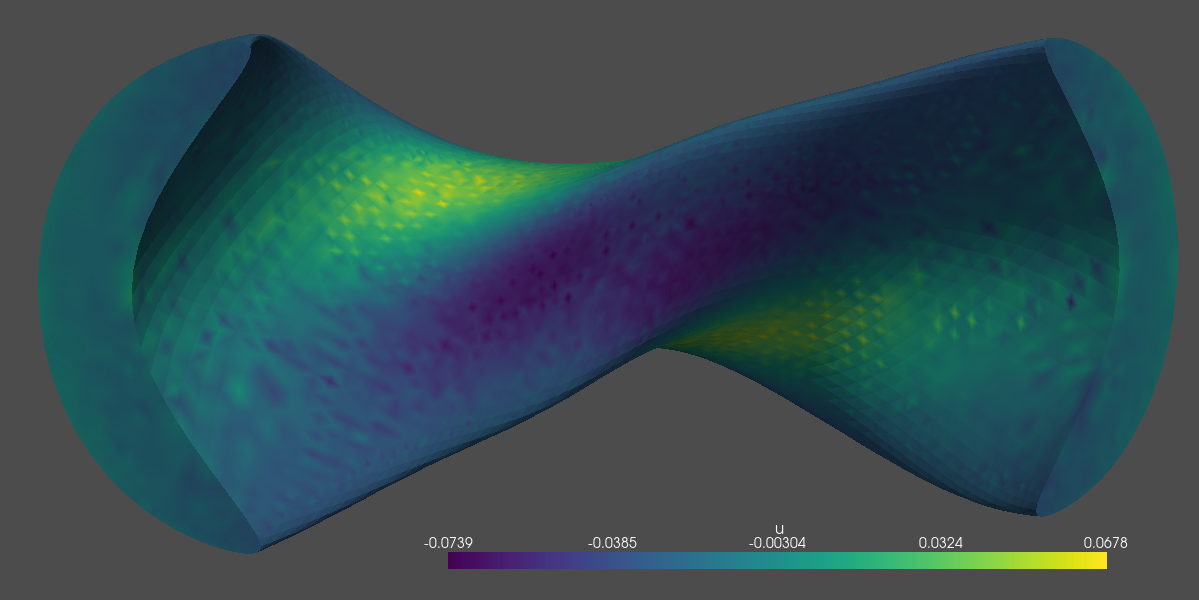}
    \caption{One section of NCSX plasma. The plot shows the shape gradient $2A^2_h\cdot B_h^\curl$ on the boundary. NCSX plasma configuration has a negative harmonic helicity. Improving our criterion implies following the opposite of the shape gradient.}
    \label{fig:shape_gradient}
\end{figure}

The shape gradient is shown in \cref{fig:shape_gradient}. As mentioned in \cref{rmk:pb_cvg_shape_grad}, we have not proven that this numerical approximation converges toward the correct shape gradient. However, we performed numerical tests using finite differences on the Fourier coefficients of the surface, which were consistent with the numerical shape gradient computed with $A^2_h$ and $B_h^\curl$.

\subsection{Two optimization programs}

{ 
As the harmonic helicity of NCSX plasma is a negative quantity\footnote{Using a planar symmetry, we could also choose to take a positive value, see \cref{prop:vol}.}, we are interested in minimizing it.  
Indeed, we recall that in a stellarator, one generates a magnetic field that is in first approximation a harmonic field in the plasma domain $\Omega$. As the stability of the plasma is related to the twisting of the magnetic field, we are trying to increase the absolute value of the harmonic helicity to search for interesting and hopefully more stable plasma configurations.

As proved in \cref{prop:vol}, the harmonic helicity scales as $\mathcal{H}(\lambda \Omega) =\lambda \mathcal{H}(\Omega)$ for $\lambda >0$. For this reason we impose an upper bound either on the volume of $\Omega$ or on the area of $\partial \Omega$. Under one of these constraints, it is unclear whether an optimal shape exists. Hence, we add the condition that $\partial \Omega$ satisfies a uniform ball condition. This imposes a strong regularity constraint on the surface and ensures existence of an optimal shape among these regular shapes, we refer for example to \cite{privatExistenceSurfacesOptimizing2022,gernerExistenceOptimalDomains2023}. Numerically, we set the uniform ball condition by imposing a lower limit on the minimal curvature radius everywhere on the surface $\partial \Omega$. These three constraints, namely upper bound on the volume $\Omega$ or on the area of $\partial \Omega$ and lower bound on the minimal curvature radii on $\partial \Omega$, along with their corresponding shape gradients, can be efficiently computed using the smooth parametrization of the surface outlined in \cref{eq:R_in_fourier,eq:Z_in_fourier} and performing quadrature across either $\partial \Omega$ or $\Omega$.}

Subsequently, we introduce smooth non-linear costs that blow up to guarantee the fulfillment of the specified constraints. Once all the costs are assembled, we apply the Broyden--Fletcher--Goldfarb--Shanno (BFGS) minimization algorithm from the scipy library \cite{virtanenSciPyFundamentalAlgorithms2020} on the Fourier components of the description of the surface.

\begin{figure}
    \includegraphics[width=\textwidth]{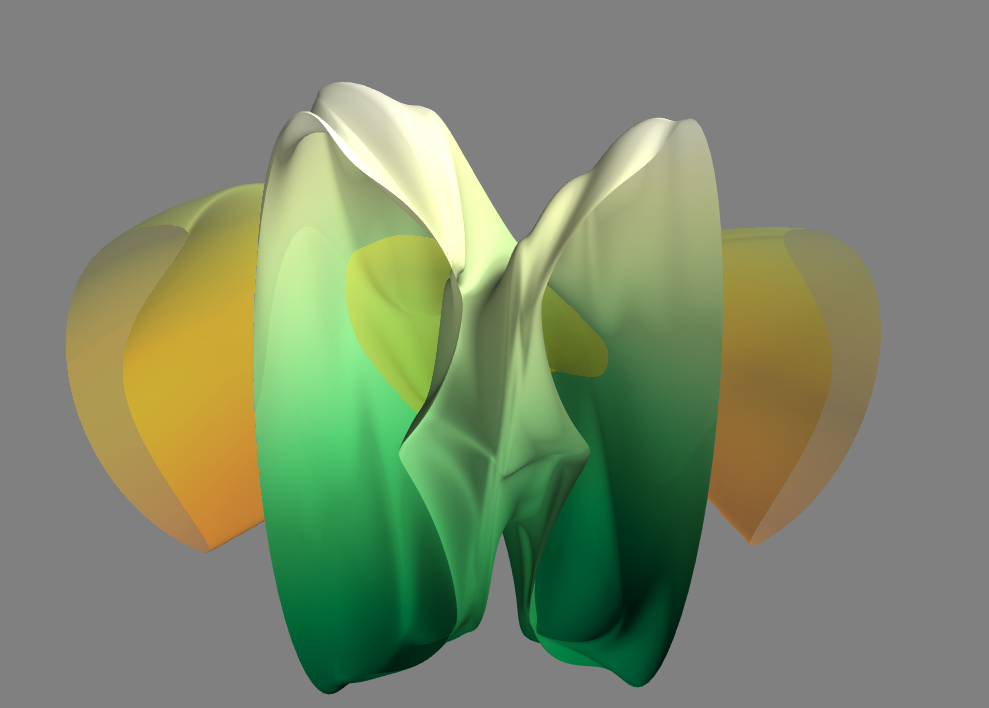}
    \includegraphics[width=\textwidth]{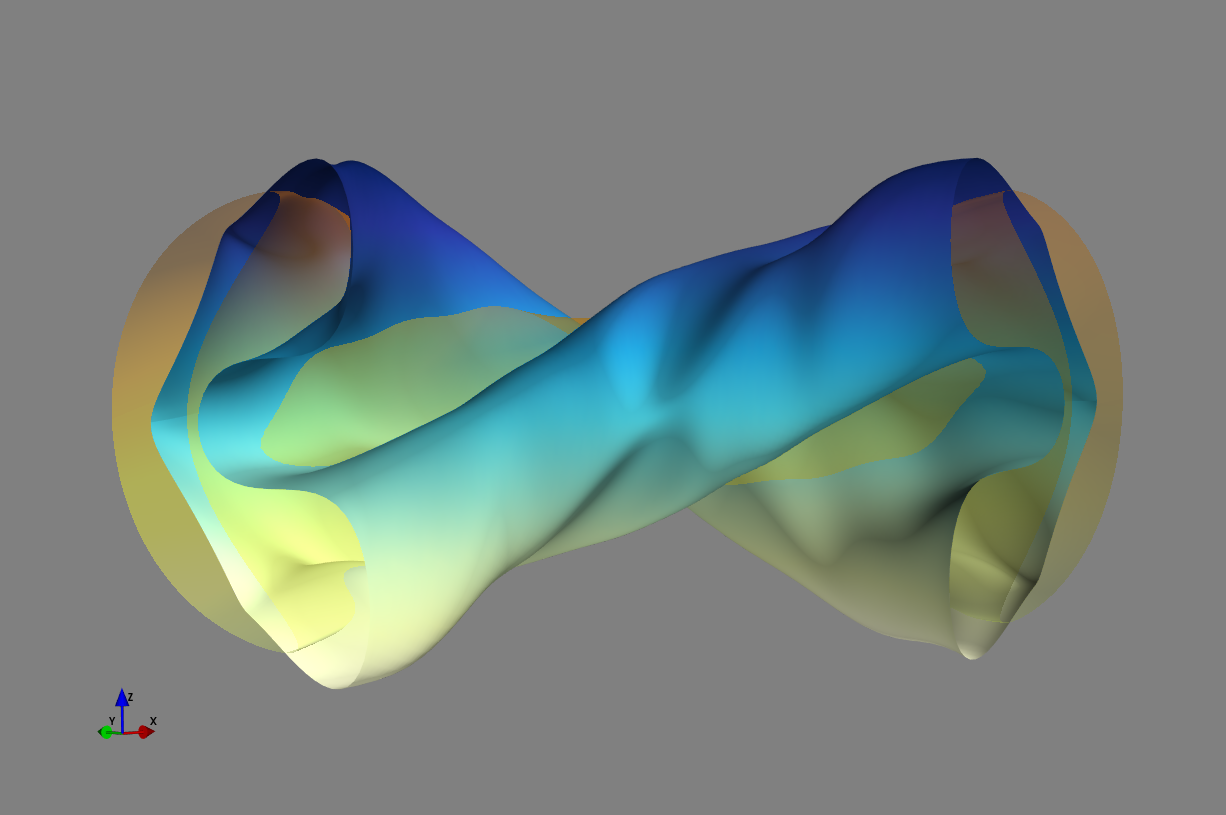}
    \caption{The results of the two optimization problems, BPC (shown in green above) and BVC (shown in blue below). The associated costs are reported in \cref{tab:res}. Additionally, the reference NCSX plasma is plotted in transparent orange.}
    \label{fig:optimized_shapes}
\end{figure}

\paragraph*{Bounded Perimeter and Curvature (BPC)}
We initially focus on optimization with a bounded perimeter, employing a soft threshold set at 25~$m^2$ to ensure comparability with the NCSX plasma shape. Our goal also involves achieving a minimal inverse curvature radius of 50 $m^{-1}$, surpassing the regularity of the NCSX reference. The results of this simulation are depicted in \cref{fig:optimized_shapes} and summarized in \cref{tab:res}. Notably, we observe a remarkable 30-fold improvement in harmonic helicity. While the shape appears to collapse toward the Oz axis, this effect is restrained by the imposed curvature constraints. This behavior may be attributed to our normalization choice for the harmonic field; as proximity to the Oz axis increases, the magnitude rises, given our fixed toroidal circulation.

\paragraph*{Bounded Volume and Curvature (BVC)}
We also perform the optimization using a constraint on the volume instead of the perimeter. We apply an upper bound at 3~$m^3$ along with the same curvature constraints. This time, we achieve more than a 3-fold improvement, and the shape does not seem to collapse toward the Oz axis. The shape is more intricate than the initial plasma shape of NCSX.

\begin{table}
    \centering
    \begin{tabular}{|c|c|c|c|c|}
        \hline
        Plasma shape & Perimeter $(m^2)$ & Volume $(m^3)$ & min curvature radius $(cm)$& Harmonic helicity \\
        \hline
        NCSX (ref) & 24.52 & 2.96 & 0.83& -0.0220 \\
        BPC & 25.02 & 5.45 & 1.89& -0.6427 \\
        BVC & 33.5 & 2.99& 2.00 & -0.0726 \\
        \hline
    \end{tabular}
    \caption{harmonic helicity and geometric properties of the optimized shapes in \cref{fig:optimized_shapes}. { As we aim to compare the optimized shapes with the NCSX reference, for BPC we penalized the area of the surface $\partial \Omega$ above $25\, m^2$, whereas for BVC we penalized the volume of $\Omega$ above $3\, m^3$. In both cases, we targeted a minimal curvature radius of $2\,cm$; thus a more regular shape than the reference NCSX plasma shape. }}
    \label{tab:res}
\end{table}

\section{Conclusion and perspectives}
We introduced a new shape functional for toroidal domains capturing the linkage of the corresponding harmonic field. Then, using carefull differential form pullbacks, we have been able to compute its shape derivative. We also showed that this functional can be efficiently computed numerically. We illustated this using finite elements exterior calculus and applied it to a stellarator device for magnetic confinement in nuclear fusion.

It seems to us that the pure optimization of harmonic helicity provides degenerate forms, which appear to be less usable for applications in magnetic confinement. However, given the very significant potential gain (factors of 3 and 30 in our two simulations), studying multi-objective optimization with more conventional costs could yield interesting results.

We also believe that the following perspectives are of interest: 

\begin{itemize}
    \item The convergence of the numerical shape gradient is currently incomplete, as mentioned in \cref{rmk:pb_cvg_shape_grad}. The tools to solve this problem, could be of importance outside the scope of this article.

    \item Our numerical results strongly indicate that removing the reach constraints might result in the non-existence of an optimal shape, However this is only a conjecture.
\end{itemize}

\paragraph{Data availability statement} The code created for this article is openly accessible through the GitLab repository \url{https://plmlab.math.cnrs.fr/rrobin/helicity}.

\paragraph*{{Acknowledgements}}The authors express sincere gratitude to Mario Sigalotti, Yannick Privat, and Alexandre Ern for their valuable advice and thoughtful discussions, which greatly enriched this work. Additionally, they thank all the members of the StellaCage collaborations\footnote{\url{https://www.ljll.math.upmc.fr/~sigalotti/cage/stellacage.html}} for their insights on the physical aspects of this work.
This research received support from Inria AEX StellaCage.

\appendix
\section{Translation from differential forms and Hodge decomposition}
\label{sec:appendix_hodge}
It is common in the literature to find the types of problems we are studying in the language of differential forms. Although this approach deals with more abstract mathematical objects, it allows the use of a unique framework. This approach can notably be found in \cite{arnoldFiniteElementExterior2006}, and is relatively common in the finite elements exterior calculus literature.

Since the aforementioned book is quite often referenced in this article, we chose to write a small appendix to explain how to translate problems from the language of differential forms to the one of functions and vector fields. This is done using usual identifications from Riemannian geometry and Hodge theory, that is, the musical isomorphisms and Hodge star operator. In the case of a three-dimensional manifold, these identifications work perfectly to equate 0 and 3 forms with functions, 1 and 2 forms with vector fields, and to translate the exterior derivatives and coderivatives with the usual differential operators of electromagnetism. These identifications are then summed up in the following commutative diagram, and table.
\\
\begin{equation}\label{diagram}
    \begin{tikzcd}
        {H\Lambda^0(\Omega)} & {H\Lambda^1(\Omega)} & {H\Lambda^2(\Omega)} & {H\Lambda^3(\Omega)} \\
        {H^1(\Omega)} & {H(\mathrm{curl}, \Omega)} & {H(\mathrm{div}, \Omega)} & {L^2(\Omega)}
        \arrow["{\mathrm{d}}", from=1-1, to=1-2]
        \arrow["{\mathrm{d}}", from=1-2, to=1-3]
        \arrow["{\mathrm{d}}", from=1-3, to=1-4]
        \arrow["{\mathrm{id}}"', from=1-1, to=2-1]
        \arrow["{\#}"', from=1-2, to=2-2]
        \arrow["{\# *}"', from=1-3, to=2-3]
        \arrow["{*}"', from=1-4, to=2-4]
        \arrow["\nabla", from=2-1, to=2-2]
        \arrow["{\mathrm{curl}}", from=2-2, to=2-3]
        \arrow["{\mathrm{div}}", from=2-3, to=2-4]
    \end{tikzcd}
\end{equation}

\begin{figure}[H]\label{fig:table_correspondences}
    \center
    \begin{tabular}{|c||c|c|c|c|c|}
        \hline
        & $H\Lambda^k(\Omega)$ & $H^*\Lambda^k(\Omega)$ & $d$ & $\delta$ & $\mathrm{Tr}(\omega)$ \\
        \hline
        \hline
        $k=0$ & $H^1(\Omega)$ & $L^2(\Omega)$ & $\nabla$ & 0 & $\omega_{|\partial\Omega}$ \\
        \hline 
        $k=1$ & $H(\curl, \Omega)$ & $H(\div, \Omega)$ & $\curl$ & $-\div$ & $\omega \times n$ \\
        \hline 
        $k=2$ & $H(\div, \Omega)$ & $H(\curl, \Omega)$ & $\div$ & $\curl$ & $\omega \cdot n$ \\
        \hline 
        $k=3$ & $L^2(\Omega)$ & $H^1(\Omega)$ & $0$ & $-\nabla$ & -- \\
        \hline
    \end{tabular}
    \caption{Table of correspondence between differential forms language and vector calculus language}
\end{figure}

Here, $H\Lambda^k(\Omega)$ denotes the set of square integrable $k$-forms whose exterior derivative is square integrable, $H^*\Lambda^k(\Omega)$ the set of square integrable $k$-forms whose exterior coderivative is square integrable, $d$ the exterior derivative, $\delta$ the exterior coderivative, $*$ the Hodge star operator, $\#$ the musical isomorphism taking 1-forms to vector fields, and $\mathrm{Tr}(\omega)$ the trace of a differential form $\omega$ defined by the pullback of $\omega$ onto $\partial \Omega$ by the inclusion map $i : \partial \Omega \rightarrow \bar{\Omega}$.

Once these identifications are given, a tool which is used quite often throughout the paper is the Hodge decomposition. It is given by the following proposition.
\begin{proposition}\label{pr:hodge_decomposition}
    Let $\Omega$ be a Lipschitz toroidal domain as defined in \cref{sec:state_of_the_art}. We have the following $L^2$ orthogonal decompositions
    \begin{align}
        L^2(\Omega)^3&=\curl H_0(\curl, \Omega)\oplus^\perp\mathcal{K}(\Omega)\oplus^\perp\nabla H^1(\Omega), \label{eq:hodge_decomposition_harm}\\
        L^2(\Omega)^3&=\curl H(\curl, \Omega)\oplus^\perp\nabla H^1_{0}(\Omega). \label{eq:hodge_decomposition_no_harm}
    \end{align}
\end{proposition}
\begin{proof}
    This is a simple consequence of the Hodge decomposition given in \cite[Section 2]{arnoldFiniteElementExterior2006}, in the cases $k=1$ and $k=2$ respectively. The reason no harmonic term appears in \cref{eq:hodge_decomposition_no_harm} is that the second De Rham cohomology space vanishes in $\Omega$, and thus, that the set of harmonic two forms is trivial.
\end{proof}

\begin{proposition}\label{pr:continuous_poincare}
    There exists $C$ such that for all $V_1 \in H\left(\curl^0, \Omega\right)^\perp \cap H(\curl, \Omega)$, $V_2 \in H\left(\div^0, \Omega\right)^\perp \cap H(\div, \Omega)$
    \begin{align*}
        \left\|V_1\right\| &\leq C\left\| \curl V_1 \right\|,
        \\
        \left\|V_2\right\| &\leq C\left\| \div V_2 \right\|.
    \end{align*}
\end{proposition}
\begin{proof}
    This is given by \cite[Eq. (2.17)]{arnoldFiniteElementExterior2006} in the cases $k=1$ and $k=2$.
\end{proof}

\begin{proposition}\label{pr:orthogonality_relations}
    We have the following orthogonality relations in $L^2(\Omega)^3$
    \begin{align}
        \nabla H^1(\Omega)^\perp &= H_0\left(\div^0,\Omega\right),\label{eq:orth_gradients}
        \\
        \curl H(\curl, \Omega)^\perp &= H_0\left(\curl^0,\Omega\right),\label{eq:orth_curls}
        \\
        \nabla H^{1}_{0}(\Omega)^\perp &= H\left(\div^0,\Omega\right),\label{eq:orth_traceless_gradients}
        \\
        \curl H_0(\curl, \Omega)^\perp &= H\left(\curl^0,\Omega\right).\label{eq:orth_traceless_curls}
    \end{align}
    Furthermore, all these linear subspaces are $L^2$ closed, so that the relations are still correct by taking the orthogonal to the other side.
\end{proposition}
\begin{proof}
    \cref{eq:orth_gradients,eq:orth_curls,eq:orth_traceless_gradients,eq:orth_traceless_curls} are given by \cite[Eqs. (2.15) and (2.16)]{arnoldFiniteElementExterior2006} in the cases $k=1$ and $k=2$. The fact that all the subspaces are $L^2$ closed is given by the continuity of $\div:H(\div, \Omega) \rightarrow L^2(\Omega)$ and $\curl:H(\curl, \Omega) \rightarrow L^2(\Omega)^3$ on the right-hand sides of \cref{eq:orth_gradients,eq:orth_curls,eq:orth_traceless_gradients,eq:orth_traceless_curls}, and by \cite[Th. 2.3.]{arnoldFiniteElementExterior2006} for the left-hand side.
\end{proof}
\printbibliography
\end{document}